\documentclass[letterpaper]{article}
\usepackage{amssymb}
\usepackage{amsmath}
\usepackage{amsfonts}

\input{tcilatex}
\begin{document}

\begin{center}
$\mathstrut $\textbf{On Harmonicity Of Holomorphic Maps}

\textbf{Between Various Types Of}

\textbf{Almost Contact Metric Manifolds \ }\ 

\ \ \ \textbf{\ }

\textit{Sadettin Erdem}

e-mail: s2erdem@anadolu.edu.tr

serdemays@gmail.com

\bigskip
\end{center}

\textbf{Abstract:}{\small \ Harmonicity of holomorphic maps between various
subclasses of almost contact metric manifolds is discussed. Consequently,
some new results are obtained. Also some known results are recovered, some
of them are generalized and some of them are corrected.}

\begin{center}
\textbf{Introduction}
\end{center}

In complex analysis, it is a basic fact that any $\left( \pm \right) $%
-holomorphic map from the complex vector space $%
%TCIMACRO{\U{2102} }%
%BeginExpansion
\mathbb{C}
%EndExpansion
^{m}$ into $%
%TCIMACRO{\U{2102} }%
%BeginExpansion
\mathbb{C}
%EndExpansion
^{n}$ \ is harmonic. Considering these vector spaces as a trivial example of
a flat Kaehler manifolds, in 1964 Eells and Sampson, $\left[ 12\right] ,$
generalized this harmonicty result to the maps among Kaehler manifolds.
Further generalization given by Lichnerowicz, in 1970 which states that any $%
\left( \pm \right) $-holomorphic map from a semi-Kaehler manifold to a
quasi-Kaehler one\ is harmonic.

However, these generalizations are done within the frame work of almost
Hermitian manifolds which are necessarily of even dimensions. The odd
dimensional counterparts, namely almost contact metric manifolds, were
included in the picture for the first time in 1995 by Ianus and Pastore, $%
\left[ 19\right] .$ Since then harmonicity of \ "holomorphic" maps among
manifolds of mixture of even and odd dimensions has been attracting the
attention of researches. Each work done in this line so far involves only
one or a few specific sub-classes of almost contact metric manifolds, (e.g.
Sasakian, cosymplectic, Kenmotsu manifolds ect). Consequently, certain tools
and arguments are developed for each cases in order to prove the results
claimed. However, in this article we developed the tools and the arguments
which can work almost all sub-classes that are delt with. One of the main
source appealed for their definitions and some of their properties of these
sub-classes is Oubina's works; $\left[ 25\right] ,$ $\left[ 26\right] .$ Our
work acomplishes mainly the following:

$\ i)$ \ It gives new harmonicity results of "holomorphic" maps among
manifolds.

$ii)$ \ It generalizes or recovers some results known (see \textbf{Remaks }$%
\left( 4.2\right) ,$ $\left( 4.3\right) ,$ $\left( 4.4\right) ,$ $\left(
4.5\right) ,$ $\left( 4.6\right) ,$ $\left( 4.7\right) ,$ $\left( 4.8\right)
,$ $\left( 4.9\right) ,\left( 4.10\right) $)$.$

$iii)$ \ It provides some corrections to some works of others (see \textbf{%
Remaks} $\left( 4.5\right) ,\left( 4.8\right) )$.

At the end of the work,we provide some tables and diagrams as a quick
reference in order to give the readers an opportunity to have a globle-look
at definitions, some properties needed of manifolds involved and their
inclusion relations ect.

\textbf{1) Preliminaries }

Let\ $\ (M^{2m+1}$, $g)$ \ be a Riemannian manifold of dimension $(2m+1)$.
If there exist a globally defined $(1,1)$-tensor field $\varphi $, a vector
field $\ \xi $\ \ and a $1$-form $\ \eta $\ \ such that for every local
sections $\ X,Y\in \Gamma (TM)$%
\begin{equation*}
\varphi ^{2}\left( X\right) =-X+\eta \left( X\right) \xi \text{\ ,\ \ \ \ \
\ \ \ }\eta \left( \xi \right) =1
\end{equation*}%
and 
\begin{equation*}
g\left( \varphi X,\ \varphi Y\right) =g\left( X,\ Y\right) -\eta \left(
X\right) \eta \left( Y\right)
\end{equation*}%
hold, then $\ M=\left( M^{2m+1};\ \varphi ,\ \xi ,\ \eta ,\ g\right) $ is
called \textit{an almost contact metric manifold. }In this case one has:%
\begin{equation*}
g\left( X,\ \xi \right) =\eta \left( X\right) \ ,\ \ \ \ \ \varphi \left(
\xi \right) =0\ \ \ \ \ \ \ and\ \ \ \ \ \ \ \ \ \eta \circ \varphi =0
\end{equation*}%
and that \ $rank(\varphi )=2m=\left( \dim M\right) -1.$

Unless otherwise stateted, the letters $X,\ Y,\ Z$ \ will be reserved for
local vector fields or local sections and $\nabla =\nabla ^{M}$ \ will
denote the Levi-Civta connection on the indicated manifold $M$ \ throughout.

\begin{center}
$\underline{Some\ Notations,Definitions\ and\ Basic\ Facts}$
\end{center}

$\ i)\ \ $%
\begin{equation*}
D=D_{\varphi }^{M}:=\varphi \left( TM\right) =\ker \left( \eta \right) \ \ \
\ \ \ \ \ \ 
\end{equation*}%
\ \ \ \ \ \ \ \ \ \ \ \ \ \ \ \ \ \ \ \ \ \ \ \ \ \ \ \ \ \ \ \ \ \ \ \ \ \
\ \ \ \ \ \ \ \ \ \ \ \ \ \ \ \ \ \ \ \ \ \ \ \ \ \ \ \ \ \ \ \ \ \ \ \ \ \
\ \ \ \ \ \ \ \ \ \ \ \ \ \ \ \ \ \ \ \ \ \ \ \ \ \ \ \ \ \ \ \ \ \ \ \ \ \
\ \ \ \ \ \ \ \ \ \ \ \ \ \ \ \ \ \ \ \ \ \ \ \ \ \ \ \ \ \ \ \ \ \ \ \ \ \
\ \ \ \ \ \ \ \ \ \ \ \ \ \ \ \ \ \ \ \ \ \ \ \ \ \ \ \ \ \ \ \ \ \ \ \ \ \
\ \ \ \ \ \ \ \ \ \ \ \ \ \ \ \ \ \ \ \ \ \ \ \ \ \ \ \ \ \ \ \ \ \ \ \ \ \
\ \ \ \ \ \ \ \ \ \ \ \ \ \ \ \ \ \ \ \ \ \ \ \ \ \ \ \ \ \ \ \ \ \ \ \ \ \
\ \ \ \ \ \ \ \ \ \ \ \ \ \ \ \ \ which is called \textit{contact
distribution on }$\ M.$ Note that \ $rank(D)=2m$ \ and \ $TM=D\oplus 
\overline{\xi },$ where $\overline{\xi }$\ denotes the line bundle
determined by $\xi .$

$ii)\ \ $%
\begin{equation*}
\Omega \left( X,\ Y\right) =\ \Omega _{M}\left( X,\ Y\right) :=g\left( X,\
\varphi Y\right) \ \ \ \ \ \ 
\end{equation*}%
\ \ defines a 2-form on $M$, that is, \ $\Omega \in $ $\wedge ^{2}\left(
M\right) .$

$iii)$\ \ \ \ \ \ \ \ \ \ \ \ \ \ \ \ \ \ \ \ \ \ \ \ \ \ \ \ \ \ \ \ \ \ \
\ \ \ \ \ \ \ \ \ \ \ \ \ \ \ \ \ \ \ \ \ \ \ \ \ \ \ \ \ \ \ \ \ \ \ \ \ \
\ \ \ \ \ \ \ \ \ \ \ \ \ \ \ \ \ \ \ \ \ \ \ \ \ \ \ \ \ \ \ \ \ \ \ \ \ \
\ \ \ \ \ \ \ \ \ \ \ \ \ \ \ \ \ \ \ \ \ \ \ \ \ \ \ \ \ \ \ \ \ \ \ \ \ \
\ \ \ \ \ \ \ \ \ \ \ \ \ \ \ \ \ \ \ \ \ \ \ \ \ \ \ \ \ \ \ \ \ \ \ \ \ \
\ \ \ \ \ \ \ \ \ \ \ \ \ \ \ \ \ \ \ \ \ \ \ \ \ \ \ \ \ \ \ \ \ \ \ \ \ \
\ \ \ \ \ \ \ \ \ \ \ \ \ \ \ \ \ \ \ \ \ \ \ \ \ \ \ \ \ \ \ \ \ \ \ \ \ \
\ \ \ \ \ \ \ \ $\ \ $%
\begin{equation*}
P\left( X,\ Y\right) =P_{\varphi }^{M}\left( X,\ Y\right) :=\left( \nabla
_{X}\varphi \right) Y+\left( \nabla _{\left( \varphi X\right) }\varphi
\right) \varphi Y\ \ 
\end{equation*}

\bigskip $iv)$%
\begin{equation*}
Q\left( X,\ Y\right) =Q_{\varphi }^{M}\left( X,\ Y\right) :=\left( \nabla
_{X}\varphi \right) Y+\left( \nabla _{Y}\varphi \right) X\ \ 
\end{equation*}

\bigskip $v)$%
\begin{equation*}
S\left( X,\ Y\right) =S_{\varphi }^{M}\left( X,\ Y\right) :=\left( \nabla
_{X}\varphi \right) Y-\left( \nabla _{Y}\varphi \right) X\ \ 
\end{equation*}

It is easy to see that

$a^{\circ })$ \ 
\begin{equation*}
P\left( X,\ \varphi X\right) =S\left( X,\ \varphi X\right) +\eta \left(
X\right) \left( \nabla _{\left( \varphi X\right) }\varphi \right) \xi
\end{equation*}%
so that $\forall \ X\in \Gamma \left( D^{M}\right) $%
\begin{equation*}
P\left( X,\ \varphi X\right) =S\left( X,\ \varphi X\right) .
\end{equation*}

$b^{\circ })$ 
\begin{equation*}
2P\left( X,\ X\right) =Q\left( X,\ X\right) +Q\left( \varphi X,\ \varphi
X\right)
\end{equation*}

$c^{\circ })$ \ $P,\ Q,\ S:TM\times TM\rightarrow TM$ are all tensor fields.

$vi)\ \ $%
\begin{equation*}
2\mathfrak{h:=\tciLaplace }_{\xi }\varphi \mathfrak{\ \ \ }
\end{equation*}%
$\mathfrak{\ \ \ \ \ \ \ \ \ \ \ \ \ \ \ \ \ \ \ \ \ \ \ \ \ \ \ \ \ \ \ \ \
\ \ \ \ \ \ \ \ \ \ \ \ \ \ \ \ \ \ \ \ \ \ \ \ \ \ \ \ \ \ \ \ \ \ \ \ \ \
\ \ \ \ \ \ \ \ \ \ \ \ \ \ \ \ \ \ \ \ \ \ \ \ \ \ \ \ \ \ \ \ \ \ \ \ \ \
\ \ \ \ \ \ \ \ \ \ \ \ \ \ \ \ \ \ \ \ \ \ \ \ \ \ \ \ }$\ \ \ \ \ \ \ \ \
\ \ \ \ \ \ \ \ \ \ \ \ \ \ \ \ \ \ \ \ \ \ \ \ \ \ \ \ \ \ \ \ \ \ \ \ \ \
\ \ \ \ \ \ \ \ \ \ \ \ \ \ \ \ \ \ \ \ \ \ \ \ \ \ \ \ \ \ \ \ \ \ \ \ \ \
\ \ \ \ \ \ \ \ \ \ \ \ \ \ \ \ \ \ \ \ \ \ \ \ \ \ \ \ \ \ \ \ \ \ \ \ \ \
\ \ \ \ \ \ \ \ \ \ \ \ \ \ \ \ \ \ \ \ \ \ \ \ \ \ \ \ \ \ \ \ \ \ \ \ \ \
\ \ \ \ \ \ \ \ \ \ \ \ \ \ where \ \tciLaplace\ \ denotes the Lie
derivative, so that $2\mathfrak{h}\left( X\right) =\left[ \xi \ ,\ \varphi X%
\right] -\varphi \left[ \xi \ ,\ X\right] .$

$vii)\ \ $The set $\left\{ e_{1},\cdots \ e_{m};\varphi e_{1},\cdots ~\
\varphi e_{m},\ \xi \right\} $ \ will denote a local orthonormal frame field
over $M$ \ throughout.

$viii)\ \ W\left( X,Y\right) =^{\varphi }W^{M}\left( X,Y\right) :=\nabla
_{X}Y+\nabla _{\left( \varphi X\right) }\varphi Y$ \ \ and \ \ $W\left(
X,X\right) =W_{X}$

$ix)\ \ $The exterior differentials 
\begin{equation*}
d\eta \in \wedge ^{2}\left( M\right) \ \ \ \ \ \ and\ \ \ \ \ \ \ \ d\Omega
\in \wedge ^{3}\left( M\right)
\end{equation*}%
of $\eta $ \ and \ $\Omega $ \ are given by 
\begin{equation*}
2d\eta \left( X,\ Y\right) =X\eta \left( Y\right) -Y\eta \left( X\right)
-\eta \left[ X,\ Y\right]
\end{equation*}%
and%
\begin{equation*}
d\Omega \left( X,\ Y,\ Z\right) =\left( \nabla _{X}\Omega \right) \left( Y,\
Z\right) \ +\left( \nabla _{Y}\Omega \right) \left( Z,\ X\right) +\left(
\nabla _{Z}\Omega \right) \left( X,\ Y\right) \ \ 
\end{equation*}%
where 
\begin{equation*}
\left( \nabla _{X}\Omega \right) \left( Y,\ Z\right) =\nabla _{X}\Omega
\left( Y,\ Z\right) \ -\Omega \left( \nabla _{X}Y,\ Z\right) -\Omega \left(
Y,\ \nabla _{X}Z\right) .
\end{equation*}

$x)\ \ $The codifferentials 
\begin{equation*}
\delta \eta \in \wedge ^{0}\left( M\right) =C^{\infty }\left( M\right) ,\ \
\ \ \ \delta \varphi \in \Gamma \left( TM\right) \ \ \ \ \ and\ \ \ \ \ \ \
\ \delta \Omega \in \wedge ^{1}\left( M\right)
\end{equation*}%
of $\eta ,\ \varphi $ \ and \ $\Omega $ \ are given by 
\begin{equation*}
\delta \eta =-\underset{i=1}{\overset{m}{\sum }}\left\{ \left( \nabla
_{e_{i}}\eta \right) e_{i}+\left( \nabla _{\varphi e_{i}}\eta \right)
\varphi e_{i}\right\} =\eta \left( \underset{i=1}{\overset{m}{\sum }}%
W_{e_{i}}\right) ,
\end{equation*}%
\begin{equation*}
\delta \varphi =\left( \nabla _{\xi }\varphi \right) \xi +\underset{i=1}{%
\overset{m}{\sum }}\left\{ \left( \nabla _{e_{i}}\varphi \right)
e_{i}+\left( \nabla _{\varphi e_{i}}\varphi \right) \varphi e_{i}\right\}
=-\varphi \left( \nabla _{\xi }\xi \right) +\underset{i=1}{\overset{m}{\sum }%
}P\left( e_{i},\ e_{i}\right)
\end{equation*}%
and%
\begin{eqnarray*}
\delta \Omega \left( X\right) &=&-\left( \nabla _{\xi }\Omega \right) \left(
\xi ,\ X\right) \ -\underset{i=1}{\overset{m}{\sum }}\left\{ \left( \nabla
_{e_{i}}\Omega \right) \left( e_{i},\ X\right) +\left( \nabla _{\varphi
e_{i}}\Omega \right) \left( \varphi e_{i},\ X\right) \right\} \\
&=&g\left( \delta \varphi ,\ X\right) \ 
\end{eqnarray*}%
\ \ \ \ \ \ \ \ We\ say that \textit{the contact distribution \ }$D$\textit{%
\ \ is minimal }if \ \ $\delta \eta =0.$

$xi)$%
\begin{equation*}
\aleph ^{\left( 1\right) }\left( X,\ Y\right) :=\aleph _{\varphi }\left( X,\
Y\right) -2d\eta \left( X,\ Y\right) \xi ,
\end{equation*}%
\ \ \ \ \ \ \ \ \ \ \ \ \ \ \ \ \ \ \ \ \ \ \ \ \ \ \ \ \ \ \ \ \ \ \ \ \ \
\ \ \ \ \ \ \ \ \ \ \ \ \ \ \ \ \ \ \ \ \ \ \ \ \ \ \ \ \ \ \ \ \ \ \ \ \ \
\ \ \ \ \ \ \ \ \ \ \ \ \ \ \ \ \ \ \ \ \ \ \ \ \ \ \ \ \ \ \ \ \ \ \ \ \ \
\ \ \ \ \ \ \ \ \ \ \ \ \ \ \ \ \ \ \ \ \ \ \ \ \ \ \ \ \ \ \ \ \ \ \ \ \ \
\ \ \ \ \ \ \ \ \ \ \ \ \ \ \ \ \ \ \ \ \ \ \ \ \ \ \ \ \ \ \ \ \ \ \ \ \ \
\ \ \ \ \ \ \ \ \ \ \ \ \ \ \ \ \ \ \ \ \ \ \ \ \ \ \ \ \ \ \ \ \ \ \ \ \ \
\ \ \ \ \ \ \ \ \ \ \ \ \ \ \ \ \ \ \ \ \ \ \ \ \ \ \ \ \ \ \ \ where \ $%
\aleph _{\varphi }$\ \ is the Nijenhuis torsion of \ $\varphi .$ \ Here $%
\aleph ^{\left( 1\right) }$\ \ is called \textit{the normality condition. }%
An almost contact metric manifold is called \textit{normal} if $\ \aleph
^{\left( 1\right) }$ \ vanishes.

One easily gets

\textbf{Lemma }$\mathbf{(1.1)}$\textbf{:}\textit{\ }

\bigskip $1^{\circ })$ \ \textit{For every}{\normalsize \ }\ $X\in \Gamma
\left( TM\right) $%
\begin{equation*}
W_{X}=W\left( X,\ X\right) =S\left( \varphi X,\ X\right) +\varphi \left[
\varphi X\ ,\ X\right] +\nabla _{X}\left( \eta \left( X\right) \xi \right)
\end{equation*}

$\ 2^{\circ })\ \ $\textit{For every}{\normalsize \ }\ $X\in \Gamma \left(
D\right) $

$i)$ $\ $%
\begin{equation*}
W_{X}=S\left( \varphi X,\ X\right) +\varphi \left[ \varphi X\ ,\ X\right]
\end{equation*}

$ii)\ \ \ $%
\begin{equation*}
P\left( X,\ X\right) =\varphi S\left( X,\ \varphi X\right) +\left( \eta 
\left[ X\ ,\ \varphi X\right] \right) \xi
\end{equation*}

$iii)\ \ $%
\begin{equation*}
S\left( \varphi X,\ X\right) =\varphi P\left( X,\ X\right) +\eta \left(
W_{X}\right) \xi
\end{equation*}

$iv)\ \ $%
\begin{equation*}
\left[ X\ ,\ \varphi X\right] =\varphi \ W_{X}+P\left( X,\ X\right)
\end{equation*}

$3^{\circ })\ $ \ \ \ $\forall \ X,\ Y\in \Gamma \left( TM\right) $ \ and \ $%
k,\ r\in C^{\infty }\left( M\right) $

\begin{equation*}
\left( \nabla _{\left( kX\right) }\varphi \right) rY=kr\left( \nabla
_{X}\varphi \right) Y
\end{equation*}%
and therefore 
\begin{equation*}
P\left( kX,\ rY\right) =krP\left( X,\ Y\right) .
\end{equation*}

\textbf{Remark }$\mathbf{(1.1)}$\textbf{\ : }\textit{\ }Note that if \textit{%
\ \ }$d\Omega =0\ \ $then $\ S\left( \varphi X,\ X\right) =0,\ \ \ \forall \
X\in \Gamma \left( D\right) $ \ (see $\ \left[ 13\right] ,$ Proposition
(2.6)). \ So, from Lemma $(1.1)/\left( (2^{\circ })/\left( i\right) \right) $
\ we get\textit{\ }\ $W_{X}=\varphi \left[ X\ ,\ \varphi X\right] ,$ \ hence 
$\ \eta \left( W_{X}\right) =0,\ \ \ \forall \ X\in \Gamma \left( D\right) $
\ and therefore $\ \ \delta \eta =0.$ \ That is, the contact distribution \ $%
D$ \ is minimal when \ $d\Omega =0$.

We shall now be giving definitions of some subclasses of almost contact
metric manifolds \ $M=\left( M^{2m+1};\ \varphi ,\ \xi ,\ \eta ,\ g\right) $
with some of their properties for latter use:

\textbf{Definition }$\mathbf{(1.1)}$\textbf{:} \ An almost contact metric
manifold $\ M=\left( M^{2m+1};\ \varphi ,\ \xi ,\ \eta ,\ g\right) $ \ is
called

$1^{\circ })$ \ $\left( \left[ 9\right] \right) ,\ \alpha $-\textit{contact
metric (or almost }$\alpha $-\textit{Sasakian})\textit{\ manifold \ }if%
\begin{equation*}
d\eta =\alpha \Omega ,
\end{equation*}%
where \ $\alpha \in 
%TCIMACRO{\U{211d} }%
%BeginExpansion
\mathbb{R}
%EndExpansion
-\left\{ 0\right\} =%
%TCIMACRO{\U{211d} }%
%BeginExpansion
\mathbb{R}
%EndExpansion
_{0}.\ 1$-{\normalsize contact metric (or almost }$1$-Sasakian)\textit{\ }%
manifold \textit{\ }is simply called \textit{contact metric}{\normalsize \
(or}\textit{\ almost Sasakian) manifold. }

Here on, $\alpha \in 
%TCIMACRO{\U{211d} }%
%BeginExpansion
\mathbb{R}
%EndExpansion
_{0}$ \ and we omit the letter \ $\alpha $ \ when \ $\alpha =1.$

$2^{\circ })$ \ \textit{quasi contact metric manifold \ }if%
\begin{equation*}
P\left( X,\ Y\right) :=2g\left( X,\ Y\right) \xi -\eta \left( Y\right)
\left\{ X+\eta \left( X\right) \xi +\mathfrak{h}\left( X\right) \right\}
\end{equation*}

$\left( 2^{\circ }.1\right) :\ \ \left( \left[ 3\right] ,\text{ Lemma }%
7.3\right) $; every contact metric manifold is quasi contact metric.

$3^{\circ })$ \ \textit{nearly }$\alpha $-\textit{contact metric manifold \ }%
if%
\begin{equation*}
2\alpha \Omega \left( X,\ Y\right) :=d\eta \left( X,\ Y\right) +d\eta \left(
\varphi X,\ \varphi Y\right)
\end{equation*}

$\left( 3^{\circ }.1\right) :\ \ $($\left[ 23\right] ,$ Lemma $2$); every
quasi contact metric manifold is nearly contact metric.

$4^{\circ })\ \ \left( \left[ 9\right] \right) ,\ $ $\alpha $-\textit{%
Sasakian manifold }if it is almost $\alpha $-Sasakian (or $\alpha $-contact
metric) and normal, or equivalently%
\begin{equation*}
\left( \nabla _{X}\varphi \right) Y:=\alpha \left\{ g\left( X,\ Y\right) \xi
-\eta \left( Y\right) \left( X\right) \right\}
\end{equation*}

$5^{\circ })\ \ $\textit{strongly pseudoconvex }$CR$\textit{-manifold} if it
is a contact metric manifold satisfying:%
\begin{equation*}
\left( \nabla _{X}\varphi \right) Y:=g\left( X+\mathfrak{h}X,\ Y\right) \xi
-\eta \left( Y\right) \left( X+\mathfrak{h}X\right) .
\end{equation*}

$\left( 5^{\circ }.1\right) :\ \ ${\normalsize Every Sasakian manifold is a
strongly pseudoconvex \ }${\normalsize CR}${\normalsize -manifold}\textit{. }%
Conversally,\textit{\ }{\normalsize every strongly pseudoconvex \ }$%
{\normalsize CR}${\normalsize -manifold with}\textit{\ \ \ }$\mathfrak{h}=0$%
\textit{\ \ }{\normalsize is \ Sasakian.}

{\normalsize \ }$6^{\circ })\ \ $\textit{nearly }$\alpha $\textit{-Sasakian
manifold}{\normalsize \ \ if }%
\begin{equation*}
Q\left( X,Y\right) :=\alpha \left\{ 2g\left( X,\ Y\right) \xi -\eta \left(
Y\right) X-\eta \left( X\right) Y\right\}
\end{equation*}

$\left( 6^{\circ }.1\right) :\ \ $Every $\alpha $-Sasakian manifold is
nearly $\alpha $-Sasakian.

$\left( 6^{\circ }.2\right) :\ (e.g.\ \left[ 23\right] ,$ Lemma $\ 2);$
every nearly $\alpha $-Sasakian manifold is a nearly $\alpha $-contact
metric manifold.

$7^{\circ })\ \ (c.f.\ \left[ 4\right] )$, almost \textit{quasi-Sasakian
manifold of rank }\ $2r+1,\ \ \left( 0\leq r\leq m\right) $ \ if%
\begin{equation*}
d\Omega =0\ \ \ \ \ and\ \ \ \ \ \ \eta \wedge \left( d\eta \right) ^{r}\neq
0\ \ \ \ \ with\ \ \ \ \ \ \ \left( d\eta \right) ^{r+1}\ =0.\ \ \ \ 
\end{equation*}

However, in this work, we restrict ourselves to the cases where $1\leq r\leq
\left( m-1\right) ,$ as the other cases covered by other classes we shall be
dealing with.

$\bullet $ \ A \textbf{normal} almost$\ $quasi-Sasakian manifold will be
called \textit{quasi-Sasakian}

$\left( 7^{\circ }.1\right) :\ \ $Every almost $\alpha $-Sasakian (i.e.
contact metric) manifold \ $M^{2m+1}$\ \ is an almost quasi-Sasakian one of
full rank ( i.e. of rank $2m+1$ )

$\left( 7^{\circ }.2\right) :\ \ $On a quasi-Sasakian manifold one has

$\ \ \ \ i)\ \ \delta \eta =0.$

$\ \ \ ii)\ \ \delta \Omega \left( \varphi X\right) =0$\ \ \ \ and so \ $%
\delta \Omega \left( X\right) =\eta \left( X\right) \delta \Omega \left( \xi
\right) .$

$8^{\circ })\ \ $\textit{\ almost }$\beta $-\textit{Kenmotsu manifold} if \ 
\begin{equation*}
\ d\Omega =2\beta \eta \wedge \Omega \ ,\ \ \ \ \ d\eta =0\ \ \ \ \text{and\ 
}\ \ \ \ \ \ d\beta \wedge \eta =0\ \ \ \ \ \ \ \ \ \ 
\end{equation*}%
where \ $\beta \in C^{\infty }\left( M\right) $ \ \ with \ $\beta \left(
p\right) \neq 0,\ \ \forall \ p\in M.$

$\bullet $ \ A \textbf{normal} almost$\ \beta $-Kenmotsu manifold will be
called $\beta $\textit{-Kenmotsu.}

$\left( 8^{\circ }.1\right) :(c.f.\ \left[ 27\right] ,\ $Theorem $\left(
3.3\right) ):$ \ An almost\ $\beta $-Kenmotsu manifold is $\beta $-Kenmotsu
if and only if \ 
\begin{equation*}
\left( \nabla _{X}\varphi \right) Y:=\beta \left\{ g\left( \varphi X,\
Y\right) \xi -\eta \left( Y\right) \left( \varphi X\right) \right\} .
\end{equation*}

$\left( 8^{\circ }.2\right) :$ \ An almost $\beta $-Kenmotsu manifold is
said to be \ $CR$\textit{-integrable} if it satisfies 
\begin{equation*}
\left( \nabla _{X}\varphi \right) Y:=g\left( \mathfrak{h}\left( X\right)
,Y\right) \xi -\eta \left( Y\right) \mathfrak{h}\left( X\right) +\beta
\left\{ g\left( \varphi X,Y\right) \xi -\eta \left( Y\right) \left( \varphi
X\right) \right\} .
\end{equation*}

$\left( 8^{\circ }.3\right) :$ \ An almost $\beta $-Kenmotsu manifold is $%
\beta $-Kenmotsu if and only if it is$\ \ CR$-integrable$\ $and $\ \mathfrak{%
h}=0.$

$\left( 8^{\circ }.4\right) :\ (c.f.\ \left[ 27\right] ,\ $Lemma $\left(
3.2\right) ):$ \ An \textbf{almost}$\ \beta $-Kenmotsu manifold satisfies 
\begin{equation*}
P\left( X,\ Y\ \right) =\beta \left\{ 2g\left( \varphi X,\ Y\right) \xi
-\eta \left( Y\right) \left( \varphi X\right) \right\} -\eta \left( Y\right) 
\mathfrak{h}\left( X\right) .
\end{equation*}

$9^{\circ })\ \ (\left[ 26\right] ,\ \left[ 4\right] );$ \ \textit{%
trans-Sasakian manifold of type} \ $\left( \wp ,\theta \right) $ \ if it is 
\textbf{normal} and 
\begin{equation*}
\left( \nabla _{X}\varphi \right) Y=\wp \left\{ g\left( X,\ Y\right) \xi
-\eta \left( Y\right) X\right\} +\theta \left\{ g\left( \varphi X,\ Y\right)
\xi -\eta \left( Y\right) \varphi X\right\}
\end{equation*}%
holds, where 
\begin{equation*}
2m\wp :=\delta \Omega \left( \xi \right) \ =\eta \left( \underset{i=1}{%
\overset{m}{\sum }}P\left( e_{i},\ e_{i}\right) \right) \ \ \ \ \ and\ \ \ \
2m\theta :=\func{div}\xi =-\delta \eta =-\eta \left( \underset{i=1}{\overset{%
m}{\sum }}W_{e_{i}}\right) .\ \ 
\end{equation*}

$\left( 9^{\circ }.1\right) :\ (\left[ 24\right] );$ \ on a trans-Sasakian
manifold $M$ \ with $\ dim\ \left( M\right) \geq 5,$ precisely one of the
following holds:

\ \ \ $i)$ \ $\ \wp =\alpha $\ $\in 
%TCIMACRO{\U{211d} }%
%BeginExpansion
\mathbb{R}
%EndExpansion
_{0}$ \ and \ $\theta =0;$ \ in which case $M$ is $\alpha $-Sasakian

$\ \ ii)\ \ \ \wp =0$ \ and \ $\theta \neq 0;$ \ in which case $M$ is $%
\theta $-Kenmotsu .

$\ iii)\ \ \ \wp =0$ \ and \ $\theta =0;$ \ in which case $M$ is
cosymplectic.\ 

However, in this work throughout, trans-Sasakian manifolds $\ M$ \ of type \ 
$\left( \wp ,\theta \right) $ \ would be the ones with $\ \wp \left(
p\right) \neq 0\neq \theta \left( p\right) ,$ \ $\forall \ p\in M$ \ and
therefore, by the virtue of the above result, \ $dim\ \left( M\right) =3.$
The other cases where \ $dim\ \left( M\right) \geq 5,$ would be delt in the
relevant parts.

$10^{\circ })\ \ (\left[ 16\right] );$ \textit{nearly trans-Sasakian
manifold }$\ $\textit{of type} \ $\left( \wp ,\theta \right) $ 
\begin{eqnarray*}
Q\left( X,Y\right) &:&=\wp \left\{ 2g\left( X,\ Y\right) \xi -\eta \left(
X\right) Y\ -\eta \left( Y\right) X\right\} \\
&&-\theta \left\{ \eta \left( X\right) \varphi Y+\eta \left( Y\right)
\varphi X\right\} .
\end{eqnarray*}

$\left( 10^{\circ }.1\right) :\ \ ${\normalsize Every }trans-Sasakian
manifold is a {\normalsize nearly trans-Sasakian manifold of the same type.}

$11^{\circ })\ \ (c.f.$ $\left[ 26\right] ,\left[ 8\right] );$ \textit{almost%
}\ \textit{semi-}$\alpha $\textit{-Sasakian manifold} if 
\begin{equation*}
2m\alpha \eta :=\delta \Omega .
\end{equation*}

$\bullet $ \ A \textbf{normal} almost$\ $\textit{semi-}$\alpha $\textit{-}%
Sasakian manifold will be called \textit{semi-}$\alpha $-\textit{Sasakian}

$\left( 11^{\circ }.1\right) :\ \ (c.f.$ $\left[ 26\right] ,$ Theorem$\
\left( 3.5\right) );$ \ every {\normalsize nearly-}$\alpha $-{\normalsize \
Sasakian manifold }is almost s{\normalsize emi-}$\alpha ${\normalsize %
-Sasakian.}

$\left( 11^{\circ }.2\right) :\ \ $On an almost semi-$\alpha $-Sasakian
manifold (and therefore on its subclasses)$\ \ \ \delta \eta =0.$

$12^{\circ })\ \ (\left[ 26\right] ,$ Theorem\ $2.9\ );\ \ $\textit{quasi-}$%
\mathcal{K}$\textit{-Sasakian manifold} if 
\begin{equation*}
P\left( X,\ Y\right) =2g\left( X,\ Y\right) \xi -2\eta \left( Y\right)
X+\eta \left( Y\right) \nabla _{\varphi X}\xi
\end{equation*}

$\left( 12^{\circ }.1\right) :\ $\ $\left( \left[ 26\right] ,\text{ Theorem}%
\ \left( 3.3\right) ,\ \text{and\ \ Proposition }\left( 3.12\right) \right)
;\ $\ {\normalsize every }quasi\textit{-}$\mathcal{K}$\textit{-}Sasakian
manifold is almost s{\normalsize emi-Sasakian and nearly contact.}$\ $

$\left( 12^{\circ }.2\right) :\ \ (\left[ 26\right] ,$ Corollary\ $3.13);$ \ 
{\normalsize every contact metric manifold is }quasi-$\mathcal{K}$-Sasakian.

$13^{\circ })\ \ $\textit{almost cosymplectic manifold} if 
\begin{equation*}
d\eta =0\ \ \ \ and\ \ \ \ \ d\Omega =0
\end{equation*}

$\bullet $ \ A \textbf{normal} almost$\ $cosymplectic manifold will be
called \textit{cosymplectic}

$\left( 13^{\circ }.1\right) :$ \ A manifold is cosymplectic if and only if

\begin{equation*}
\left( \nabla _{X}\varphi \right) Y=0.
\end{equation*}

$\left( 13^{\circ }.2\right) :\ \ ${\normalsize Every cosymplectic manifold
is a }quasi\textit{-}Sasakian manifold of rank one.

$14^{\circ })\ \ $\textit{nearly-}$\mathcal{K}$\textit{-cosymplectic manifold%
} if%
\begin{equation*}
\left( \nabla _{X}\varphi \right) X=0\ \ \ \ \ \ and\ \ \ \ \ \ \ \ \ \
\left( \nabla _{X}\varphi \right) \xi =0\ \ \ \ .
\end{equation*}

$15^{\circ })$ \ \ \textit{nearly-cosymplectic manifold} if%
\begin{equation*}
\left( \nabla _{X}\varphi \right) X=0.
\end{equation*}

$\ \left( 15^{\circ }.1\right) :\ $Clearly every nearly-$\mathcal{K}$\textit{%
-}cosymplectic manifold is nearly-cosymplectic.

$16^{\circ })\ \ \ \left( \left[ 26\right] ,\text{Theorem}\ 2.9\right) $ \ 
\textit{quasi-}$\mathcal{K}$\textit{-cosymplectic manifold} if%
\begin{equation*}
P\left( X,\ Y\right) :=\eta \left( Y\right) \nabla _{\varphi X}\xi
\end{equation*}

$\left( 16^{\circ }.1\right) :\ \ \left( \left( \left[ 26\right] ,\text{%
Theorem}\ \left( 3.3\right) \right) ,\ \left[ 8\right] \right) ;$ every
almost cosymplectic and every nearly-$\mathcal{K}$-cosymplectic manifolds
are also quasi-$\mathcal{K}$-cosymplectic

$17^{\circ })\ \ $\textit{almost semi-cosymplectic manifold} if%
\begin{equation*}
\delta \eta =0\ ,\ \ \ \ \ \ \ \ \ \ \ \ \ \ \ \ \ \ \delta \Omega =0
\end{equation*}

$\bullet $ \ A \textbf{normal} almost$\ $semi-cosymplectic manifold will be
called \textit{semi-cosymplectic}

$\ \left( 17^{\circ }.1\right) :\ \left( \left( \left[ 26\right] ,\text{%
Theorem}\ \left( 3.3\right) ,\ \text{and\ Theorem}\ \left( 3.4\right)
\right) ,\ \left[ 8\right] \right) ;$ every nearly-cosymplectic and every
quasi-$\mathcal{K}$-cosymplectic manifolds are also almost semi-cosymplectic

$18^{\circ })\ \ $\textit{quasi-symplectic manifold} if%
\begin{equation*}
S\left( \varphi X,\ X\right) =0,\ \ \ \ \forall \ X\in \Gamma \left( D\right)
\end{equation*}

$19^{\circ })\ \ $\textit{semi-symplectic manifold} if%
\begin{equation*}
\underset{i=1}{\overset{m}{\sum }}S\left( \varphi e_{i},\ e_{i}\right) =0
\end{equation*}

We also recall here some subclasses of almost Hermitian manifolds \ \ \ $%
H=\left( H^{2m};\ J,\ \ G\right) .$ \ Note here that the $(1,1)$-tensor
field $\ \varphi =J$ \ is of full rank \ $2m$ \ and therefore \ $%
D=D_{H}^{J}:=J\left( TH\right) =TH.$ \ Set 
\begin{equation*}
\Phi \left( X,\ Y\right) =\ \Phi _{H}\left( X,\ Y\right) :=G\left( X,\
JY\right) \ \ \ \ \ \ 
\end{equation*}

\textbf{Definition }$\mathbf{(1.2)}$: \ An almost Hermitian manifold $%
H=(H^{2m};\ J,\ G)$ \ is called

$1^{\circ })\ \ \ $\textit{almost Kaehler if }

\begin{equation*}
d\Phi =0\text{ \ \ \ \ }
\end{equation*}

$2^{\circ })\ \ \ $\textit{Kaehler if }

\begin{equation*}
d\Phi =0\text{ \ \ \ \ \ and \ \ \ \ }\aleph _{J}=0
\end{equation*}%
or equivalently

$\ \ \ $%
\begin{equation*}
\nabla J=0.
\end{equation*}

$3^{\circ })$ \ \textit{\ nearly-Kaehler} if 
\begin{equation*}
\left( \nabla _{X}J\right) X=0,\ \ 
\end{equation*}

$4^{\circ })$ \ \textit{\ quasi-Kaehler (or }$\left( 1,2\right) $-\textit{%
symplectic}) if%
\begin{equation*}
S\left( JX,\ X\right) =0.
\end{equation*}

$\left( 4^{\circ }.1\right) :$ \ \ Unlike the almost contact metric case, on
an almost Hermitian manifold one has%
\begin{equation*}
S\left( JX,\ X\right) =P\left( JX,\ X\right) =JP\left( X,\ X\right) .
\end{equation*}%
$\ \ \ \ \ \ \left( 4^{\circ }.2\right) :$ \ \ The following are equivalent:

$\ \ \bullet $ \ \ $H$ \ is quasi-Kaehler i.e. \ $S\left( JX,\ X\right)
=P\left( JX,\ X\right) =0$

$\ \ \bullet $ \ \ \ $P\left( X,\ X\right) =0,\ \ \forall \ X\in \Gamma
\left( TH\right) $

$\ \ \bullet $ \ \ \ $P\left( X,\ Y\right) =0,\ \ \forall \ X,\ Y\in \Gamma
\left( TH\right) $

$5^{\circ })$ \ \textit{\ semi-Kaehler }{\normalsize if }%
\begin{equation*}
\overset{m}{\underset{i=1}{\sum }}S\left( Je_{i},\ e_{i}\right) =0,
\end{equation*}%
where \ \ $\left\{ e_{1},\cdots \ e_{m};Je_{1},\cdots \ Je_{m}\right\} $ \
is a local orthonormal frame field over $H.$

\textbf{2) \ Some Auxiliary Results}

\textbf{Definition (2.1): }\textit{An almost contact metric manifold \ }$M$%
\textit{\ \ is said to satisfy geodesic condition \ }$\left( GC\right) $%
\textit{\ \ if \ }$\nabla _{\xi }\xi =0.$

\textbf{Lemma }$\mathbf{(2.1)}$\textbf{: }$\left( c.f.\ \left( \left[ 13%
\right] ,\text{Proposition }\left( 2.6\right) \right) ,\left( \left[ 6\right]
,\text{Lemma }\left( 3.5\right) \right) \right) $: \textit{\ Let} \ $M$%
\textit{\ \ be an almost contact metric manifold. \ Then \ }$\forall \ X\in
\Gamma \left( D\right) $ \ and $\ \forall \ Z\in \Gamma \left( TM\right) $ \ 
\textit{\ \ }%
\begin{equation*}
g\left( S\left( \varphi X,\ X\right) ,Z\right) =-d\Omega \left( \varphi X,\
X,\ Z\right)
\end{equation*}%
and therefore%
\begin{equation*}
\overset{m}{\underset{i=1}{\sum }}g\left( S\left( \varphi e_{i},\
e_{i}\right) ,\ Z\right) =-\overset{m}{\underset{i=1}{\sum }}d\Omega \left(
\varphi e_{i},e_{i},Z\right) .
\end{equation*}

\textbf{Lemma}$\mathbf{(2.2)}$\textbf{: }\ $\left( e.g.\text{ \ }\left[ 13%
\right] ,\text{ \ Lemma }(2.1)\right) ;$ \ \textit{For an almost contact
metric manifold \ }$M$\textit{\ \ }

$1^{\circ })$ \ $\left[ X,\ \xi \right] \in \Gamma \left( D\right) $\textit{%
,\ }$\forall $\textit{\ }$X\in \Gamma \left( D\right) $\textit{\ \ if and
only if \ }$M$\textit{\ satisfies} $\left( GC\right) .$

$2^{\circ })$ \ \textit{If \ }$M$\textit{\ \ is normal then it satisfies }$%
\left( GC\right) .$

\textbf{Lemma }$\mathbf{(2.3)}$\textbf{:}

$1^{\circ }\mathbf{\ )}:$ \ \textit{For the subclasses of almost contact
metric manifolds listed below the following identity }%
\begin{equation*}
S\left( \varphi X,\ X\right) =0,\ \ \ \ \ \ \ \forall \ X\in \Gamma \left(
D\right)
\end{equation*}%
\textit{holds and therefore they are quasi-symplectic:}

$\ \ i)$ \ \ \textit{Quasi contact metric manifolds }

$ii)$ \ \ \textit{Quasi }$\mathcal{K}$\textit{-Sasakian manifolds}

$iii)$ \ \ \textit{Almost }$\alpha $-Sasakian\textit{\ ( }$\alpha $-\textit{%
contact metric})\textit{\ manifolds and thus its subclasses:}

$\ \ \ \bullet $ \ \ \textit{Strongly pseudoconvex }$CR$\textit{-manifolds}

$\ \ \ \bullet $ \ \ $\alpha $-\textit{\ Sasakian manifolds}

$iv)$ \ \ \textit{Quasi-Sasakian manifolds }

$v)$ \ \ Quasi-$\mathcal{K}$-\textit{cosymplectic manifolds and thus its
subclass}

$\ \ \bullet $ \ \ \ \ Almost \textit{cosymplectic manifolds.}

$\ \ \bullet $ \ \ \ \ \textit{Cosymplectic manifolds.}

$vi)$ \ \ \textit{Nearly}-\textit{cosymplectic manifolds and thus its
subclass: Nearly}-$\mathcal{K}$-\textit{cosymplectic manifolds.}

$2^{\circ }\mathbf{:}$ \ 

$i)$ \ \textit{For an almost semi-}$\alpha $\textit{-Sasakian and almost
semi-cosymplectic manifolds the following identity }%
\begin{equation*}
\underset{i=1}{\overset{m}{\sum }}S\left( \varphi e_{i},\ e_{i}\right)
=-\nabla _{\xi }\xi
\end{equation*}%
\textit{holds.}

$ii)$ \ \textit{For semi-}$\alpha $\textit{-Sasakian and semi-cosymplectic
manifolds} \textit{the following identity }%
\begin{equation*}
\underset{i=1}{\overset{m}{\sum }}S\left( \varphi e_{i},\ e_{i}\right) =0
\end{equation*}%
\textit{holds and therefore they are semi-symplectic.}

\textbf{Proof:}

$1^{\circ }\mathbf{/}\left( (i),\ (ii)\text{ and}\ \ (v)\right) \ \ $ follow
directly from their definitions. \ For $1^{\circ }\ \mathbf{/}\left( (iii)\ 
\text{and}\ (iv)\right) ,$ \ note that \ $d\Omega =0.$ So, by Lemma $\left(
2.1\right) ,$ the result follows.

For $1^{\circ }$\ $\mathbf{/}(vi)$ \ note that by definition 
\begin{equation*}
\left( \nabla _{X}\varphi \right) X=0,\ \ \ \ \ \forall \ X\in \Gamma \left(
TM\right)
\end{equation*}%
and so $\ $

\begin{equation*}
P\left( X,\ X\right) =0.\ \ \ \ \ 
\end{equation*}%
But then this, together with Lemma $(1.1)/\left( \left( 2^{\circ }\right)
/\left( iii\right) \right) $, \ gives that 
\begin{equation}
S\left( \varphi X,\ X\right) =\eta \left( W_{X}\right) \xi \ \ \ \ \ \
\forall \ X\in \Gamma \left( D\right) .  \tag{$2.1$}
\end{equation}%
On the other hand, on a nearly-cosymplectic manifold, the vector field \ $%
\xi $ \ is Killing, $\left( \left[ 2\right] ,\text{Proposition }6.1\right) .$
So we have 
\begin{equation}
\left( \mathcal{\tciLaplace }_{\xi }g\right) \left( X,\ X\right) =2g\left(
\nabla _{X}\xi ,\ X\right) =0,\text{ \ \ \ \ }\forall \ X\in \Gamma \left(
TM\right) .  \tag{$2.2$}
\end{equation}%
Now, \ $\forall X\in \Gamma \left( D\right) $ \ we have%
\begin{eqnarray*}
\eta \left( W_{X}\right) &=&g\left( \nabla _{X}X,\ \xi \right) +g\left(
\nabla _{\varphi X}\left( \varphi X\right) ,\ \xi \right) \\
&=&-g\left( \nabla _{X}\xi ,\ X\right) -g\left( \nabla _{\varphi X}\xi ,\
\varphi X\right)
\end{eqnarray*}%
This gives, using \ $\left( 2,2\right) $ \ that, \ \ $\eta \left(
W_{X}\right) =0.$ \ So, from \ $\left( 2.1\right) ,$ the result follows.

$2^{\circ }\mathbf{)/}\left( i\right) :$ \ Note that $\ \forall \ X\in
\Gamma \left( TM\right) $ 
\begin{eqnarray}
\delta \Omega \left( X\right) &=&g\left( \delta \varphi ,\ X\right)  \notag
\\
&=&g\left( -\varphi \left( \nabla _{\xi }\xi \right) +\underset{i=1}{\overset%
{m}{\sum }}P\left( e_{i},\ e_{i}\right) ,\ X\right) .  \TCItag{$2.3$}
\end{eqnarray}

If $\ M$\ \ is an almost semi-$\alpha $-\textbf{Sasakian} manifold then,
from the definition of \ $M,$ \ we have \ $\ \ $%
\begin{equation*}
2m\alpha \eta \left( X\right) =\delta \Omega \left( X\right) =g\left(
-\varphi \left( \nabla _{\xi }\xi \right) +\underset{i=1}{\overset{m}{\sum }}%
P\left( e_{i},\ e_{i}\right) ,\text{\ }X\right) =0,\text{ \ }\forall \ X\in
\Gamma \left( D\right)
\end{equation*}%
so that 
\begin{equation*}
g\left( \underset{i=1}{\overset{m}{\sum }}P\left( e_{i},\ e_{i}\right) ,\ 
\text{\ }X\right) =g\left( \varphi \left( \nabla _{\xi }\xi \right) ,\
X\right) ,\text{ \ }\forall \ X\in \Gamma \left( D\right)
\end{equation*}%
and%
\begin{equation*}
2m\alpha \eta \left( \xi \right) =\delta \Omega \left( \xi \right) =g\left( 
\underset{i=1}{\overset{m}{\sum }}P\left( e_{i},\ e_{i}\right) ,\text{\ }\xi
\right) =2m\alpha .
\end{equation*}%
The last two equations give us that 
\begin{equation}
\underset{i=1}{\overset{m}{\sum }}P\left( e_{i},\ e_{i}\right) =\varphi
\nabla _{\xi }\xi +2m\alpha \xi .  \tag{$2.4$}
\end{equation}

If $\ M$\ \ is an almost semi-\textbf{cosymplectic} manifold then, from the
definition of \ $M,$ \ we have \ 
\begin{equation*}
\delta \Omega \left( X\right) =g\left( -\varphi \left( \nabla _{\xi }\xi
\right) +\underset{i=1}{\overset{m}{\sum }}P\left( e_{i},\ e_{i}\right) ,%
\text{\ }X\right) =0,\ \ \ \ \forall ~X\in \Gamma \left( TM\right) \text{ }
\end{equation*}%
so that 
\begin{equation}
\underset{i=1}{\overset{m}{\sum }}P\left( e_{i},\ e_{i}\right) =\varphi
\left( \nabla _{\xi }\xi \right) .  \tag{$2.5$}
\end{equation}%
But then, by the virtue of\ Lemma $\left( 1.1\right) /\left( \left( 2^{\circ
}\right) /\left( iii\right) \right) $ \ and using \ $\left( 2.4\right) ,$\ $%
\left( 2.5\right) $ \ we get, on both almost semi-$\alpha $-Sasakian and
almost semi-cosymplectic manifolds, 
\begin{equation*}
\underset{i=1}{\overset{m}{\sum }}S\left( \varphi e_{i},\ e_{i}\right) =%
\underset{i=1}{\overset{m}{\sum }}\varphi P\left( e_{i},\ e_{i}\right)
+\left( \delta \eta \right) \xi =-\nabla _{\xi }\xi +\left( \delta \eta
\right) \xi .
\end{equation*}%
Hence,%
\begin{equation*}
\underset{i=1}{\overset{m}{\sum }}S\left( \varphi e_{i},\ e_{i}\right)
=-\nabla _{\xi }\xi ,
\end{equation*}%
since \ $\delta \eta =0$ \ on both almost semi $\alpha $-Sasakian and almost
semi-cosymplectic manifolds.

$2^{\circ })\mathbf{/}\left( ii\right) :$ \ On semi-$\alpha $-Sasakian and
semi-cosymplectic\textbf{\ }manifolds, \ $\nabla _{\xi }\xi =0$ \ by the
normality of them, (see Lemma $\left( 2.2\right) $). So the equality 
\begin{equation*}
\underset{i=1}{\overset{m}{\sum }}S\left( \varphi e_{i},\ e_{i}\right) =0
\end{equation*}%
follows from the part \ $\left( 2^{\circ }\right) /\left( i\right) .$

\textbf{Lemma }$\mathbf{(2.4)}$\textbf{:} T\textit{he subclasses of almost
contact metric manifolds listed below satisfy }$\left( GC\right) :$

$1^{\circ })$\textit{\ \ Nearly-}$\alpha $-\textit{contact metric manifolds
and thus its subclasses:}

$\ \bullet $\textit{\ \ Quasi-contact metric manifolds.}

$\ \bullet $\textit{\ \ Quasi-}$\mathcal{K}$\textit{-Sasakian manifolds.}

$2^{\circ })$\textit{\ }$\ \alpha $\textit{-Contact metric manifolds and
thus its subclasses:}

$\ \bullet $\textit{\ \ Strongly pseudoconvex CR-manifolds,}

$\ \bullet $\textit{\ \ }$\alpha $\textit{-Sasakian manifolds.}

$3^{\circ })$\textit{\ \ Nearly-}$\alpha $-\textit{Sasakian manifolds.}

$4^{\circ })$\textit{\ \ Quasi-Sasakian manifolds.}

$5^{\circ })$\textit{\ \ Nearly-trans-Sasakian manifolds of type \ }$\left(
\wp ,\theta \right) $\textit{\ \ and thus its subclass:}

$\ \bullet $\textit{\ \ Trans-Sasakian manifolds of type\ }$\left( \wp
,\theta \right) .$

$6^{\circ })$\textit{\ \ Semi-Sasakian manifolds }

$7^{\circ })$\textit{\ \ Quasi-}$\mathcal{K}$\textit{-cosymplectic manifolds
and thus its subclasses:}

$\ \bullet $\textit{\ \ Nearly-}$\mathcal{K}$\textit{-cosymplectic manifolds.%
}

$\ \bullet $\textit{\ \ Almost cosymplectic manifolds.}

$\ \bullet $\textit{\ \ Cosymplectic manifolds.}

$8^{\circ })$\textit{\ \ Nearly-cosymplectic manifolds.}

$9^{\circ })$\textit{\ \ Semi-cosymplectic manifolds.}

$10^{\circ })$\textit{\ \ almost }$\beta $\textit{-Kenmotsu manifolds and
thus its subclasses:}

$\ \bullet $\textit{\ \ }$CR$\textit{-integrable almost }$\beta $-\textit{%
Kenmotsu manifolds.}

$\ \bullet $\textit{\ \ }$\beta $\textit{-Kenmotsu manifolds.}

\textbf{Proof:}

$1^{\circ })$\textbf{\ \ }Let\textbf{\ }$M$ \ be a nearly contact metric
manifold. Note that \ $\forall \ X\in \Gamma \left( D\right) $ \ we have 
\begin{equation*}
-\eta \left( \left[ X,\ \xi \right] \right) =2d\eta \left( X,\ \xi \right)
=2\left\{ d\eta \left( X,\ \xi \right) +d\eta \left( \varphi X,\ \varphi \xi
\right) \right\}
\end{equation*}%
But then, since \ $M$ \ is a nearly contact metric manifold, this gives%
\begin{equation*}
-\eta \left( \left[ X,\ \xi \right] \right) =4\alpha \Omega \left( X,\ \xi
\right) =4\alpha g\left( X,\ \varphi \xi \right) =0.
\end{equation*}%
So, from Lemma $(2.2)$, the result follows.

$2^{\circ })$\textbf{\ \ }Let\textbf{\ \ }$M$ \ be an $\alpha $-contact
metric manifold. Then the result follows by the same argument used above.

$3^{\circ })$\textbf{\ \ }Let \textbf{\ }$M$ \ be a nearly $\alpha $%
-Sasakian manifold. From its definition, we have%
\begin{equation*}
-2\varphi \left( \nabla _{\xi }\xi \right) =Q\left( \xi ,\xi \right) =\alpha
\left\{ 2g\left( \xi ,\ \xi \right) \xi -2\eta \left( \xi \right) \xi
\right\} =0,
\end{equation*}

which gives the result.

$4^{\circ })$\textbf{\ \ }Let \textbf{\ }$M$ \ be a quasi-Sasakian manifold.
Since \ $M$ is normal, from Lemma $(2.2)$, the result follows.

$5^{\circ })$\textbf{\ \ }Let \textbf{\ }$M$ \ be a nearly-trans-Sasakian
manifold $\ $of type \ $\left( \wp ,\theta \right) .\ \ $Using its
definition, we get

\begin{eqnarray*}
-2\varphi \left( \nabla _{\xi }\xi \right) &=&Q\left( \xi ,\xi \right)
=2\left( \nabla _{\xi }\varphi \right) \xi \\
&=&\wp \left\{ 2g\left( \xi ,\ \xi \right) \xi -2\eta \left( \xi \right) \xi
\right\} -\theta \left\{ \eta \left( \xi \right) \varphi \xi +\eta \left(
\xi \right) \varphi \xi \right\} \\
&=&0.
\end{eqnarray*}

From which the result follows.

$6^{\circ })$\textit{\ \ }Let \textbf{\ }$M$ \ be a semi-Sasakian normal
manifold.\textit{\ \ }Since it is normal, from Lemma $(2.2)$, the result
follows

$7^{\circ })$\textit{\ \ }Let \textbf{\ }$M$ \ be a quasi-$K$-cosymplectic
manifold.\textit{\ }From its definition, we have%
\begin{equation*}
-\varphi \left( \nabla _{\xi }\xi \right) =P\left( \xi ,\ \xi \right) =\eta
\left( \xi \right) \left( \nabla _{\left( \varphi \xi \right) }\xi \right)
=0.
\end{equation*}%
From which the result follows.

$8^{\circ })$\textit{\ \ }Let \textbf{\ }$M$ \ be a nearly-cosymplectic
manifold. Then directly from its definition the result follows.

$9^{\circ })$\textit{\ \ }Let \textbf{\ }$M$ \ be a semi-cosymplectic
manifold. Since it is normal, from Lemma $(2.2)$, the result follows.

$10^{\circ })$\textit{\ \ }Let \textbf{\ }$M$ \ be an almost $\beta $%
-Kenmotsu manifold, then by definition \ $d\eta =0.$ So, one gets%
\begin{equation*}
-\eta \left( \left[ X,\ \xi \right] \right) =2d\eta \left( X,\ \xi \right)
=0.
\end{equation*}%
Then, from Lemma $(2.2)$, the result follows.

\textbf{Definition }$\mathbf{(2.2)}$\textbf{: }\textit{An almost contact
metric manifold \ }$\left( M;\varphi ,\xi ,\eta ,g\right) $\textit{\ \ is
said to be}

$1^{\circ })$\textit{\ \ }$\varphi $-\textit{involutive if }%
\begin{equation*}
\left[ X,\ \varphi X\right] \in \Gamma \left( D\right) ,\text{ \ \ \ }%
\forall \text{\ }X\in \Gamma \left( D\right) .
\end{equation*}

$2^{\circ })$\textit{\ \ semi }$\varphi $-\textit{involutive if }%
\begin{equation*}
\underset{i=1}{\overset{m}{\sum }}\left[ e_{i},\ \varphi e_{i}\right] \in
\Gamma \left( D\right) ,
\end{equation*}%
\textit{for every local orthonormal frame field }$\ \left\{ e_{1},\cdots \
e_{m};\varphi e_{1},\cdots \ \varphi e_{m}\right\} $\textit{\ \ for} \ $D.$

$3^{\circ })$\textit{\ }

$i)$ \ \textit{non-}$\varphi $-\textit{involutive }if\textit{\ }for every
given $p\in M$%
\begin{equation*}
\left[ X,\ \varphi X\right] _{p}\notin D_{p},\text{ \ \ \ }\forall \text{\ }%
X\in \Gamma \left( D\right) \text{ \ \ \ \ with \ \ }X_{p}\neq 0
\end{equation*}

$ii)$ \textit{non-involutive\ }if for every given $p\in M$%
\begin{equation*}
\left[ X,\ Y\right] _{p}\notin D_{p},\text{ \ \ \ for some \ \ }X,\ Y\in
\Gamma \left( D\right) \text{ \ \ \ \ with \ \ }X_{p}\neq 0,\ \ Y_{p}\neq 0
\end{equation*}

$4^{\circ })$\textit{\ \ non-semi-}$\varphi $-\textit{involutive }if for
every given $p\in M$%
\begin{equation*}
\underset{i=1}{\overset{m}{\sum }}\left[ e_{i},\ \varphi e_{i}\right]
_{p}\notin D_{p},
\end{equation*}%
\textit{for every local orthonormal frame field }$\ \left\{ e_{1},\cdots \
e_{m};\varphi e_{1},\cdots \ \varphi e_{m}\right\} $\textit{\ \ for} \ $D.$

\textbf{Remark }$\mathbf{(2.1)}$\textbf{: }\textit{\ }Note that

$i)$ \ if \ $M$ \ is \ $\varphi $-involutive then clearly it is also semi-$%
\varphi $-involutive.

$ii)$ \ if \ $M$ \ is non-$\varphi $-involutive then it is also
non-involutive.\ 

However,

$iii)$ \ a non-$\varphi $-involutive $M$ \ need not be non-semi-$\varphi $%
-involutive in general.

\textbf{Lemma }$\mathbf{(2.5)}$\textbf{:}

$\mathbf{1^{\circ })}$ \ \textit{The subclasses of almost contact metric
manifolds listed below are }\ $\varphi $-\textit{involutive}:

$\ i)$\textit{\ \ Quasi-}$\mathcal{K}$\textit{-cosymplectic manifolds and
thus its subclasses:}

$\ \bullet $\textit{\ \ Nearly-}$\mathcal{K}$\textit{-cosymplectic manifolds.%
}

$\ \bullet $\textit{\ \ Almost cosymplectic manifolds.}

$\ \bullet $\textit{\ \ Cosymplectic manifolds.}

$ii)$\textit{\ \ Nearly-cosymplectic manifolds.}

$iii)$\textit{\ \ Almost }$\beta $\textit{-Kenmotsu manifolds and thus its
subclasses:}

$\ \bullet $\textit{\ \ }$\beta $\textit{-Kenmotsu manifolds.}

$\ \bullet $\textit{\ \ }$CR$\textit{-integrable almost }$\beta $-\textit{%
Kenmotsu manifolds.}

$\mathbf{2^{\circ })}$\textit{\ \ Almost semi-cosymplectic manifolds and
therefore}$\ $\textit{semi-cosymplectic manifolds are semi-}$\varphi $-%
\textit{involutive.}

\textbf{Proof:}

$\ \mathbf{1^{\circ })/}\left( i\right) :$ \ Let \ $M$ \ be a quasi-$%
\mathcal{K}$-cosymplectic manifold. Then directly from its definition one
gets 
\begin{equation*}
P\left( X,\ X\right) =0,\text{ \ \ \ }\forall \ X\in \Gamma \left( D\right) .
\end{equation*}%
But then Lemma $(1.1)/\left( ii\right) $-$\left( d^{\circ }\right) \ \ $%
gives that%
\begin{equation*}
\left[ X,\ \varphi X\right] =\varphi \left( W_{X}\right) \in \Gamma \left(
D\right) ,\text{ \ \ \ }\forall \text{\ }X\in \Gamma \left( D\right) .
\end{equation*}

$\mathbf{1^{\circ })/}\left( ii\right) :$ \ For a nearly cosymplectic
manifold $M$ \ the same argument used in \ $\left( i\right) $ \ gives the
result.

$\mathbf{1^{\circ })/}\left( iii\right) :$ \ Let \ $M$ \ be an almost $\beta 
$-Kenmotsu manifold. Since \ $d\eta =0$ \ on \ $M$ \ one gets%
\begin{equation*}
\eta \left( \left[ X,\ \varphi X\right] \right) =-2d\eta \left( X,\ \varphi
X\right) =0
\end{equation*}%
and hence the result follows.

$\mathbf{2^{\circ }):}$\textit{\ \ }Let \ $M$ \ be an almost
semi-cosymplectic manifold. Since \ $\delta \Omega =0$ \ on \ $M,$ \ the
equation \ $\left( 2.3\right) $ gives%
\begin{equation*}
\delta \Omega \left( \xi \right) =g\left( -\varphi \left( \nabla _{\xi }\xi
\right) +\underset{i=1}{\overset{m}{\sum }}P\left( e_{i},\ e_{i}\right) ,\ 
\text{\ }\xi \right) =g\left( \underset{i=1}{\overset{m}{\sum }}P\left(
e_{i},\ e_{i}\right) ,\ \text{\ }\xi \right) =0,
\end{equation*}%
which means 
\begin{equation*}
\underset{i=1}{\overset{m}{\sum }}P\left( e_{i},\ e_{i}\right) \in \Gamma
\left( D\right) .
\end{equation*}%
But then, Lemma $(1.1)/\left( \left( 2^{\circ }\right) /\left( iv\right)
\right) \ \ $gives that%
\begin{equation*}
\underset{i=1}{\overset{m}{\sum }}\left[ e_{i},\ \varphi e_{i}\right] \in
\Gamma \left( D\right) ,
\end{equation*}%
which completes the proof.

\textbf{Lemma }$\mathbf{(2.6)}$\textbf{: \ }

$\ \mathbf{1}^{\circ }\mathbf{)}$\textit{\ \ Nearly contact metric manifolds
and thus its subclasses:}

$\ \bullet $\textit{\ \ Quasi contact metric manifolds}

$\ \bullet $\textit{\ \ Quasi-}$\mathcal{K}$\textit{-Sasakian manifolds}

$\ \bullet $\textit{\ \ Nearly-Sasakian manifolds}

$\ \bullet \ $\textit{Contact metric manifolds and thus its subclasses:}

$\ \bullet $\textit{\ \ Strongly pseudoconvex CR-manifolds}

$\ \bullet $\textit{\ \ Sasakian manifolds}

\textit{satisfy that }%
\begin{equation*}
\eta \left( \left[ X,\ \varphi X\right] \right) =2g\left( X,\ X\right) ,\ \
\ \ \ \forall \ X\in \Gamma \left( D\right) .
\end{equation*}%
\textit{So, they are \textbf{non}-}$\varphi $\textit{-involutive and also 
\textbf{non}-\textbf{semi}-}$\varphi $\textit{-involutive.}

$\mathbf{2^{\circ })}$\textit{\ \ Nearly-trans-Sasakian manifolds of type \ }%
$\left( \wp ,\ \theta \right) $ \ \textit{and thus its subclass:}

$\ \bullet $\textit{\ \ Trans-Sasakian manifolds} \textit{of type \ }$\left(
\wp ,\ \theta \right) $

\textit{satisfy that}%
\begin{equation*}
\eta \left( \left[ X,\ \varphi X\right] \right) =2\wp g\left( X,\ X\right)
,\ \ \ \ \ \forall \ X\in \Gamma \left( D\right) .
\end{equation*}%
\textit{So, they are non-}$\varphi $\textit{-involutive and also non-semi-}$%
\varphi $\textit{-involutive.}

$\mathbf{3^{\circ })}$ \ \textit{Almost} \textit{semi-Sasakian manifolds and
thus its subclass:}

$\ \bullet $\textit{\ \ Semi-Sasakian\ manifolds}

\textit{satisfy that}%
\begin{equation*}
\eta \left( \underset{i=1}{\overset{m}{\sum }}\left[ e_{i},\ \varphi e_{i}%
\right] \right) =2m.
\end{equation*}%
\textit{So, they are \textbf{non-semi}-}$\varphi $\textit{-involutive.}

\textbf{Proof:}

$\mathbf{\mathbf{1}^{\circ }\mathbf{)}}$ \ Let \ $M$ \ be a nearly contact
metric manifold. Then by its definition, \ $\forall \ X\in \Gamma \left(
D\right) $ \ $\ \ $ \ \ 
\begin{eqnarray*}
2\Omega \left( X,\ \varphi X\right) &=&d\eta \left( X,\ \varphi X\right)
+d\eta \left( \varphi X,\ \varphi ^{2}X\right) \\
&=&d\eta \left( X,\ \varphi X\right) -d\eta \left( \varphi X,\ X\right) \\
&=&2d\eta \left( X,\ \varphi X\right) =-2g\left( X,\ X\right) .
\end{eqnarray*}%
So we get, $\forall \ p\in M$ \ 
\begin{equation*}
2g\left( X_{p},\ X_{p}\right) =-2d\eta \left( X_{p},\ \varphi X_{p}\right)
=\eta \left( \left[ X,\ \varphi X\right] _{p}\right) \text{ \ \ \ \ \ \ \
for \ }\ X_{p}\neq 0.
\end{equation*}%
Thus we have%
\begin{equation*}
\eta \left( \underset{i=1}{\overset{m}{\sum }}\left[ e_{i},\ \varphi e_{i}%
\right] _{p}\right) =2m
\end{equation*}%
Then the result follows.

$\mathbf{\mathbf{2}^{\circ }\mathbf{)}}$ \ Let \ $M$ \ be a
nearly-trans-Sasakian manifolds of type \ $\left( \wp ,\ \theta \right) .$ \
Then, directly from its definition\ one gets, 
\begin{equation*}
Q\left( X,\ X\right) =\wp 2g\left( X,\ X\right) \xi ,\ \ \ \ \forall \ X\in
\Gamma \left( D\right) .
\end{equation*}%
Noting that 
\begin{equation*}
2P\left( X,\ X\right) =Q\left( X,\ X\right) +Q\left( \varphi X,\ \varphi
X\right)
\end{equation*}%
we get%
\begin{equation*}
P\left( X,\ X\right) =\wp 2g\left( X,\ X\right) \xi ,\ \ \ \ \forall \ X\in
\Gamma \left( D\right) .
\end{equation*}%
But then, using Lemma $(1.1)/\left( \left( 2^{\circ }\right) /\left(
iv\right) \right) ,\ \ \forall \ p\in M$ \ we\ have 
\begin{equation*}
\left[ X,\ \varphi X\right] _{p}=\varphi \left( W_{X_{p}}\right) +\wp
2g\left( X_{p},\ X_{p}\right) \xi
\end{equation*}%
which gives 
\begin{equation*}
\eta \left( \left[ X,\ \varphi X\right] _{p}\right) =\wp 2g\left( X_{p},\
X_{p}\right) \text{\ \ \ \ \ \ for \ }\ X_{p}\neq 0.
\end{equation*}%
Thus we have%
\begin{equation*}
\eta \left( \underset{i=1}{\overset{m}{\sum }}\left[ e_{i},\ \varphi e_{i}%
\right] _{p}\right) =\wp 2m
\end{equation*}%
Noting that $\wp \left( p\right) \neq 0,\ \forall \ p\in M,$ this gives the
result.

$\mathbf{3^{\circ }\mathbf{)}}$\ \ Let \ $M$ \ be an almost semi-Sasakian
manifold\textit{. }Then the equation \ $\left( 2.4\right) $ \ is valid \
(see \ the proof of Lemma $\left( 2.3\right) /\left( \mathbf{2}^{\circ
}\right) )$, that is,\textit{\ }%
\begin{equation*}
\underset{i=1}{\overset{m}{\sum }}P\left( e_{i},\ e_{i}\right) =\varphi
\nabla _{\xi }\xi +2m\xi .
\end{equation*}%
On the other hand, by the virtue of Lemma $\left( 1.1\right) /\left( \left(
2^{\circ }\right) /\left( iv\right) \right) ,$ \ we have

\begin{equation*}
\underset{i=1}{\overset{m}{\sum }}\left[ e_{i},\ \varphi e_{i}\right]
_{p}=\varphi \left( \underset{i=1}{\overset{m}{\sum }}\left(
W_{e_{i}}\right) _{p}\right) +\underset{i=1}{\overset{m}{\sum }}P\left(
e_{i},\ e_{i}\right) _{p},
\end{equation*}%
So, using these last two equations we get%
\begin{equation*}
\underset{i=1}{\overset{m}{\sum }}\left[ e_{i},\ \varphi e_{i}\right]
_{p}=\varphi \left( \underset{i=1}{\overset{m}{\sum }}\left(
W_{e_{i}}\right) _{p}\right) +\varphi \left( \nabla _{\xi }\xi \right)
_{p}+2m\xi _{p}
\end{equation*}%
so that%
\begin{equation*}
\eta \left( \underset{i=1}{\overset{m}{\sum }}\left[ e_{i},\ \varphi e_{i}%
\right] _{p}\right) =2m.
\end{equation*}%
\textit{\ }

\textbf{Lemma }$\mathbf{(2.7)}$\textbf{: \ }

$\mathbf{1}^{\circ }\mathbf{)}$ \ \textit{For a\ \textbf{nearly\textit{-}%
trans}-Sasakian manifold of type }$\left( \wp ,\ \theta \right) $ \textit{%
the following holds:}%
\begin{equation*}
S\left( \varphi X,\ X\right) =\eta \left( W_{X}\right) \xi ,\ \ \ \ \forall
\ X\in \Gamma \left( D\right)
\end{equation*}%
\textit{and therefore}%
\begin{equation*}
\underset{i=1}{\overset{m}{\sum }}S\left( \varphi e_{i},\ e_{i}\right)
=\left( \delta \eta \right) \xi
\end{equation*}

$\mathbf{2}^{\circ }\mathbf{)}$ \ \textit{For a nearly-}$\alpha $-\textit{%
Sasakian manifold the following }%
\begin{equation*}
S\left( \varphi X,\ X\right) =\eta \left( W_{X}\right) \xi ,\ \ \ \ \forall
\ X\in \Gamma \left( D\right)
\end{equation*}%
\textit{and therefore} 
\begin{equation*}
\underset{i=1}{\overset{m}{\sum }}S\left( \varphi e_{i},\ e_{i}\right) =0
\end{equation*}%
\textit{hold.}

$3^{\circ })$ \ \textit{For a trans-Sasakian manifold of type \ }$\left( \wp
,\ \theta \right) $ \ \textit{the following }%
\begin{equation*}
S\left( \varphi X,\ X\right) =\frac{1}{m}\left( \delta \eta \right) g\left(
X,\ X\right) \xi ,\ \ \ \ \forall \ X\in \Gamma \left( D\right)
\end{equation*}%
\textit{and therefore}%
\begin{equation*}
\underset{i=1}{\overset{m}{\sum }}S\left( \varphi e_{i},\ e_{i}\right)
=\left( \delta \eta \right) \xi
\end{equation*}%
\textit{hold}

$4^{\circ })$ \ \textit{For an \textbf{almost}\ }$\beta $-\textit{Kenmotsu
manifold and its subclasses: }

$\bullet $ \ $\beta $\textit{-Kenmotsu manifolds.}

$\bullet $\textit{\ \ }$CR$\textit{-integrable almost }$\beta $-\textit{%
Kenmotsu manifolds }

\textit{the following }%
\begin{equation*}
S\left( \varphi X,\ X\right) =-2\beta g\left( X,\ X\right) \xi ,\ \ \ \
\forall \ X\in \Gamma \left( D\right)
\end{equation*}%
\textit{and therefore}%
\begin{equation*}
\underset{i=1}{\overset{m}{\sum }}S\left( \varphi e_{i},\ e_{i}\right)
=-2\beta m\xi
\end{equation*}%
\textit{hold.}

\textbf{Proof:}

$1^{\circ })$ \ Let \ $M$ \ be nearly trans-Sasakian manifold of type \ $%
\left( \wp ,\ \theta \right) .$ \ Then, as in the proof of Lemma $\left(
2.6\right) /\left( (\mathbf{\mathbf{1}^{\circ }\mathbf{)}/(}iii\mathbf{)}%
\right) $, one gets 
\begin{equation*}
P\left( X,\ X\right) =2\wp g\left( X,\ X\right) \xi ,\ \ \ \ \forall \ X\in
\Gamma \left( D\right) .
\end{equation*}%
This gives, using Lemma $\left( 1.1\right) /\left( \left( 2\mathbf{^{\circ }%
\mathbf{)}/(}iii\right) \right) $, that

\begin{equation*}
S\left( \varphi X,\ X\right) =\eta \left( W_{X}\right) \xi ,\ \ \ \ \forall
\ X\in \Gamma \left( D\right) .
\end{equation*}

$2^{\circ })$ \ Let \ $M$ \ be\textit{\ }a\textit{\ }nearly-$\alpha $%
-Sasakian manifold. Then, by the same argument used in $(\mathbf{\mathbf{1}%
^{\circ }\mathbf{)}}$, one gets that

\begin{equation*}
S\left( \varphi X,\ X\right) =\eta \left( W_{X}\right) \xi ,\ \ \ \ \forall
\ X\in \Gamma \left( D\right) .
\end{equation*}%
and that we get 
\begin{equation*}
\underset{i=1}{\overset{m}{\sum }}S\left( \varphi e_{i},\ e_{i}\right)
=\left( \delta \eta \right) \xi .
\end{equation*}%
However, since \ $\delta \eta =0$ \ on an almost semi-$\alpha $-Sasakian and
therefore on its subclass nearly-$\alpha $-Sasakian manifolds, the result
follows.

$3^{\circ })$ \ Let \ $M$ \ be a trans-Sasakian manifold of type \textit{\ }$%
\left( \wp ,\ \theta \right) .$ Then the required results follow directly
from its definition.

$4^{\circ })$ \ Let \ $M$ \ be\textit{\ }an almost \ $\beta $-\ Kenmotsu
manifold. Then, from the statement \ $\left( 8^{\circ }.4\right) ,$ we get%
\begin{equation*}
S\left( \varphi X,\ X\right) =-2\beta g\left( X,\ X\right) \xi ,\ \ \ \
\forall \ X\in \Gamma \left( D\right) .
\end{equation*}

\textbf{Lemma }$\left( 2.8\right) $\textbf{\ }\ $\left( \left[ 2\right] ,\ 
\text{page 70, equation }\left( 3.2.3\right) \right) $; For a smooth map $F$
\ between Riemannian manifolds the following holds: Fer every local vector
field \ $X,\ Y$%
\begin{equation*}
dF\left( \left[ X,\ Y\right] \right) =\left[ dF\left( X\right) ,\ dF\left(
Y\right) \right] .
\end{equation*}

Blow we give a grouping and abbreviations of the manifolds
considered\bigskip :

Group $\mathcal{A}:$

$\mathcal{A}1:$ $\alpha $-Contact metric manifolds, ($c$); Quasi-contact
metric manifolds, ($q$-$c$); Quasi-$\mathcal{K}$-Sasakian manifold, ($q$-$%
\mathcal{K}$-$S$); Strongly Pseudoconvex $CR$-manifold, ($p$-$CR$); $\alpha $%
-Sasakian manifold,($S$).

$\mathcal{A}2:$ Nearly trans-Sasakian manifold of type ($p$,$\theta $), ($n$-%
$t$-$S$).

$\mathcal{A}3:$ Trans-Sasakian manifold of type ($p$,$\theta $), ($t$-$S$).

$\mathcal{A}4:$ Quasi Sasakian manifold, ($q$-$S$).

$\mathcal{A}_{4}^{a}:$ Almost Quasi Sasakian manifold, ($a$-$q$-$S$).

$\mathcal{A}5:$ Nearly $\alpha $-Contact metric manifold, ($n$-$c$).

$\mathcal{A}6:$ Nearly $\alpha $-Sasakian manifold, ($n$-$S$).

$\mathcal{A}7:$ Semi-Sasakian manifold, ($s$-$S$).

$\mathcal{A}_{7}^{a}:$ Almost semi-Sasakian manifold, ($a$-$s$-$S$).

Group $\mathcal{B}:$

$\mathcal{B}1:$ Cosymplectic manifold, ($Cs$); Almost cosymplectic manifold,
($a$-$Cs$); Nearly cosymplectic manifold, ($n$-$Cs$); Nearly-$\mathcal{K}$%
-cosymplectic manifold, ($n$-$\mathcal{K}$-$Cs$); Quasi- $\mathcal{K}$
-cosymplectic manifold, ($q$-$\mathcal{K}$-$Cs$).

$\mathcal{B}2:$ Semi-cosymplectic manifold, ($s$-$Cs$).

$\mathcal{B}_{2}^{a}:$ Almost semi-cosymplectic manifold, ($a$-$s$-$Cs$).

$\mathcal{B}3:$ $\beta $-Kenmotsu manifold, ($Ksu$); $CR$-integrable almost $%
\beta $-Kenmotsu manifold, ($CR$-$a$-$Ksu$); Almost $\beta $-Kenmotsu
manifold, ($a$-$Ksu).$

Group $\mathcal{C}:$

$\mathcal{C}1:$ Kaehler manifold, ($K$); Nearly-Kaehler manifold, ($n$-$K$);
Almost Kaehler manifold, ($a$-$K$); Quasi-Kaehler manifold, ($q$-$K$).

$\mathcal{C}2:$ Semi-Kaehler manifold, ($s$-$K$).

We also include separately \textbf{tables ( I, II and III )} and \textbf{%
diagrams ( I, II and III ), }written in Microsoft Word, which list and group
some type of manifolds considered, with some of their basic properties and
inclusion relations. We shall be reffering those tables and diagrams
frequently throughout.

\textbf{3) \ Harmonicity}

Let 
\begin{equation*}
f:M=\left( M^{2m+1};\ \varphi ,\ \xi ,\ \eta ,\ g\right) \rightarrow
N=\left( N^{2n+1};\ \phi ,\ \gamma ,\ \sigma ,\ h\right)
\end{equation*}%
be a smooth map between almost contact metric manifolds. We set:

\ 
\begin{equation*}
\ \ f_{\ast }=df\text{ \ \ \ \ and \ }\ E_{i}=f_{\ast }\left( e_{i}\right) ,%
\text{\ \ }
\end{equation*}%
\begin{equation*}
W^{f}\left( X,Y\right) =^{\varphi }W^{f}\left( X,Y\right) =\nabla _{\left(
f_{\ast }X\right) }^{N}\left( f_{\ast }Y\right) +\nabla _{\left( f_{\ast
}\varphi X\right) }^{N}\left( f_{\ast }\varphi Y\right) \ \ \ 
\end{equation*}%
$\ $and%
\begin{equation*}
W^{f}\left( X,X\right) =W_{X}^{f},
\end{equation*}%
\begin{equation*}
U\left( X,Y\right) =^{\varphi }U\left( X,Y\right) =W^{f}\left( X,Y\right)
-f_{\ast }W^{M}\left( X,Y\right) \ \ \ \text{\ \ \ }
\end{equation*}%
and%
\begin{equation*}
^{\varphi }U\left( X,X\right) =U_{X}=W_{X}^{f}-f_{\ast }\left(
W_{X}^{M}\right) .
\end{equation*}%
Recall that $W^{M}\left( X,Y\right) =\nabla _{X}^{M}Y+\nabla _{\varphi
X}^{M}\left( \varphi Y\right) $

The second fundemental form of \ $f$ \ is given by \ 
\begin{equation*}
\left( \mathcal{\bigtriangledown }f_{\ast }\right) \left( X,Y\right) =\nabla
_{\left( f_{\ast }X\right) }^{N}\left( f_{\ast }Y\right) -f_{\ast }\left(
\nabla _{X}^{M}Y\right) .
\end{equation*}%
It is well known that the second fundemental form is symmetric, that is, \ $%
\left( \mathcal{\bigtriangledown }f_{\ast }\right) \left( X,\ Y\right)
=\left( \mathcal{\bigtriangledown }f_{\ast }\right) \left( Y,\ X\right) .$ \
Note that \ 
\begin{equation*}
U\left( X,\ Y\right) =\left( \mathcal{\bigtriangledown }f_{\ast }\right)
\left( X,\ Y\right) +\left( \mathcal{\bigtriangledown }f_{\ast }\right)
\left( \varphi X,\ \varphi Y\right) ,
\end{equation*}%
so that 
\begin{equation*}
^{\varphi }U_{X}=\left( \mathcal{\bigtriangledown }f_{\ast }\right) \left(
X,X\right) +\left( \mathcal{\bigtriangledown }f_{\ast }\right) \left(
\varphi X,\varphi X\right) \ \ \ \ \ \text{and \ \ }U\left( X,\ Y\right)
=U\left( Y,\ X\right) .
\end{equation*}%
The tension field \ $\mathcal{T}\left( f\right) $ (or\ the harmonicity
equation) of \ $f$ \ is given by 
\begin{equation}
\mathcal{T}\left( f\right) :=\ ^{\varphi }U_{\xi }+\overset{m}{\underset{i=1}%
{\sum }}\ ^{\varphi }U_{{\Large e}_{i}}=\left( \mathcal{\bigtriangledown }%
f_{\ast }\right) \left( \xi ,\ \xi \right) +\overset{m}{\underset{i=1}{\sum }%
}\left\{ \left( \mathcal{\bigtriangledown }f_{\ast }\right) \left( {\Large e}%
_{i},{\Large e}_{i}\right) +\left( \mathcal{\bigtriangledown }f_{\ast
}\right) \left( \varphi {\Large e}_{i},\varphi {\Large e}_{i}\right) \right\}
\tag{$3.1$}
\end{equation}%
for any local orthonormal frame field $\ \left\{ e_{1},\cdots \
e_{m};\varphi e_{1},\cdots ~\ \varphi e_{m},\ \xi \right\} $ \ over $M.$ In
cases where \ $M=H^{2m}=\left( H^{2m};J,G\right) $, an almost Hermitian
manifold, \ the equation $\left( 3.1\right) $ becomes 
\begin{equation}
\mathcal{T}\left( f\right) =\overset{m}{\underset{i=1}{\sum }}\ ^{J}U_{%
{\Large e}_{i}},  \tag{$3.2$}
\end{equation}%
where $\ \left\{ e_{1},\cdots \ e_{m};Je_{1},\cdots ~\ Je_{m},\right\} $ \
is any local orthonormal frame field over $H$. \ Note that \ 
\begin{equation*}
U\left( \mu X,\omega Y\right) =\mu \omega U\left( X,Y\right) ,\text{ \ \ \
for any \ }\mu ,\omega \in C^{\infty }\left( M\right)
\end{equation*}

\textbf{Definition }$\left( 3.1\right) :$ The map \ 
\begin{equation*}
f:M=\left( M^{2m+1};\ \varphi ,\ \xi ,\ \eta ,\ g\right) \rightarrow
N=\left( N^{2n+1};\ \phi ,\ \gamma ,\ \sigma ,\ h\right)
\end{equation*}%
is said to be

$\ 1\mathbf{^{\circ }\mathbf{)}}$ $\left( \varphi ,\phi \right) $-\textit{%
holomorphic} (resp: $\left( \varphi ,\phi \right) $-\textit{antiholomorphic}
) if 
\begin{equation*}
f_{\ast }\circ \varphi =\phi \circ f_{\ast }\text{ \ \ \ }\left( \text{resp: 
}f_{\ast }\circ \varphi =-\phi \circ f_{\ast }\right) .
\end{equation*}

$2\mathbf{^{\circ })}/\left( i\right) $ \ \textit{pluriharmonic} if 
\begin{equation*}
U_{X}=0,\ \ \ \forall \ X\in \Gamma \left( TM\right) \text{ \ \ or
equivalently \ \ }U\left( X,\ Y\right) =0;\ \ \ \forall \ X,\ Y\in \Gamma
\left( TM\right) .
\end{equation*}

$2\mathbf{^{\circ })}/\left( ii\right) $ $\ D$-\textit{pluriharmonic} if%
\begin{equation*}
U_{X}=0,\ \ \ \ \ \forall \ X\in \Gamma \left( D^{M}\right) \text{ \ \ or
equivalently \ \ }U\left( X,\ Y\right) =0;\ \ \ \ \ \forall \ X,\ Y\in
\Gamma \left( D^{M}\right) .
\end{equation*}

$3\mathbf{^{\circ }})$ $\ $\textit{harmonic} if 
\begin{equation*}
\mathcal{T}\left( f\right) =0.
\end{equation*}%
We write $\pm \left( \varphi ,\phi \right) $-holomorphic (in short,\ $\left(
\pm \right) $-holomorphic) to mean either $\left( \varphi ,\phi \right) $%
-holomorphic or $\left( \varphi ,\phi \right) $-antiholomorphic.

$4\mathbf{^{\circ }}/\left( i\right) $\ \ \textit{weakly conformal }if there
exists a non-negative function\textit{\ \ }$\mu :M\rightarrow 
%TCIMACRO{\U{211d} }%
%BeginExpansion
\mathbb{R}
%EndExpansion
$\textit{\ \ }such that\textit{\ \ }%
\begin{equation*}
h\left( \left( f_{\ast }\right) _{p}X,\ \left( f_{\ast }\right) _{p}Y\right)
=\mu ^{2}\left( p\right) g\left( X,Y\right) ,\text{ \ \ \ \ }\forall \ X,\
Y\in T_{p}M.
\end{equation*}%
Here the function $\mu $ \ is called the \textit{conformality factor} of \ $%
f.$

We call $\ f$ \ 

$\bullet $ \ \textit{conformal }if \ $\mu \left( p\right) >0,\ \forall \
p\in M$

$\bullet $ \ \textit{homothetic} if $\mu $ is a constant function.

If \ $\dim M>\dim N$ \ then the weakly conformal $\ f$ \ is constant. (see $%
\left[ 2\right] $ Proposition $\left( 2.3.4\right) $)

$4\mathbf{^{\circ }}/\left( ii\right) $ \textit{\ horizontally weakly
conformal }if there exists a non-negative function\textit{\ \ }$\mu
:M\rightarrow 
%TCIMACRO{\U{211d} }%
%BeginExpansion
\mathbb{R}
%EndExpansion
$\textit{\ \ }such that\textit{\ \ }%
\begin{equation*}
h\left( \left( f_{\ast }\right) _{p}X,\ \left( f_{\ast }\right) _{p}Y\right)
=\mu ^{2}\left( p\right) g\left( X,Y\right) ,\text{ \ \ \ \ }\forall \ X,\
Y\in \mathcal{H}_{p}=\left( Ker\left( f_{\ast }\right) _{p}\right) ^{\perp }
\end{equation*}%
and $\ $the differential map $\ \left( \ f_{\ast }\right) _{p}:\mathcal{H}%
_{p}\rightarrow T_{f\left( p\right) }N$ \ \ is surjective at points where $%
\mu \left( p\right) >0.$

Here, the distributions $\mathcal{V=V}^{f}\mathcal{=}Ker\left( f_{\ast
}\right) $ \ and $\ \mathcal{H=H}^{f}\mathcal{=}\left( Ker\left( f_{\ast
}\right) \right) ^{\perp }$ determined by $f$ are called \textit{vertical}
and \textit{horizontal distributions }respectively. Vector fields belonging
to $\mathcal{V}$ (resp: belonging to $\mathcal{H}$ ) are also called \textit{%
vertical} (resp: \textit{horizontal) vector fields. }

We call $\ f$ \ 

$\bullet $ \ \textit{horizontally} \textit{conformal }if \ $\mu \left(
p\right) >0,\ \forall \ p\in M$

$\bullet $ \ \textit{horizontally homothetic} if $\mu $ is a constant
function along horizontal curves, i.e. $\left( d\mu \right) \left( X\right)
=0,\ \ \forall \ X\in \Gamma \left( \mathcal{H}\right) .$ Here the function $%
\mu $ \ is called \textit{dilation} of \ $f.$

$\bullet $ \ \textit{harmonic morphism} \ if it is harmonic and horizontally
weakly conformal.

If \ $\dim M<\dim N$ \ then the horizontally weakly conformal $\ f$ \ is
constant. (see $\left[ 2\right] $ Proposition $\left( 2.4.3\right) $)\ 

We call the points $p$ \ at which

$\bullet $ \ \ $\mu \left( p\right) =0,$ \ critical points of \ $f.$

$\bullet $ \ \ $\mu \left( p\right) >0,$ \ regular points of \ $f.$

\textbf{Remark }$\left( 3.1\right) :$

$1^{\circ })$ \ Clearly, pluriharmonicity implies harmonicity. The converse
is not true in general.

$2^{\circ })$ \ Pluriharmonicity and $D$-Pluriharmonicity do coincide when \ 
$M=H,$ that is, $M$ is an almost Hermitian manifold.

It is not diffucult to prove the following:

\textbf{Lemma }$\left( 3.1\right) :$ \ 

$1^{\circ })$ \ \textit{Let \ \ }$f:\left( M^{2m+1};\ \varphi ,\ \xi ,\ \eta
,\ g\right) \rightarrow \left( N^{2n+1};\ \phi ,\ \gamma ,\ \sigma ,\
h\right) $\textit{\ \ be a\ }$\left( \pm \right) $-\textit{holomorphic map.
Then}

$\ i)$\textit{\ \ }%
\begin{equation*}
f_{\ast }\xi =\lambda \gamma \ \ \ \text{and \ \ }f^{-1}\sigma :=\sigma
f_{\ast }=\lambda \eta
\end{equation*}%
\textit{for some }$\lambda \in C^{\infty }\left( M\right) .$

$ii)$ 
\begin{equation*}
\ f_{\ast }\left( D^{M}\right) \subset D^{N}\ \ \ \text{and\ }\ \ \phi
f_{\ast }\left( D^{M}\right) =\ f_{\ast }\left( D^{M}\right)
\end{equation*}

$2^{\circ })$ \ \textit{Let \ \ }$f:\left( M^{2m+1};\ \varphi ,\ \xi ,\ \eta
,\ g\right) \rightarrow \left( H^{2n};\ J,\ G\right) $\textit{\ \ be a\ }$%
\left( \pm \right) \left( \varphi ,J\right) $-\textit{holomorphic map. Then
\ }$f_{\ast }\xi =0.$\textit{\ \ }

Now on, $\ \lambda $ \ will denote\ throughout the function $\ \lambda \in
C^{\infty }\left( M\right) $ described above in $\left( 1^{\circ }/i\right) $%
.

\textbf{Lemma }$\left( 3.2\right) :$ \ \ 

$1^{\circ })$ \textit{\ For a\ }$\left( \pm \right) \left( \varphi ,\phi
\right) $-\textit{holomorphic map \ }$f:\left( M^{2m+1};\ \varphi ,\ \xi ,\
\eta ,\ g\right) \rightarrow \left( N^{2n+1};\ \phi ,\ \gamma ,\ \sigma ,\
h\right) $ \textit{between almost contact metric manifolds we have: }

$\ i)\ \ \forall \ X\in \Gamma \left( TM\right) $%
\begin{equation*}
f_{\ast }\left[ X,\ \varphi X\right] =\left[ \ f_{\ast }X,\ f_{\ast }\varphi
X\right] =\pm \left[ \ f_{\ast }X,\ \phi f_{\ast }X\right] ,\ 
\end{equation*}

$ii)$ \ $\forall \ X\in \Gamma \left( TM\right) ,$%
\begin{eqnarray*}
U_{X} &=&S^{N}\left( \phi f_{\ast }X,\ f_{\ast }X\right) -\ f_{\ast }\left(
S^{M}\left( \varphi X,X\right) \right) \\
&&+\nabla _{f_{\ast }X}^{N}\left( \sigma \left( f_{\ast }X\right) \gamma
\right) -f_{\ast }\left( \nabla _{X}^{M}\left( \eta \left( X\right) \xi
\right) \right) .\ \ \ 
\end{eqnarray*}%
\ $\ \ $

\textit{In\ particular},\ \ $\forall \ X\in \Gamma \left( D\right) ,$%
\begin{equation*}
U_{X}=S^{N}\left( \phi f_{\ast }X,\ f_{\ast }X\right) -\ f_{\ast }\left(
S^{M}\left( \varphi X,X\right) \right)
\end{equation*}%
\textit{and that} 
\begin{equation*}
\overset{m}{\underset{i=1}{\sum }}U_{e_{i}}=\overset{m}{\underset{i=1}{\sum }%
}\left\{ S^{N}\left( \phi E_{i},\ E_{i}\right) -\ f_{\ast }S^{M}\left(
\varphi e_{i},e_{i}\right) \right\} ,\ 
\end{equation*}

$2^{\circ })$ \ \textit{For a\ }$\left( \pm \right) \left( \varphi ,J\right) 
$-\textit{holomorphic map \ }$f:\left( M^{2m+1};\ \varphi ,\ \xi ,\ \eta ,\
g\right) \rightarrow \left( H^{2n};\ J,\ G\right) $ \textit{from an almost
contact metric manifold into an almost Hermitian maifold we have: }

$i)\ \ \forall \ X\in \Gamma \left( TM\right) $%
\begin{equation*}
f_{\ast }\left[ X,\ \varphi X\right] =\left[ \ f_{\ast }X,\ f_{\ast }\varphi
X\right] =\pm \left[ \ f_{\ast }X,\ Jf_{\ast }X\right] ,\ 
\end{equation*}

$ii)$ \ $\forall \ X\in \Gamma \left( TM\right) ,$%
\begin{equation*}
U_{X}=S^{N}\left( Jf_{\ast }X,\ f_{\ast }X\right) -\ f_{\ast }\left(
S^{M}\left( \varphi X,X\right) \right) -\eta \left( X\right) f_{\ast }\left(
\nabla _{X}\xi \right)
\end{equation*}%
\ 

\textit{In\ particular},\ \ $\forall \ X\in \Gamma \left( D\right) ,$%
\begin{equation*}
U_{X}=S^{N}\left( Jf_{\ast }X,\ f_{\ast }X\right) -\ f_{\ast }\left(
S^{M}\left( \varphi X,X\right) \right)
\end{equation*}%
\textit{and that} 
\begin{equation*}
\overset{m}{\underset{i=1}{\sum }}U_{e_{i}}=\overset{m}{\underset{i=1}{\sum }%
}\left\{ S^{N}\left( JE_{i},\ E_{i}\right) -\ f_{\ast }S^{M}\left( \varphi
e_{i},e_{i}\right) \right\} ,\ 
\end{equation*}

$3^{\circ })$ \ \textit{For a\ }$\left( \pm \right) \left( J,\phi \right) $-%
\textit{holomorphic map \ }$f:\left( H^{2m};\ J,\ G\right) \rightarrow
\left( N^{2n+1};\ \phi ,\ \gamma ,\ \sigma ,\ h\right) $ \textit{from an
almost Hermitian maifold into an almost contact metric maifold} we have:

$i)\ \ \forall \ X\in \Gamma \left( TM\right) $%
\begin{equation*}
f_{\ast }\left[ X,\ JX\right] =\left[ \ f_{\ast }X,\ f_{\ast }JX\right] =\pm %
\left[ \ f_{\ast }X,\ \phi f_{\ast }X\right] ,\ 
\end{equation*}

$ii)$ \ $\forall \ X\in \Gamma \left( TM\right) ,$%
\begin{equation*}
U_{X}=S^{N}\left( \phi f_{\ast }X,\ f_{\ast }X\right) -\ f_{\ast }\left(
S^{M}\left( JX,X\right) \right) \ \ \ 
\end{equation*}%
\textit{so that}%
\begin{equation*}
\overset{m}{\underset{i=1}{\sum }}U_{e_{i}}=\overset{m}{\underset{i=1}{\sum }%
}\left\{ S^{N}\left( \phi E_{i},\ E_{i}\right) -\ f_{\ast }S^{M}\left(
Je_{i},e_{i}\right) \right\} ,\ 
\end{equation*}%
\textit{for any orthonormal frame field }$\left\{ e_{1},\cdots ,e_{m};\
Je_{1},\cdots ,Je_{m}\right\} $\textit{\ over} $H^{2m}.$

\textbf{Proof:}

$1^{\circ })/\left( i\right) $ :\ Note that, by Lemma $\left( 2.8\right) $,
\ 
\begin{equation*}
f_{\ast }\left[ X,\ \varphi X\right] =\left[ \ f_{\ast }X,\ f_{\ast }\varphi
X\right] .
\end{equation*}%
Then, $\left( \pm \right) $-holomorphicity of $\ M$ \ gives the result.

$1^{\circ })/\left( ii\right) $ \ :\ Note that, by the $\left( \pm \right) $%
-holomorphicty of \ $f,$ \ 
\begin{equation*}
W_{X}^{f}=\nabla _{\left( f_{\ast }X\right) }\left( f_{\ast }X\right)
+\nabla _{\left( \phi f_{\ast }X\right) }\left( \phi f_{\ast }X\right) .
\end{equation*}%
This gives, by Lemma $\left( 1.1\right) /\left( 1^{\circ }\right) ,$ \ 
\begin{equation}
W_{X}^{f}=S^{N}\left( \phi f_{\ast }X,\ f_{\ast }X\right) +\phi \left[ \
f_{\ast }X,\ \phi f_{\ast }X\right] +\nabla _{f_{\ast }X}^{N}\left( \sigma
\left( \ f_{\ast }X\right) \gamma \right) .  \tag{$3.3$}
\end{equation}%
On the other hand 
\begin{equation*}
f_{\ast }\left( W_{X}^{M}\right) =f_{\ast }S^{M}\left( \varphi X,\ X\right)
+f_{\ast }\varphi \left[ X,\ \varphi X\right] +f_{\ast }\left( \nabla
_{X}^{M}\left( \eta \left( X\right) \xi \right) \right) .
\end{equation*}%
This gives, using the \ $\left( \pm \right) $-holomorphicty of \ $f$ \ and
part $(i),$ \ \ 
\begin{equation*}
f_{\ast }\left( W_{X}^{M}\right) =f_{\ast }S^{M}\left( \varphi X,\ X\right)
+\phi \left[ f_{\ast }X,\ \phi f_{\ast }X\right] +f_{\ast }\left( \nabla
_{X}\left( \eta \left( X\right) \xi \right) \right) .
\end{equation*}%
Cobining this with \ $(3.3)$ we get 
\begin{eqnarray*}
U_{X} &=&W_{X}^{f}-f_{\ast }\left( W_{X}^{M}\right) =S^{N}\left( \phi
f_{\ast }X,\ f_{\ast }X\right) -\ f_{\ast }\left( S^{M}\left( \varphi X,\
X\right) \right) \\
&&\ \ \ \ \ \ \ \ \ \ \ \ \ \ \ \ \ \ \ \ \ \ \ \ +\nabla _{f_{\ast
}X}\left( \sigma \left( \ f_{\ast }X\right) \gamma \right) -f_{\ast }\left(
\nabla _{X}\left( \eta \left( X\right) \xi \right) \right)
\end{eqnarray*}%
which is the required result.

In particular, observe that \ \ $\forall \ X\in \Gamma \left( D^{M}\right) $
\ we have \ $\ f_{\ast }X\in \Gamma \left( D^{N}\right) $ \ by the $\left(
\pm \right) $-holomorphicty of \ $f\ \ \ $and therefore $\ \eta \left(
X\right) =0$ \ and $\ \sigma \left( f_{\ast }X\right) =0.$ So we get 
\begin{equation*}
U_{X}=S^{N}\left( \phi f_{\ast }X,\ f_{\ast }X\right) -\ f_{\ast }\left(
S^{M}\left( \varphi X,\ X\right) \right) ,\ \ \forall \ X\in \Gamma \left(
D^{M}\right)
\end{equation*}%
\ \ ,

$2^{\circ })/\left( i\right) $ and $\left( ii\right) $ : Noting that $%
f_{\ast }\xi =0$ and adopting the proof in $(1^{\circ })/\left( \left(
i\right) ,\left( ii\right) \right) ,$ \ for $\left( \pm \right) \left(
\varphi ,J\right) $-holomorphicity will give the result\textit{.}

$3^{\circ })/\left( i\right) $ and $\left( ii\right) $ : Adopting the proof
in $(1^{\circ })/\left( \left( i\right) ,\left( ii\right) \right) ,$ \ for $%
\left( \pm \right) \left( J,\phi \right) $-holomorphicity will give the
result\textit{.}

\textbf{Proposition }$\left( 3.1\right) :$

$1^{\circ })$ \ \textit{For a\ }$\left( \pm \right) $\textit{-holomorphic
map \ }$f:\left( M^{2m+1};\ \varphi ,\ \xi ,\ \eta ,\ g\right) \rightarrow
\left( N^{2n+1};\ \phi ,\ \gamma ,\ \sigma ,\ h\right) $\textit{\ between
almost contact metric maifolds, the tension field \ }$\mathcal{T}\left(
f\right) $\textit{\ \ of \ }$f$\textit{\ \ takes the form:}%
\begin{equation}
\mathcal{T}\left( f\right) =U_{\xi }+\overset{m}{\underset{i=1}{\sum }}%
\left\{ S^{N}\left( \phi E_{i},\ E_{i}\right) -\ f_{\ast }S^{M}\left(
\varphi e_{i},e_{i}\right) \right\}  \tag{$3.4$}
\end{equation}%
\textit{or equivalently }%
\begin{equation*}
\mathcal{T}\left( f\right) =\nabla _{\left( f_{\ast }\xi \right) }\left(
f_{\ast }\xi \right) -\ f_{\ast }\left\{ \varphi \left( \delta \varphi
\right) +\left( \delta \eta \right) \xi \right\} +\overset{m}{\underset{i=1}{%
\sum }}S^{N}\left( \phi E_{i},E_{i}\right) .
\end{equation*}

$2^{\circ })$\textit{\ \ For a\ }$\left( \pm \right) $\textit{-holomorphic
map \ }$f:\left( M^{2m+1};\ \varphi ,\ \xi ,\ \eta ,\ g\right) \rightarrow
\left( H_{1}^{2n};\ J_{1},\ G_{1}\right) $\textit{\ from an almost contact
metric manifold into an almost Hermitian manifold, the tension field \ }$%
\mathcal{T}\left( f\right) $\textit{\ \ of \ }$f$\textit{\ \ takes the form: 
}%
\begin{equation*}
\mathcal{T}\left( f\right) =-f_{\ast }\left( \nabla _{\xi }\xi \right) +%
\overset{m}{\underset{i=1}{\sum }}\left\{ S^{H_{1}}\left( J_{1}E_{i},\
E_{i}\right) -\ f_{\ast }S^{M}\left( \varphi e_{i},e_{i}\right) \right\}
\end{equation*}

$3^{\circ })$\textit{\ \ For a\ }$\left( \pm \right) $\textit{-holomorphic
map \ }$f:\left( H^{2m};\ J,\ G\right) \rightarrow \left( N^{2n+1};\ \phi ,\
\gamma ,\ \sigma ,\ h\right) $\textit{\ from an almost Hermitian manifold
into an almost contact metric manifold, the tension field \ }$\mathcal{T}%
\left( f\right) $\textit{\ \ of \ }$f$\textit{\ \ takes the form:}%
\begin{equation*}
\mathcal{T}\left( f\right) =\overset{m}{\underset{i=1}{\sum }}\left\{
S^{N}\left( \phi E_{i},\ E_{i}\right) -\ f_{\ast }S^{H}\left(
Je_{i},e_{i}\right) \right\} .
\end{equation*}

\textbf{Proof: }

$1^{\circ })$\textit{\ \ }Using\textbf{\ }Lemma $(3.2)/\left( \left(
1^{\circ }\right) /\left( ii\right) \right) $ and the harmoicity equation $%
(3.1)$\ gives the first part. For the second part, note that 
\begin{equation*}
\delta \varphi =-\varphi \left( \nabla _{\xi }\xi \right) +\overset{m}{%
\underset{i=1}{\sum }}P^{M}\left( e_{i},e_{i}\right)
\end{equation*}%
and therefore, since $\ \nabla _{\xi }\xi \in \Gamma \left( D\right) ,$ one
gets%
\begin{equation}
\varphi \left( \delta \varphi \right) =\nabla _{\xi }\xi +\varphi \overset{m}%
{\underset{i=1}{\sum }}P^{M}\left( e_{i},e_{i}\right) .  \tag{$3.5$}
\end{equation}%
On the other hand, using Lemma $(1.1)/\left( \left( 2^{\circ }\right)
/\left( iii\right) \right) ,$ we get 
\begin{equation*}
\overset{m}{\underset{i=1}{\sum }}S^{M}\left( \varphi e_{i},e_{i}\right)
=\left( \varphi \overset{m}{\underset{i=1}{\sum }}P^{M}\left(
e_{i},e_{i}\right) \right) +\left( \delta \eta \right) \xi .
\end{equation*}%
So, this gives ( by using $(3.5)),$ 
\begin{equation*}
\overset{m}{\underset{i=1}{\sum }}S^{M}\left( \varphi e_{i},e_{i}\right)
=\varphi \left( \delta \varphi \right) +\left( \delta \eta \right) \xi
-\nabla _{\xi }\xi .
\end{equation*}%
Inserting this in $(3.4),$ we get 
\begin{equation}
\mathcal{T}\left( f\right) =U_{\xi }+\ f_{\ast }\left( \nabla _{\xi }\xi
\right) -\ f_{\ast }\left\{ \varphi \left( \delta \varphi \right) +\left(
\delta \eta \right) \xi \right\} +\overset{m}{\underset{i=1}{\sum }}%
S^{N}\left( \phi E_{i},\ E_{i}\right) .  \tag{$3.6$}
\end{equation}%
But note that 
\begin{equation*}
U_{\xi }+\ f_{\ast }\left( \nabla _{\xi }\xi \right) =\nabla _{\left(
f_{\ast }\xi \right) }\left( f_{\ast }\xi \right) ,
\end{equation*}%
so, using this in the equation $(3.6),$\ gives the result.

$2^{\circ })$ \ Noting that $\ f_{\ast }\xi =0,$ we get $U_{\xi }=f_{\ast
}\left( \nabla _{\xi }\xi \right) .$ So, by the same argument used in $%
\left( 1^{\circ }\right) $ above, the result follows.

$3^{\circ })$ \ Note here that the set $\left\{ e_{1},\cdots ,e_{m};\
Je_{1},\cdots ,Je_{m}\right\} $ constitutes an orthonormal frame field on $%
\left( H^{2m};\ J,\ G\right) .$ So, the result follows from the harmoicity
equation $(3.1)$ .

\textbf{Lemma }$\left( 3.3\right) :$\textit{\ \ Let }$\ \ f:M=\left(
M^{2m+1};\ \varphi ,\ \xi ,\ \eta ,\ g\right) \rightarrow N$\textit{\ \ be a 
\textbf{non-constant}, }$\left( \pm \right) $\textit{-holomorphic, \textbf{%
horizontally weakly conformal} map with dilation \ }$\mu $\textit{\ \ from
an almost contact metric manifold \ }$M$\textit{\ \ into either an almost
contact metric manifold \ }$N=\left( N^{2n+1};\ \phi ,\ \gamma ,\ \sigma ,\
h\right) $\textit{\ or an almost Hermitian manifold }$\ N=\left(
H_{1}^{2n};\ J_{1},\ G_{1}\right) $\textit{. \ Then}

$1^{\circ })$\textit{\ \ }$\ m\geq n$

$2^{\circ })$\textit{\ \ When \ }$N=\left( N^{2n+1};\ \phi ,\ \gamma ,\
\sigma ,\ h\right) $\textit{\ \ is an almost contact metric manifold }

$\ i)\ \ \xi $ $\ $is a horizotal vector field, that is, $\ \xi $ $\in 
\mathcal{H}^{f}$ so that the vertical component $\xi _{\mathcal{V}}$\ of $%
\xi $ is identically zero or equivalently \ $\mathcal{V}^{f}\subset D^{M}.$%
\textit{\ \ }

$ii)\ \ \mu =\left\vert \lambda \right\vert .$\textit{\ }

$3^{\circ })$\textit{\ \ \ When \ }$N=\left( H_{1}^{2n};\ J_{1},\
G_{1}\right) $\textit{\ is an\ almost Hermitian manifold }$\ $

$\ \bullet \ \ \xi $ $\ $is a vertical vector field, that is, $\ \xi $ $\in 
\mathcal{V}^{f}$ so that the horizontal component $\xi _{\mathcal{H}}$ of $%
\xi $ is identically zero or equivalently \ $\mathcal{H}^{f}\subset D^{M}.$%
\textit{\ \ }

$4^{\circ })$\textit{\ \ }$\left( c.f.\ \left[ 2\right] ,\ \text{Lemma }%
\left( 2.4.4\right) \right) $ \textit{There is a local orthonormal frame
field \ }%
\begin{equation*}
\left\{ e_{_{1}},\cdots ,e_{_{m}};\ \varphi e_{_{1}},\cdots ,\varphi
e_{_{m}}\right\}
\end{equation*}%
\textit{for \ }$D^{M}$\textit{\ \ such that the set}%
\begin{equation*}
\left\{ \ v_{_{1}},\cdots ,v_{n};\ \psi v_{_{1}},\cdots ,\psi v_{n}\right\}
\end{equation*}%
\textit{forms a local orthnormal frame field for}

$\bullet $\textit{\ \ }$D^{N}$\textit{\ when \ }$N=\left( N^{2n+1};\ \phi ,\
\gamma ,\ \sigma ,\ h\right) $\textit{\ in which case }$\psi =\phi ,$

$\bullet $\textit{\ \ }$TH$\textit{\ when \ }$N=H_{1}=\left( H_{1}^{2n};\
J_{1},\ G_{1}\right) $\textit{\ in which case }$\psi =J_{1},$

\textit{where }%
\begin{equation*}
f_{\ast }\left( e_{i}\right) =\left\{ 
\begin{array}{cc}
\mu v_{_{i}}, & i=1,\cdots ,n \\ 
0 & i>n%
\end{array}%
\right\} ,
\end{equation*}

Further

$i)$ \ \textit{\ \ }%
\begin{equation*}
\overset{m}{\underset{i=1}{\sum }}S^{N}\left( \psi E_{i},\ E_{i}\right) =\mu
^{2}\overset{n}{\underset{i=1}{\sum }}S^{N}\left( \psi v_{i},\ v_{i}\right)
\end{equation*}

$ii)$ \ 
\begin{equation*}
\overset{m}{\underset{i=1}{\sum }}h\left( E_{i},\ E_{i}\right) =\overset{n}{%
\underset{i=1}{\sum }}h\left( E_{i},\ E_{i}\right) =n\mu ^{2}
\end{equation*}

\textbf{Proof:}

$\bigskip 1^{\circ }):$ \ \ For \ $\ m<n,$ the differential $\ f_{\ast }$
cannot be surjective at regular points. Therefore \ $f$ \ cannot be
horizontally weakly conformal unless it is constant.

$2^{\circ }):$ \ We may write

\begin{equation*}
\xi =\xi _{\mathcal{V}}+\xi _{\mathcal{H}}\ \ \ \ with\ \ \xi _{\mathcal{V}%
}\perp \xi _{\mathcal{H}}
\end{equation*}%
$\ \ $where $\xi _{\mathcal{V}}\in \mathcal{V}$ \ and \ $\xi _{\mathcal{H}%
}\in \mathcal{H}$ \ are the vertical and horizontal components of $\xi .$

\ 

$\left( i\right) :\ \ $For horizontalitiy of $\xi $ \ it is enough to show
that $\xi _{\mathcal{V}}=0.$ We do this$\ $in a few steps:

$\ \bullet $ \ $\xi $ \ cannot be vertical at any point, $i.e.$ $\ \xi
_{p}\notin \mathcal{V}$ \ (so that $\left( \xi _{\mathcal{H}}\right)
_{p}\neq 0)$ at any point $p\in M.$

For otherwise, suppose $f_{\ast }\xi _{p}=0$ for some $p\in M.$ But \ then,
by the virtue of Lemma $\left( 3.1\right) /1^{\circ },$ we see that $\gamma
_{f\left( p\right) }\notin f_{\ast }\left( T_{p}M\right) $ and therefore $%
f_{\ast }$ can not be surjective. So, we must have that 
\begin{equation*}
\ \xi _{p}\notin \mathcal{V}\ ,\ \ \forall \ p\in M.
\end{equation*}

$\bullet $ $\ \mathcal{V}^{f}$\ \ is closed under $\varphi ,\ \ i.e.\ \
\varphi \mathcal{V\subseteq V}.$

Indeed, for any \ $X\in \mathcal{V}$ \ we have%
\begin{equation*}
f_{\ast }\left( \varphi X\right) =\phi \left( f_{\ast }\left( X\right)
\right) =0
\end{equation*}%
so that \ $\varphi \left( X\right) \in \mathcal{V}.$

$\ \bullet $ $\ \varphi \left( \xi _{\mathcal{H}}\right) \in \mathcal{H}.$

Indeed, observe that for any $\ X\in \mathcal{V}$

\begin{equation*}
g\left( \varphi \left( \xi _{\mathcal{H}}\right) ,\ X\right) =-g\left( \xi _{%
\mathcal{H}},\ \varphi X\right) =0
\end{equation*}%
since $\varphi X\in \mathcal{V}.$

$\bullet $ \ Finely, suppose that $\left( \xi _{\mathcal{V}}\right) _{p}\neq
0,$ for$~$\ some $\ p\in M.$ \ Since \ 
\begin{equation*}
\varphi \left( \xi \right) =\varphi \left( \xi _{\mathcal{V}}\right)
+\varphi \left( \xi _{\mathcal{H}}\right) =0
\end{equation*}%
with \ 
\begin{equation*}
\varphi \left( \xi _{\mathcal{V}}\right) \in \mathcal{V}\text{ \ \ and \ \ }%
\varphi \left( \xi _{\mathcal{H}}\right) \in \mathcal{H},
\end{equation*}%
we get $\ \varphi \left( \xi _{\mathcal{V}}\right) _{p}=0.$ This is possible
only when $\left( \xi _{\mathcal{V}}\right) _{p}=s\xi _{p}$ \ for some
number $s\neq 0$ $,$ \ that is, $\xi _{p}$ \ is vertical. But $\ \xi $ \
cannot be vertical at any point $p\in M$. \ So we must have 
\begin{equation*}
\left( \xi _{\mathcal{V}}\right) _{p}=0,\ \ \forall \ p\in M
\end{equation*}%
that is, $\xi $ \ is horizontal.

$2^{\circ })/\ \left( ii\right) :\ $To show $\ \mu =|\lambda |,$ note that $%
\xi $\ is a horizontal vector field. So we have 
\begin{equation*}
h\left( \ f_{\ast }\xi ,\ \ f_{\ast }\xi \right) =\mu ^{2}g\left( \xi ,\ \xi
\right) =\mu ^{2}.
\end{equation*}%
\ On the other hand, 
\begin{equation*}
h\left( \ f_{\ast }\xi ,\ \ f_{\ast }\xi \right) =h\left( \lambda \gamma ,\
\lambda \gamma \right) =\lambda ^{2}
\end{equation*}%
since $\ f_{\ast }\xi =\lambda \gamma .$ So the result follows.

$3^{\circ })$ \ Note that 
\begin{equation*}
0=\ f_{\ast }\left( \varphi \left( \xi \right) \right) =\ J_{1}\left(
f_{\ast }\left( \xi \right) \right)
\end{equation*}%
which gives that $\ f_{\ast }\left( \xi \right) =0.$ \ That is, $\xi $ $\in 
\mathcal{V}^{f}.$

$4^{\circ })$ \ Since the set 
\begin{equation*}
\left\{ v_{1},\cdots ,v_{n};\ \psi v_{1},\cdots ,\psi v_{n}\right\}
\end{equation*}%
is a local orthonormal frame field for $\ D^{N}$\ with \ $E_{i}=f_{\ast
}\left( e_{i}\right) =\mu v_{i}$ and $S^{N}$ is a tensor field, we have%
\begin{equation*}
S^{N}\left( \psi E_{i},\ E_{i}\right) =\mu ^{2}S^{N}\left( \psi v_{i},\
v_{i}\right)
\end{equation*}%
and%
\begin{equation*}
h\left( E_{i},\ E_{i}\right) =\mu ^{2}g\left( e_{i},\ e_{i}\right) =\mu ^{2},
\end{equation*}%
from which $4^{\circ })/\left( (ii),(iii)\right) $\ follow.

\textbf{Proposition }$\left( 3.2\right) :$\textit{\ }

$1^{\circ })$\textit{\ \ }$\left( c.f.\ \left[ 14\right] ,\text{\textit{%
Theorem} }\left( 2.2\right) \right) $ \textit{If \ }$M=\left( M^{2m+1};\
\varphi ,\ \xi ,\ \eta ,\ g\right) \in $\textit{\ }$\mathcal{B}$\textit{\ \
and \ }$N=\left( N^{2n+1};\ \phi ,\ \gamma ,\ \sigma ,\ h\right) \in 
\mathcal{A}_{1}\cup \mathcal{A}_{2}\cup \mathcal{A}_{3}\cup \mathcal{A}%
_{5}\cup \mathcal{A}_{6}$\textit{\ \ then any }$\left( \pm \right) $\textit{%
-holomorphic map }$f:M\rightarrow N$\textit{\ \ \ \ is constant along }$%
\varphi \left( TM\right) =D^{M},$\textit{\ \ (i.e. \ }$\ f_{\ast }\left(
X\right) =0,\ \forall \ X\in \Gamma \left( D^{M}\right) )$\textit{\ and
hence }

$i)$ \ \textit{the tension field \ }$\mathcal{T}\left( f\right) $\textit{\ \
of \ }$f$\textit{\ \ takes the form:}%
\begin{equation}
\mathcal{T}\left( f\right) =\nabla _{\left( f_{\ast }\xi \right) }^{N}\left(
f_{\ast }\xi \right) -\left( \delta \eta \right) \lambda \gamma . 
\tag{$3.7$}
\end{equation}

$ii)$ \ the non-constant $\ f$ \textit{\ can be neither weakly conformal nor 
\textbf{horizontally} weakly conformal.}

$2^{\circ })$\textit{\ \ If \ }$M=\left( M^{2m+1};\ \varphi ,\ \xi ,\ \eta
,\ g\right) \in $\textit{\ }$\mathcal{B}$\textit{\ \ and \ }$N=\left(
N^{2n+1};\ \phi ,\ \gamma ,\ \sigma ,\ h\right) \in \mathcal{A}_{7}^{a}$%
\textit{\ then there is no \textbf{non-constant} }$\left( \pm \right) $%
\textit{-holomorphic \textbf{horizontally} weakly conformal map }$%
f:M\rightarrow N\ .$\ 

$3^{\circ })$\textit{\ If \ }$M=H=\left( H^{2m},J\ ,G\right) $\textit{\ \ is
an almost Hermitian manifold and \ }$N=\left( N^{2n+1};\ \phi ,\ \gamma ,\
\sigma ,\ h\right) \in \mathcal{A}_{1}\cup \mathcal{A}_{2}\cup \mathcal{A}%
_{3}\cup \mathcal{A}_{5}\cup \mathcal{A}_{6}$\textit{\ \ \ then there is no
non-constant }$\left( \pm \right) $\textit{-holomorphic map }$f:H\rightarrow
N$

$4^{\circ })$\ \textit{\ \ If \ }$M=H=\left( H^{2m},J\ ,G\right) $\textit{\
\ is an almost Hermitian manifold and \ }$N=\left( N^{2n+1};\ \phi ,\ \gamma
,\ \sigma ,\ h\right) \ \ $\textit{is an\ almost contact metric manifold\
then there is no non-constant }$\left( \pm \right) $\textit{-holomorphic 
\textbf{horizontally }weakly conformal\textbf{\ }map }$f:H\rightarrow N$

$5^{\circ })$\ \textit{\ \ If \ }$M=\left( M^{2m+1};\ \varphi ,\ \xi ,\ \eta
,\ g\right) $\textit{\ \ is an\ almost contact metric manifold and }$%
N=H=\left( H^{2m},J\ ,G\right) $ \ \textit{is an almost Hermitian manifold
then there is no non-constant }$\left( \pm \right) $\textit{-holomorphic
weakly conformal map }$f:M\rightarrow H$

({\small Throughout we shall be writing, such as }$N\in A_{1}\cup A_{2}\cup
A_{3}\cup A_{5}\cup A_{6}${\small \ (as it is done above for example), in
stead of listing all type of manifolds }$N${\small \ might be, for the sake
of simplicity. For the meaning of }$N\in A_{1}\cup A_{2}\cup A_{3}\cup
A_{5}\cup A_{6}${\small \ we often refer the tables).}

\textbf{Proof : }

\bigskip $1^{\circ })$\textit{\ \ }For \textbf{any} local frame field $%
\left\{ e_{_{1}},\cdots ,e_{_{m}};\varphi e_{_{1}},\cdots ,\varphi
e_{_{m}}\right\} $ \ for\ $D^{M}$, \ we have \ 
\begin{equation*}
L_{\varphi }=\overset{m}{\underset{j=1}{\sum }}\left[ \varphi e_{j},e_{j}%
\right] \text{ }\in \Gamma \left( D^{M}\right) .
\end{equation*}%
since\ $M$ \ is \textbf{semi}-$\varphi $-involutive when $M$\textit{\ \ }$%
\in \mathcal{B}$. \ Therefore, by Lemma $\left( 3.1\right) /\left( 1^{\circ
}-ii\right) ,$ \ we get 
\begin{equation}
f_{\ast }\left( L_{\varphi }\right) \text{ }\in \Gamma \left( D^{N}\right) 
\tag{$3.8$}
\end{equation}%
since \ \textit{\ }$f$ \ is \textit{\ }$\left( \pm \right) $-holomorphic. On
the other hand, for $\ N=\left( N^{2n+1};\ \phi ,\ \gamma ,\ \sigma ,\
h\right) \in \mathcal{A}_{1}\cup \mathcal{A}_{2}\cup \mathcal{A}_{3}\cup 
\mathcal{A}_{5}\cup \mathcal{A}_{6},$ \ suppose $\ f_{\ast }\left( v\right)
\neq 0$ \ for some nonezero \ $v\in D_{p}^{M}$ \ and for some $\ p\in M\ .$
\ W.l.o.g. choose a local frame field \ $\left\{ e_{_{1}},\cdots
,e_{_{m}};\varphi e_{_{1}},\cdots ,\varphi e_{_{m}}\right\} $ \ for\ $D^{M}$
\ with \ $e_{_{1}}\left( p\right) =v.$ Noting$\ $that $\ f_{\ast
}e_{1}=E_{_{1}}\neq 0$ $\ $with $\ E_{_{1}}\in \Gamma \left( D^{N}\right) $
and $N$ \ is \textbf{non}-$\phi $-involutive, \ (see the TABLE-\textbf{II} )
we get 
\begin{equation*}
f_{\ast }\left( \left[ \varphi e_{1},e_{1}\right] _{p}\right) =\pm \left[
\phi E_{_{1}},\ E_{_{1}}\right] _{f\left( p\right) }\notin \Gamma \left(
D^{N}\right) .\text{ }
\end{equation*}%
But then, we may write 
\begin{equation*}
\left[ \phi E_{j},\ E_{j}\right] _{f\left( p\right) }=\left(
Z_{j}+r_{j}\gamma _{f\left( p\right) }\right) \text{ }
\end{equation*}%
where \ $Z_{j}\in D_{f\left( p\right) }^{N},\ \ $\ and (by the virtue of
Lemma $\left( 2.6\right) $), \ 
\begin{equation*}
\begin{array}{ll}
r_{j}= & \left\{ 
\begin{array}{ll}
2h\left( E_{_{j}},\ E_{j}\right) , & N\in A_{1}\cup A_{5}\cup A_{6}\mathit{%
.\ } \\ 
2\wp h\left( E_{j},\ E_{_{j}}\right) & N\in A_{2}\cup A_{3}%
\end{array}%
\right\}%
\end{array}%
\end{equation*}%
So,%
\begin{equation*}
\overset{m}{\underset{j=1}{\sum }}\left[ \phi E_{j},\ E_{j}\right] _{f\left(
p\right) }=\overset{m}{\underset{j=1}{\sum }}\left( Z_{j}\right) _{_{f\left(
p\right) }}+\left( \overset{m}{\underset{j=1}{\sum }}\ r_{j}\right) \gamma
_{f\left( p\right) }
\end{equation*}%
with \ $\overset{m}{\underset{j=1}{\sum }}\left( Z_{j}\right) _{_{f\left(
p\right) }}\in D_{_{f\left( p\right) }}^{N}$ \ and \ $r_{j}\geq 0,\ \
\forall \ j.$ Since $\wp \left( q\right) \neq 0$ \ $\forall \ q\in N,$ \ and 
$\ h\left( E_{_{1}},\ E_{_{1}}\right) >0$ \ we get%
\begin{equation*}
\overset{m}{\underset{j=1}{\sum }}\ r_{j}>0.
\end{equation*}%
This means that \ \ 
\begin{equation*}
f_{\ast }\left( L_{\varphi }\right) _{p}=\pm \overset{m}{\underset{i=1}{\sum 
}}\left[ \phi E_{i},\ E_{i}\right] _{f\left( p\right) }\notin D_{_{f\left(
p\right) }}^{N}
\end{equation*}%
\ \ which contradicts with $\left( 3.8\right) .$ So $\forall \ p\in M$ \ we
must have $f_{\ast }\left( X\right) =0,$ \ $\forall \ X\in D_{p}^{M}$\ \
which completes the proof of first part of the assertion. $\ \ $

$1^{\circ })\left( i\right) :$ \ In order to show that 
\begin{equation*}
\mathcal{T}\left( f\right) =\nabla _{\left( f_{\ast }\xi \right) }\left(
f_{\ast }\xi \right) -\left( \delta \eta \right) \left( \ f_{\ast }\xi
\right) ,
\end{equation*}%
recall from Proposition $\left( 3.1\right) /\left( 1^{\circ }\right) $ that 
\begin{equation}
\mathcal{T}\left( f\right) =\nabla _{\left( f_{\ast }\xi \right) }\left(
f_{\ast }\xi \right) -\ f_{\ast }\varphi \left( \delta \varphi \right)
-\left( \delta \eta \right) f_{\ast }\xi +\overset{m}{\underset{i=1}{\sum }}%
S^{N}\left( \phi E_{i},E_{i}\right) .  \tag{$3.9$}
\end{equation}%
But then since $\varphi \left( \delta \varphi \right) $ and$\ e_{i}\in
\Gamma \left( D^{M}\right) ,$ by first part of the assertion, one gets $\
f_{\ast }\varphi \left( \delta \varphi \right) =0$ \ and $\ E_{i}=\ f_{\ast
}e_{i}=0$ \ and therefore 
\begin{equation*}
\overset{m}{\underset{i=1}{\sum }}S^{N}\left( \phi E_{i},E_{i}\right) =0.
\end{equation*}%
So, the equation $\left( 3.9\right) $ gives the result required.

$1^{\circ })\left( ii\right) :$ \ $f_{\ast }:TM\rightarrow TN$\ \ can be
neither injective and therefore $f$ \ can not be weakly conformal nor
surjective and therefore $f$ \ can not be horizontally weakly conformal
since $f_{\ast }\left( X\right) =0,\ \forall \ X\in \Gamma \left(
D^{M}\right) .$

$2^{\circ })$ \ For $M=\left( M^{2m+1};\ \varphi ,\ \xi ,\ \eta ,\ g\right)
\in $\textit{\ }$\mathcal{B}$\textit{\ \ and }$\ N=\left( N^{2n+1};\ \phi ,\
\gamma ,\ \sigma ,\ h\right) \in \mathcal{A}_{7}^{a}$ \ suppose \ $f$ \ is
non-constant $\left( \pm \right) $\textit{-}holomorphic horizontally weakly
conformal map with the dilation $\mu $.

The equation $\left( 3.8\right) $, that is,%
\begin{equation*}
f_{\ast }\left( L_{\varphi }\right) \text{ }\in \Gamma \left( D^{N}\right)
\end{equation*}%
is still valid by the same argument used in $\left( 1^{\circ }\right) .$
Now, choose a local orthonormal frame field \ $\left\{ e_{_{1}},\cdots
,e_{_{m}};\varphi e_{_{1}},\cdots ,\varphi e_{_{m}}\right\} $ \ for\ $D^{M}$
\ such that 
\begin{equation*}
\left\{ v_{1},\cdots ,v_{n};\ \phi v_{1},\cdots ,\phi v_{n}\right\}
\end{equation*}%
forms a local orthonormal frame field for $\ D^{N},$\ with \ $E_{i}=f_{\ast
}\left( e_{i}\right) =\mu v_{i};$ $i=1,\cdots ,n$\ . (This can be done,%
\textit{\ }by the virtue of Lemma $\left( 3.3\right) /\left( 4^{\circ
}/i\right) $ since $f$ \ is non-constant $\left( \pm \right) $\textit{-}%
holomorphic horizontally weakly conformal map). So 
\begin{eqnarray*}
f_{\ast }\left( \left( L_{\varphi }\right) _{p}\right) &=&\pm \overset{m}{%
\underset{i=1}{\sum }}\left[ \phi E_{i},\ E_{i}\right] _{f\left( p\right) }
\\
&=&\pm \overset{n}{\underset{i=1}{\sum }}\left[ \phi \mu \ v_{i},\ \mu \
v_{i}\right] _{f\left( p\right) }=\pi _{f\left( p\right) }+\omega _{f\left(
p\right) }
\end{eqnarray*}%
where \ 
\begin{equation*}
\pi _{f\left( p\right) }=\pm \mu \overset{n}{\underset{i=1}{\sum }}\left\{
\left( d\mu \left( v_{i}\right) \right) \left( \phi \ v_{i}\right) -\left(
d\mu \left( \phi v_{i}\right) \right) \left( \ v_{i}\right) \right\}
_{f\left( p\right) }
\end{equation*}%
and 
\begin{equation*}
\omega _{f\left( p\right) }=\pm \mu ^{2}\overset{n}{\underset{i=1}{\sum }}%
\left[ \ \phi v_{i},\ \ v_{i}\right] _{f\left( p\right) }.
\end{equation*}%
But note that \ $\pi _{f\left( p\right) }\in D_{f\left( p\right) }^{N}$ \
while $\ \omega _{f\left( p\right) }\notin D_{f\left( p\right) }^{N}$ $\ $%
since \ $N$ \ is \ non-\textbf{semi}-$\phi $-involutive $\ $as $N\in 
\mathcal{A}_{7}^{a}$. \ So we have \ $f_{\ast }\left( \left( L_{\varphi
}\right) _{p}\right) \notin D_{f\left( p\right) }^{N}.\ \ $But then this
contradicts with $\left( 3.8\right) ,$ so \ $f$ \ must be constant.

$3^{\circ })$\textit{\ \ }For some \ $p\in H$ \ suppose $\ f_{\ast }\left(
v\right) \neq 0$ \ for some nonezero \ $v\in T_{p}H.$ W.l.o.g. let \ $%
\left\{ e_{_{1}},\cdots ,e_{_{m}};Je_{_{1}},\cdots ,Je_{_{m}}\right\} $ \ be
a local frame field for\ $TH$ \ with \ $e_{_{1}}\left( p\right) =v.$ Since $%
M=\left( H^{2m},J\ ,G\right) $\textit{\ }is an almost Hermitian manifold we
have $f_{\ast }\left( TH\right) \subset D^{N}$ \ by the \ $\left( \pm
\right) $-holomorphicity of \ $f,$ \ so that 
\begin{equation}
f_{\ast }\left( L_{J}\right) \text{ }\in \Gamma \left( D^{N}\right) . 
\tag{$3.10$}
\end{equation}%
On the other hand, by mimicing the argument used in part $\left( 1^{\circ
}\right) ,$ we get that 
\begin{equation*}
f_{\ast }\left( L_{J}\right) \text{ }\notin \Gamma \left( D^{N}\right) .
\end{equation*}%
which contradicts with $\left( 3.10\right) .$ So we must have $f_{\ast
}\left( T_{p}H\right) =0,$ \ $\forall \ p\in H,$ that is \ $f$ \ is constant.

$4^{\circ })$\textit{\ \ }\ By the $\left( \pm \right) $-holomorphicity of \ 
$f$ \ we have that $f_{\ast }\left( T_{p}H\right) \subset D_{f_{\ast }\left(
p\right) }^{N}.$ \ This means that $f_{\ast }$ \ can not be surjective and
therefore can not be horizontally weakly conformal.

$5^{\circ })$\textit{\ \ \ }By the $\left( \pm \right) $-holomorphicity of \ 
$f$ \ we have that $f_{\ast }\left( \xi \right) =0.$ \ This means that $%
f_{\ast }$ \ can not be injective and therefore can not be weakly conformal.

\textbf{Proposition }$\left( 3.3\right) :$ \ $\left( c.f.\ \left[ 13\right] ,%
\text{Proposition }\left( 3.1\right) \right) :$ \textit{\ Let\ }%
\begin{equation*}
f:M=\left( M^{2m+1};\ \varphi ,\ \xi ,\ \eta ,\ g\right) \rightarrow
N=\left( N^{2n+1};\ \phi ,\ \gamma ,\ \sigma ,\ h\right)
\end{equation*}%
\textit{\ be a \textbf{non-constant} }$\left( \pm \right) $\textit{%
-holomorphic map between almost contact metric manifolds and }$\lambda $%
\textit{\ be as in Lemma }$\left( 3.1\right) /\left( i\right) .$\textit{\ }

$1^{\circ })$\textit{\ \ If }$\ M$\textit{\ \ and \ }$N$\textit{\ \ both
satisfy \ }$\left( GC\right) $\textit{\ \ then }$d\lambda =d\lambda \left(
\xi \right) \eta $\textit{\ }$\ $\textit{and \ therefore }$\ d\lambda \left(
X\right) =0,\ \ \forall \ X\in \Gamma \left( D^{M}\right) $\textit{, that is
\ }$\lambda $\textit{\ \ is constant along }$D^{M}$\textit{.\ }

$2^{\circ })$\textit{\ \ If }$\ M$\textit{\ \ and \ }$N$\textit{\ \ both
satisfy \ }$\left( GC\right) $\textit{\ and further \ }$M$\textit{\ \ is
also \textbf{non-semi}-}$\varphi $\textit{-involutive or \textbf{non}-}$\phi 
$-\textit{involutive\ then \ }$\lambda $\textit{\ \ is constant and that \ }$%
U_{\xi }=0.$\textit{\ Therefore, if the map }$f$\ \textit{\ is weakly
conformal with the conformal factor (resp: horizontally weakly conformal
with the dilation) }$\mu $ \textit{then it\ is homothetic (resp:
horizontally homothetic).\ }

$3^{\circ })$\textit{\ \ If \ \ }$M$\textit{\ \ is \textbf{non}-}$\varphi $%
-involutive\textit{\ or \textbf{non}-\textbf{semi}-}$\varphi $-involutive%
\textit{\ and \ }$N$\textit{\ \ is \ }$\phi $\textit{-involutive then \ }$%
\lambda =0.$\textit{\ Therefore, the map \ }$f$\ \ \textit{can neither be
weakly conformal nor \textbf{horizontally} weakly conformal.}

$4^{\circ })$\textit{\ \ If \ \ }$M$\textit{\ \ is \textbf{non}-\textbf{semi}%
-}$\varphi $-involutive\textit{\ and \ }$N$\textit{\ \ is \ }\textbf{semi-}$%
\phi $\textit{-involutive then the map} \textit{\ }$f$\ \ \textit{can not be 
\textbf{horizontally} weakly conformal.}

\textbf{Proof: }

$1^{\circ }):$\ Let \ $Y\in \Gamma \left( TM\right) .$ \ Writing \ $Y=X+r\xi 
$ \ (where \ $X$ $\in \Gamma \left( D^{M}\right) $ \ and $r$ $\in C^{\infty
}\left( M\right) ),$ we see that, \ 
\begin{equation*}
2d\eta \left( \xi ,Y\right) =2d\eta \left( \xi ,X\right) +2rd\eta \left( \xi
,\xi \right) =2d\eta \left( \xi ,X\right) =-\eta \left[ \xi ,X\right] .
\end{equation*}%
But then since \ $M$ \ satisfies \ $\left( GC\right) $ we have, by Lemma $%
\left( 2.2\right) /\left( 1^{\circ }\right) ,$ that

\begin{equation*}
\left[ \xi ,X\right] \in \Gamma \left( D^{M}\right)
\end{equation*}%
and therefore the above equation gives 
\begin{equation}
2d\eta \left( \xi ,Y\right) =-\eta \left[ \xi ,X\right] =0.  \tag{$3.11$}
\end{equation}

Now consider the pull back \ $1$-form \ $\widetilde{\eta }=f^{-1}\sigma
=\sigma f_{\ast }=\lambda \eta .$ Since $N$ satisfies \ $(GC)$, by Lemma $%
(2.2)/\left( 1^{\circ }\right) $, we get \ 
\begin{eqnarray*}
2d\widetilde{\eta }\left( \xi ,Y\right) &=&2d\sigma \left( f_{\ast }\xi ,\
f_{\ast }Y\right) =2d\sigma \left( \lambda \gamma ,Z\right) \\
&=&\lambda 2d\sigma \left( \gamma ,Z\right) =-\lambda \sigma \left[ \gamma ,Z%
\right] =0,
\end{eqnarray*}%
where \ $f_{\ast }Y=Z+s\gamma ,\ \ $with $\ Z\in \Gamma \left( D^{N}\right) $
\ and \ $s\in C^{\infty }\left( N\right) .$

On the other hand, \ since $d\eta \left( \xi ,Y\right) =0$ \ by \ $(3.11),$
we have%
\begin{eqnarray*}
0 &=&d\widetilde{\eta }\left( \xi ,Y\right) =\left( d\lambda \wedge \eta
\right) \left( \xi ,Y\right) +\lambda d\eta \left( \xi ,Y\right) \\
&=&\left( d\lambda \wedge \eta \right) \left( \xi ,Y\right) .
\end{eqnarray*}%
But then,%
\begin{equation*}
0=\left( d\lambda \wedge \eta \right) \left( \xi ,Y\right) =d\lambda \left(
\xi \right) \eta \left( Y\right) -d\lambda \left( Y\right) .
\end{equation*}%
Hence we get 
\begin{equation}
d\lambda =d\lambda \left( \xi \right) \eta  \tag{$3.12$}
\end{equation}%
and therefore 
\begin{equation}
d\lambda \left( X\right) =d\lambda \left( \xi \right) \eta \left( X\right)
=0,\text{ \ \ \ }\forall \ X\in \Gamma \left( D^{M}\right) .  \tag{$3.13$}
\end{equation}

$2^{\circ })$\textit{\ \ }Now $\left( 3.12\right) $ \ gives us \ 
\begin{equation*}
d\lambda \wedge \eta =d\lambda \left( \xi \right) \eta \wedge \eta =0
\end{equation*}%
and that 
\begin{equation*}
d\lambda \wedge d\eta =d\left( d\lambda \wedge \eta \right) =0.
\end{equation*}%
So this, together with $\left( 3.11\right) $, gives us 
\begin{eqnarray}
0 &=&2\left( d\lambda \wedge d\eta \right) \left( \xi ,X,Y\right) =2d\lambda
\left( \xi \right) d\eta \left( X,Y\right)  \notag \\
&=&-d\lambda \left( \xi \right) \eta \left( \left[ X,Y\right] \right) ,\ \ \
\forall \ X,Y\in \Gamma \left( D^{M}\right)  \TCItag{$3.14$}
\end{eqnarray}%
and that%
\begin{equation}
d\lambda \left( \xi \right) \eta \left( \overset{m}{\underset{i=1}{\sum }}%
\left[ e_{i},\varphi e_{i}\right] \right) =0.  \tag{$3.15$}
\end{equation}

Now

$\bullet \ \ $If $\ M$ \ is \textbf{non-}involutive\textbf{\ }then at any
point $p\in M$ we have $\ $%
\begin{equation*}
\eta \left( \left[ X_{\circ },Y_{\circ }\right] _{p}\right) \neq 0,\ \text{%
for some \ }\ X_{\circ },Y_{\circ }\in \Gamma \left( D^{M}\right) .
\end{equation*}%
So, from $(3.14)$ \ we get that%
\begin{equation*}
d\lambda \left( \xi \right) =0.
\end{equation*}

$\bullet $ $\ $If $\ M$ \ is \textbf{non-semi-}$\varphi $-involutive then at
any point $p\in M$ \ we have 
\begin{equation*}
\eta \left( \overset{m}{\underset{i=1}{\sum }}\left[ e_{i},\varphi e_{i}%
\right] \right) _{p}\neq 0.\ \ \ 
\end{equation*}%
So by $(3.15)$ \ we get \ 
\begin{equation*}
d\lambda \left( \xi \right) =0.
\end{equation*}%
This, together with $\left( 3.13\right) ,$ gives that $\lambda $\textit{\ \ }%
is constant.

To show that $U_{\xi }=0$ note that, since $\ M$ \ satisfies \ $(GC),$ \ we
have 
\begin{equation*}
U_{\xi }=\nabla _{\left( f_{\ast }\xi \right) }^{N}\left( f_{\ast }\xi
\right) -f_{\ast }\left( \nabla _{\xi }^{M}\xi \right) =\nabla _{\left(
f_{\ast }\xi \right) }^{N}\left( f_{\ast }\xi \right)
\end{equation*}%
But then, since \ $\lambda $ is constant and $\ N$ \ satisfies \ $(GC)$ \ we
have%
\begin{equation*}
U_{\xi }=\nabla _{\left( f_{\ast }\xi \right) }^{N}\left( f_{\ast }\xi
\right) =\nabla _{\left( \lambda \gamma \right) }^{N}\left( \lambda \gamma
\right) =\lambda ^{2}\nabla _{\gamma }^{N}\gamma =0.
\end{equation*}

Now, if $f$ \ is weakly conformal (resp: horizontally weakly conformal) then
it has to be homothetic (resp: horizontally homothetic) since, by Lemma $%
\left( 3.3\right) /\left( 2^{\circ }\text{-}ii\right) ,$ \ $\mu =\left\vert
\lambda \right\vert $ \ which is constant.

$3^{\circ })$ \ Note that, from Lemma $\left( 2.8\right) $ \ and\ $\left(
\pm \right) $-holomorphicity of $\ f,$ \ we get 
\begin{equation}
f_{\ast }\left[ X,\varphi X\right] =\pm \left[ f_{\ast }X,\ \phi f_{\ast }X%
\right] ,\ \ \ \forall \ X\in \Gamma \left( D^{M}\right)  \tag{$3.16$}
\end{equation}%
with \ $f_{\ast }X\in \Gamma \left( D^{N}\right) .$ By the assumption that $%
N $\textit{\ \ }is \ \ $\phi $-involutive\ we get:%
\begin{equation*}
\left[ f_{\ast }X,\ \phi f_{\ast }X\right] \in \Gamma \left( D^{N}\right) ,\
\ \ \forall \ X\in \Gamma \left( D^{M}\right)
\end{equation*}%
and therefore, from $\left( 3.16\right) ,$

\begin{equation}
f_{\ast }\left[ X,\ \varphi X\right] \in \Gamma \left( D^{N}\right) ,\ \ \
\forall \ X\in \Gamma \left( D^{M}\right) \   \tag{$3.17$}
\end{equation}%
Now\ 

$\bullet $ \ If $\ M$\ \ is non-$\varphi $-involutive\ then for $\ X\in
\Gamma \left( D^{M}\right) $ we have that at any given $p\in M$ \textit{\ }%
\begin{equation*}
\left[ X,\ \varphi X\right] _{p}=\left( Z_{X}\right) _{p}+r_{_{X}}\left(
p\right) \ \xi _{p},\ \ 
\end{equation*}%
for some $Z_{X}\in \Gamma \left( D^{M}\right) $ \ and smooth funtion $%
r_{_{X}}$ with $r_{_{X}}\left( p\right) \neq 0.$ But then

\begin{equation*}
f_{\ast }\left[ X,\varphi X\right] =f_{\ast }Z_{X}+r_{_{X}}\ \lambda \gamma
\ \ 
\end{equation*}%
\textit{\ }with \ \textit{\ }$f_{\ast }Z_{X}\in \Gamma \left( D^{N}\right) .$
On the other hand, this together with $\left( 3.17\right) $ gives%
\begin{equation*}
f_{\ast }Z_{X}+r_{_{X}}\lambda \gamma \in \Gamma \left( D^{N}\right) \ \ 
\end{equation*}%
from which we get $\ \left( r_{_{X}}\ \lambda \right) \left( p\right) =0$ \
and that $\lambda \left( p\right) =0.$

\bigskip $\bullet $ \ If\ $M$\ \ is non-\textit{\textbf{semi}-}$\varphi $%
-involutive. Then at any given $p\in M$ \ we may express\textit{\ }$\ $%
\begin{equation*}
L=\overset{m}{\underset{i=1}{\sum }}\left[ e_{i},\varphi e_{i}\right] =Z+r\xi
\end{equation*}%
for some $Z\in \Gamma \left( D^{M}\right) $ \ and smooth funtion $r$ \ with $%
r\left( p\right) \neq 0.$ So we have

\begin{equation}
f_{\ast }\left( L\right) =f_{\ast }Z+\lambda r\gamma  \tag{$3.18$}
\end{equation}%
with \textit{\ }$f_{\ast }Z\in \Gamma \left( D^{N}\right) .$ On the other
hand, from $\left( 3.17\right) $ \ we get%
\begin{equation*}
f_{\ast }\left[ e_{i},\varphi e_{i}\right] =\pm \left[ E_{i},\ \phi E_{i}%
\right] \in \Gamma \left( D^{N}\right)
\end{equation*}%
and therefore 
\begin{equation*}
f_{\ast }\left( L\right) =\pm \overset{m}{\underset{i=1}{\sum }}\left[
E_{i},\ \phi E_{i}\right] \in \Gamma \left( D^{N}\right)
\end{equation*}%
But then, this together with $\left( 3.18\right) $ gives%
\begin{equation*}
f_{\ast }Z+\lambda r\gamma \in \Gamma \left( D^{N}\right) \ \ 
\end{equation*}%
from which we get $\ r\lambda =0$ and therefore $\lambda =0.$ That is, \ $%
f_{\ast }\xi =0.$ So, $f$ \ can be neither injective (and therefore weakly
conformal) nor surjective (and therefore \textbf{horizontally} weakly
conformal) unlees it is constant.

$4^{\circ })$\textit{\ }\ \ Suppose $f$\ \ is a $\left( \pm \right) $%
-holomorphic \textbf{horizontally} weakly conformal map. We show that its
dilition $\mu =\left\vert \lambda \right\vert =0,$ so that it can not be
surjective and therefore it has to be constant.

For let,%
\begin{equation*}
\left\{ e_{_{1}},\cdots ,e_{_{m}};\ \varphi e_{_{1}},\cdots ,\varphi
e_{_{m}}\right\}
\end{equation*}%
be a local orthonormal frame field\textit{\ \ }for \ $D^{M}$\ \ such that \
the set

\begin{equation*}
\left\{ v_{_{1}},\cdots ,v_{n};\ \phi v_{_{1}},\cdots ,\phi v_{n}\right\}
\end{equation*}%
forms a local orthonormal frame field for \ $D^{N},$ \ where $\ f_{\ast
}e_{i}=E_{i}=\rho v_{i}$ \ for \ $i=1,\cdots ,n$ \ and $f_{\ast
}e_{i}=E_{i}=0$\ \ for\ $i=n+1,\cdots ,m.$ \ (This can be done by the virtue
of Lemma $\left( 3.3\right) /\left( \left( 3^{\circ }\right) (i)\right) ).$
Since\ $M$\ \ is \textbf{non}-\textit{\textbf{semi}-}$\varphi $-involutive,
at any given $p\in M,$ \ we may express\textit{\ }$\ $%
\begin{equation*}
L=\overset{m}{\underset{i=1}{\sum }}\left[ e_{i},\varphi e_{i}\right] =Z+r\xi
\end{equation*}%
for some $Z\in \Gamma \left( D^{M}\right) $ \ and smooth funtion $r$ \ with $%
r\left( p\right) \neq 0.$ We see that%
\begin{eqnarray*}
f_{\ast }\left( L\right) &=&\overset{m}{\underset{i=1}{\sum }}f_{\ast }\left[
e_{i},\varphi e_{i}\right] =\pm \overset{m}{\underset{i=1}{\sum }}\left[
f_{\ast }e_{i},\phi f_{\ast }e_{i}\right] \\
&=&\pm \overset{n}{\underset{i=1}{\sum }}\left[ \rho \ v_{i},\ \rho \ \phi
v_{i}\right] =\pi +\omega
\end{eqnarray*}%
where \ 
\begin{equation*}
\pi =\pm \rho \overset{n}{\underset{i=1}{\sum }}\left\{ \left( d\rho \left(
v_{i}\right) \right) \left( \phi \ v_{i}\right) -\left( d\rho \left( \phi
v_{i}\right) \right) \left( \ v_{i}\right) \right\}
\end{equation*}%
and 
\begin{equation*}
\omega =\pm \rho ^{2}\overset{n}{\underset{i=1}{\sum }}\left[ \ v_{i},\ \phi
\ v_{i}\right] .
\end{equation*}%
Clearly \ \ $\pi \in \Gamma \left( D^{N}\right) $ \ since \ $v_{i},\ \phi \
v_{i}\in \Gamma \left( D^{N}\right) .$ Also \ $\omega \in \Gamma \left(
D^{N}\right) $ \ since \ $\left\{ v_{1},\cdots ,v_{n};\ \phi v_{1},\cdots
\phi v_{n}\right\} $ \ forms a local orthonormal frame field for $\ D^{N}$ \
and \ $N$ \ is \textbf{semi}-$\phi $-involutive. Consequently we get

\begin{equation}
f_{\ast }\left( L\right) \in \Gamma \left( D^{N}\right)  \tag{$3.19$}
\end{equation}%
On the other hand, as in the proof of \ $\left( 3^{\circ }\right) ,$ at any
given $p\in M$ \ we have

\begin{equation*}
f_{\ast }\left( L\right) =f_{\ast }Z+\lambda r\gamma
\end{equation*}%
for some $f_{\ast }Z\in \Gamma \left( D^{N}\right) $ \ and smooth funtion $r$
\ with $r\left( p\right) \neq 0.$ But then, from $\left( 3.19\right) ,$ we
get $\lambda r=0$ \ which gives \ $\lambda \left( p\right) =0$ \ at any
given $p\in M.$ \ $\blacksquare $

Using \textbf{TABLE-II } and Proposition $\left( 3.3\right) $ give us the
following:

\textbf{Corollary }$\left( 3.1\right) $ \textit{Let \ \ }%
\begin{equation*}
f:M=\left( M^{2m+1};\ \varphi ,\ \xi ,\ \eta ,\ g\right) \rightarrow
N=\left( N^{2n+1};\ \phi ,\ \gamma ,\ \sigma ,\ h\right)
\end{equation*}%
\textit{\ be a \textbf{non-constant} }$\left( \pm \right) $\textit{%
-holomorphic map between almost contact metric manifolds and }$\lambda $%
\textit{\ be as in Lemma }$\left( 3.1\right) /\left( i\right) $\textit{.}

$1^{\circ })$ $\ d\lambda =d\lambda \left( \xi \right) \eta $\textit{\ }$\ $%
\textit{and \ therefore }$\ d\lambda \left( X\right) =0,\ \ \forall \ X\in
\Gamma \left( D^{M}\right) $\textit{, that is \ }$\lambda $\textit{\ \ is
constant along }$D^{M},$\textit{\ in either one of the following cases:}

$\ \bullet $\textit{\ \ }$M$\textit{\ }$\in \mathcal{A}_{4}$ \ \textit{and} $%
N$\textit{\ \ }$\in \mathcal{A}_{1}\cup \mathcal{A}_{2}\cup \mathcal{A}%
_{3}\cup \mathcal{A}_{4}\cup \mathcal{A}_{5}\cup \mathcal{A}_{6}\cup 
\mathcal{A}_{7}\cup \mathcal{B}_{1}\cup \mathcal{B}_{2}\cup \mathcal{B}_{3}$

$\ \bullet \ M\in \mathcal{B}_{1}$\textit{\ }$\cup \ \mathcal{B}_{2}$\textit{%
\ }$\cup \ \mathcal{B}_{3}$ \ \textit{and } $N$\textit{\ \ }$\in \mathcal{A}%
_{1}\cup \ \mathcal{A}_{2}\cup \mathcal{A}_{3}\cup \mathcal{A}_{4}\cup 
\mathcal{A}_{5}\cup \mathcal{A}_{6}\cup \mathcal{A}_{7}$\textit{\ }$\cup \ 
\mathcal{B}_{1}$\textit{\ }$\cup \ \mathcal{B}_{2}$\textit{\ }$\cup \ 
\mathcal{B}_{3}$

$2^{\circ })$\textit{\ \ If \ }$M$\textit{\ \ }$\in \mathcal{A}_{1}\cup \ 
\mathcal{A}_{2}\cup \mathcal{A}_{3}\cup \mathcal{A}_{5}\cup \mathcal{A}%
_{6}\cup \mathcal{A}_{7}$\textit{\ \ and\ }$\ N$\textit{\ \ }$\in \mathcal{A}%
_{1}\cup \ \mathcal{A}_{2}\cup \mathcal{A}_{3}\cup \mathcal{A}_{4}\cup 
\mathcal{A}_{5}\cup \mathcal{A}_{6}\cup \mathcal{A}_{7}\cup \ \mathcal{B}%
_{2}\ \ $\textit{then\ }$\lambda $\textit{\ is constant\ and that \ }$U_{\xi
}=0$\textit{. \ Therefore, if the map }$f$\textit{\ \ is weakly conformal (%
\textbf{resp:} \textbf{horizontally} weakly conformal) then it\ is
homothetic (\textbf{resp}: \textbf{horizontally} homothetic).\ }

$3^{\circ })$\textit{\ \ Let \ }$M$\textit{\ \ }$\in \mathcal{A}$\textit{\ \
except that it is not \textbf{almost} quasi-Sasakian. }

\textit{\ \ }$i)$\textit{\ \ \ If }$\ N$\textit{\ \ }$\in \mathcal{B}%
_{1}\cup \ \mathcal{B}_{3}$\textit{\ then the function }$\lambda $\textit{\
vanishes. \ Therefore the map \ }$f$\textit{\ \ can neither be weakly
conformal nor \textbf{horizontally} weakly conformal.}

$ii)$\textit{\ \ \ If }\ $N$\textit{\ \ }$\in \mathcal{B}_{2}$\textit{\ \
then\ the map \ }$f$\textit{\ \ can not be \textbf{horizontally} weakly
conformal.}

\textbf{Lemma }$\left( 3.4\right) :$ $\left( \left[ 2\right] ,\text{
Proposition }\left( 3.5.1\right) /\left( ii\right) \right) ;$ \textit{A
weakly conformal map from a Riemannian manifold of dimension not equal to 2
is harmonic if and only if it is homothetic and its image is minimal.}

\textbf{Lemma }$\left( 3.5\right) :$ $\left( \left[ 2\right] ,\text{%
Corollary }\left( 4.5.5\right) \right) ;$ \textit{For a h\textbf{orizontally}
weakly conformal map }$F$\textit{\ \ from a Riemannian manifold into a
Riemannian manifold of dimension greater than 2, any two of the following
conditions imply the third:}

$\ i)$ \ $F$\textit{\ is harmonic, }

$ii)$ \ $F$\textit{\ is \textbf{horizontally} homothetic, }

$iii)$ \ \textit{the fibres of }$F$\textit{\ are minimal.}

\textbf{4) \ Main Results: }Let $f:M\rightarrow N$ denote a $\left( \pm
\right) $\textbf{-holomorphic map }from an almost contact metric manifold $\
M=\left( M^{2m+1};\ \varphi ,\ \xi ,\ \eta ,\ g\right) $ \ or an almost
Hermitian manifold \ $M=\left( H^{2m},\ J,\ G\right) $ \ into an almost
contact metric manifold \ $N=\left( N^{2n+1};\ \phi ,\ \gamma ,\ \sigma ,\
h\right) $ \ or an almost Hermitian manifold \ $N=\left( H_{1}^{2n},\
J_{1},\ G_{1}\right) $\ throughout this section. For such an \ $f$ , \ we
have the following series of results:

\textbf{Theorem} $\left( 4.1\right) :$ \ ($c.f$ $\ \left[ 13\right] $
Theorem $\left( 3.3\right) ),$ \ 

\textit{Suppose\ that}

$\bullet $\textit{\ \ }$M$\textit{\ \ is a non-\textbf{semi}-}$\varphi $%
\textit{-involutive (or non-}$\varphi $\textit{-\ involutive) and \textbf{%
semi}-symplectic manifold satisfying \ }$(GC),$

$\bullet $\textit{\ \ }$N$\textit{\ \ is a quasi-symplectic manifold
satisfying \ }$(GC).$

\textit{Then \ }$f$\textit{\ \ is harmonic. }

\textbf{Proof : }\ Since \ $f$ \ is \ $\left( \pm \right) $-holomorphic,
proposition $\left( 3.1\right) $ gives us $\ $%
\begin{equation*}
\mathcal{T}\left( f\right) =U_{\xi }+\overset{m}{\underset{i=1}{\sum }}%
\left\{ S^{N}\left( \phi E_{i},\ E_{i}\right) -\ f_{\ast }S^{M}\left(
\varphi e_{i},e_{i}\right) \right\} .
\end{equation*}%
Since \ $M$ \ and \ $N$ \ both satisfy \ $\left( GC\right) $ and $M$ \ is
non-semi-$\varphi $-involutive\ by Proposition\textbf{\ }$\left( 3.3\right)
/\left( 2^{\circ }\right) ,$ we have \ $U_{\xi }=0.$ Also,
semi-symplecticity of \ $M$ \ and quasi-symplecticity of \ $N$\ \ give that%
\begin{equation*}
\underset{i=1}{\overset{m}{\text{ }\sum }}S^{M}\left( \varphi
e_{i},e_{i}\right) =0\ \ \ \text{and\ \ \ }S^{N}\left( \phi E_{i},\
E_{i}\right) =0
\end{equation*}
\ respectively. That means that \ $f$ \ is harmonic. \ \ 

\textbf{Remark }$\left( 4.1\right) :$ \ In the above Theorem, if we impose
on \ $M$\ \ more restrictive condition "\textit{quasi-symplectic}" rather
than "\textit{semi-symplectic}" then \ $f$\ \ also becomes $D$-pluriharmonic.

\textbf{Corollary }$\left( 4.1\right) :$

\textit{\ }$\mathbf{1}^{\circ })$ \ 

$\ \ \ i)$\textit{\ \ If \ }$M\in A_{1}$\textit{\ \ and }$N\in \mathcal{A}%
_{1}\mathcal{\cup A}_{4}\mathcal{\cup \ B}_{1}\mathcal{\cup \ C}_{1}$\ 
\textit{\ then \ }$f$\textit{\ \ is harmonic and also }$D$\textit{%
-pluriharmonic.}

$\ ii)$\textit{\ \ If \ }$M\in \mathcal{A}_{6}\mathcal{\cup A}_{7}$\textit{\
and }$N\in \mathcal{A}_{1}\mathcal{\cup A}_{4}\mathcal{\cup \ B}_{1}\mathcal{%
\cup \ C}_{1}$\textit{\ then \ }$f$\textit{\ \ is harmonic.}

$\mathbf{2}^{\circ })$ \ \textit{For \ \ }$M\in \mathcal{A}_{1}\mathcal{\cup
A}_{6}\mathcal{\cup A}_{7}$ \ \textit{If }$\ N\in \mathcal{A}_{1}\mathcal{%
\cup A}_{4}$\textit{\ \ then every }$\ $\textit{weakly conformal map }$f$%
\textit{\ \ is homothetic harmonic with minimal image.}

\textbf{Proof : }\ 

$\mathbf{1}^{\circ })/\left( i\right) ,\ \left( ii\right) :$ \ As it can be
seen from the \textbf{TABLE-II} that \ $M$ \ and \ $N$ \ satisfy the
hypothesis of Theorem $\left( 4.1\right) $ and Remark\textbf{\ }$\left(
4.1\right) .$ So the results follow.

$\mathbf{2}^{\circ }):$\textit{\ \ }By part $\left( \mathbf{1}^{\circ
}\right) ,$ we see that \ $f$\ \ is harmonic. But then, from Lemma $\left(
3.4\right) ,$ the result follows

\textbf{Remark }$\left( 4.2\right) :$ \ 

$\ i)$\ \ In $\left( \left[ 23\right] ,\text{ Theorem }5\right) $ \ states
that

$\bullet $ \ Any $\left( \pm \right) $-holomorphic submersion $\ f$\ \
between quasi-contact metric manifolds\ is harmonic.

However, the condition on $\ f$\ \ that being a submersion can be dropped as
one can deduce from Corollary\ $\left( 4.1\right) /\left( \mathbf{1}^{\circ }%
\text{-}i\right) $ in particular that

$\bullet $ \ Any $\left( \pm \right) $-holomorphic map between quasi-contact
metric manifolds is harmonic.

$ii)$ Corollary\ $\left( 4.1\right) /\left( \mathbf{1}^{\circ }/i\right) $
recovers the following results, as particular cases:

$\bullet $ $\ \left( \left[ 17\right] ,\text{ Theorem }\left( 2.1\right)
\right) :$ \ Any holomorphic map between two strongly pseudoconvex $CR$%
-manifolds is harmonic.

$\bullet $ $\ \left( \left[ 17\right] ,\text{ Proposition }\left( 2.2\right)
\right) :$\ Let \ \ $f:M=\left( M^{2m+1};\ \varphi ,\ \xi ,\ \eta ,\
g\right) \rightarrow \left( H^{2n};\ J,\ G\right) $\ \ be a \ $\left(
\varphi ,J\right) $-holomorphic map from a strongly pseudoconvex $CR$%
-manifold into a Kaehler one. Then \ $f$\ \ is harmonic.

\textbf{Theorem }$\left( 4.2\right) :$ \ 

$1^{\circ })$\textit{\ \ \ If \ }$M\in \mathcal{A}_{1}\cup \mathcal{A}%
_{2}\cup \mathcal{A}_{3}\cup \mathcal{A}_{6}$\textit{\ \ and \ \ }$N\in 
\mathcal{B}_{3}$ \ \textit{then \ }$f$\textit{\ \ is harmonic if and only if
\ it is constant.}$\ $

$2^{\circ })$\textit{\ \ For \ }$M\in \mathcal{A}_{1}\cup \mathcal{A}%
_{6}\cup \mathcal{A}_{7}$\textit{\ \ and \ \ }$N\in \mathcal{A}_{3}$\textit{%
\ }

$i)$\ \ \textit{\ }$f$\textit{\ \ is harmonic if and only if \ it is
constant along \ }$D^{M}$ ( that is, $f_{\ast }X=0,\ \ \forall \ X\in \Gamma
\left( D^{M}\right) ).$

$ii)$\textit{\ \ Let }$\mathit{\ }f$\ \ \textit{further be weakly (or 
\textbf{horizontally} weakly) conformal map from }$M$\textit{\ into }$N.$ \
Then \textit{\ \ }$f$\textit{\ \ is harmonic if and only if \ it is constant.%
}

\textbf{Proof :}

$1^{\circ }):$ \ \ As it can be seen from \textbf{TABLE-II }that

$\bullet $ \ $M$ \ satisfies \ $\left( GC\right) ,$

$\bullet $ \ $M$ \ is non-semi-involutive

$\bullet $ $\ M\ \ $is semi-symplectic when $\ M\in \mathcal{A}_{1}\cup 
\mathcal{A}_{6}$

$\bullet $ \ $N$ \ is $\ \varphi $-involutive. \ \ \ \ \ \ \ \ \ \ \ \ \ \ \
\ \ \ \ \ \ \ \ \ \ \ \ \ \ \ \ \ \ \ \ \ \ \ \ \ \ \ \ \ \ \ \ \ \ \ \ \ \
\ \ \ \ \ \ \ \ \ \ \ \ \ \ \ \ \ \ \ \ \ \ \ \ \ \ \ \ \ \ \ \ \ \ \ \ \ \
\ \ \ \ \ \ \ \ \ \ \ \ \ \ \ \ \ \ \ \ \ \ \ \ \ \ \ \ \ \ \ \ \ \ \ \ \ \
\ \ \ \ \ \ \ \ \ \ \ \ \ \ \ \ \ \ \ \ \ \ \ \ \ \ \ \ \ \ \ \ \ \ \ \ \ \
\ \ \ \ \ \ \ \ \ \ \ \ \ \ \ \ \ \ \ \ \ \ \ \ \ \ \ \ \ \ \ \ \ \ \ \ \ \
\ \ \ \ \ \ \ \ \ \ \ \ \ \ \ \ \ \ \ \ \ \ \ \ \ \ \ \ \ \ \ \ \ \ \ \ \ \
\ \ \ \ \ \ \ \ \ \ \ \ \ \ \ \ \ \ \ \ \ \ \ \ \ \ \ \ \ \ \ \ \ \ \ \ \ \
\ \ \ \ \ \ \ \ So, by Proposition $\left( 3.3\right) /\left( 3^{\circ
}\right) $ \ we get \ $\lambda =0$ \ (and \ that $\ \ f_{\ast }\xi =0)$ $\ $%
and $\ U_{\xi }=0.$

On the other hand, 
\begin{equation*}
\overset{m}{\underset{i=1}{\sum }}S^{M}\left( \varphi e_{i},e_{i}\right)
=\left\{ 
\begin{tabular}{cc}
$0$, & $M\in \mathcal{A}_{1}\cup \mathcal{A}_{6}$ \\ 
$\left( \delta \eta \right) \xi ,$ & $M\in \mathcal{A}_{2}\cup \mathcal{A}%
_{3}$%
\end{tabular}%
\right\}
\end{equation*}%
and hence%
\begin{equation*}
\ f_{\ast }\overset{m}{\underset{i=1}{\sum }}S^{M}\left( \varphi
e_{i},e_{i}\right) =\left\{ 
\begin{tabular}{cc}
$0$, & $M\in \mathcal{A}_{1}\cup \mathcal{A}_{6}$ \\ 
$\left( \delta \eta \right) f_{\ast }\xi =0,$ & $M\in \mathcal{A}_{2}\cup 
\mathcal{A}_{3}$%
\end{tabular}%
\right\} .
\end{equation*}%
From this, together with the fact that \ $U_{\xi }=0,$ \ Proposition $\left(
3.1\right) $ gives that

\begin{equation*}
\mathcal{T}\left( f\right) =\overset{m}{\underset{i=1}{\sum }}S^{N}\left(
\phi E_{i},\ E_{i}\right) .
\end{equation*}%
But then, since \ $N\in \mathcal{B}_{3},$ \ we have 
\begin{equation*}
S^{N}\left( \phi E_{i},\ E_{i}\right) =2\beta h\left( E_{i,\ }E_{i}\right)
\gamma \ \ \ 
\end{equation*}%
and therefore%
\begin{equation*}
\mathcal{T}\left( f\right) =\overset{m}{\underset{i=1}{\sum }}S^{N}\left(
\phi E_{i},\ E_{i}\right) =2\beta \overset{m}{\underset{i=1}{\sum }}h\left(
E_{i,\ }E_{i}\right) \gamma ,\ \ \beta \in \mathcal{C}^{\infty }\left(
N\right) \text{ \ with \ }\beta \left( q\right) \neq 0,\text{ \ \ }\forall \
q\in N\ \ 
\end{equation*}%
So, if $\ f$ \ is harmonic if and only if \ 
\begin{equation*}
h\left( E_{i,\ }E_{i}\right) =0,\ \ \ \ \ \forall \ i=1,\cdots ,m.
\end{equation*}%
That is, \ 
\begin{equation*}
E_{i}=f_{\ast }e_{i}=0,\ \ \ \ \forall \ i=1,\cdots ,m.
\end{equation*}%
From this the result follows. \ 

$2^{\circ })/\left( i\right) $ \ $:$ \ \ As it can be seen from \textbf{%
TABLE-II \ }that

$\bullet $ \ $M$ \ satisfies \ $\left( GC\right) ,$

$\bullet $ \ $M$ \ is none-semi-involutive

$\bullet $ \ $N$ \ satisfies \ $\left( GC\right) .$ \ \ \ \ \ \ \ \ \ \ \ \
\ \ \ \ \ \ \ \ \ \ \ \ \ \ \ \ \ \ \ \ \ \ \ \ \ \ \ \ \ \ \ \ \ \ \ \ \ \
\ \ \ \ \ \ \ \ \ \ \ \ \ \ \ \ \ \ \ \ \ \ \ \ \ \ \ \ \ \ \ \ \ \ \ \ \ \
\ \ \ \ \ \ \ \ \ \ \ \ \ \ \ \ \ \ \ \ \ \ \ \ \ \ \ \ \ \ \ \ \ \ \ \ \ \
\ \ \ \ \ \ \ \ \ \ \ \ \ \ \ \ \ \ \ \ \ \ \ \ \ \ \ \ \ \ \ \ \ \ \ \ \ \
\ \ \ \ \ \ \ \ \ \ \ \ \ \ \ \ \ \ \ \ \ \ \ \ \ \ \ \ \ \ \ \ \ \ \ \ \ \
\ \ \ \ \ \ \ \ \ \ \ \ \ \ \ \ \ \ \ \ \ \ \ \ \ \ \ \ \ \ \ \ \ \ \ \ \ \
\ \ \ \ \ \ \ \ \ \ \ \ \ \ \ \ \ \ \ \ \ \ \ \ \ \ \ \ \ \ \ \ \ \ \ \ \ \
\ \ \ \ \ \ \ \ \ So, by Proposition $\left( 3.3\right) /\left( 2^{\circ
}\right) $ \ we have \ $\mathcal{U}_{\xi }=0.$ \ Therefore, from Proposition 
$\left( 3.1\right) $ we get%
\begin{equation*}
\mathcal{T}\left( f\right) =\overset{m}{\underset{i=1}{\sum }}\left\{
S^{N}\left( \phi E_{i},\ E_{i}\right) -f_{\ast }S^{M}\left( \varphi e_{i},\
e_{i}\right) \right\} .
\end{equation*}%
map\textit{\ }$f.$\textit{\ \ }But then, \ since \ $M\in \mathcal{A}_{1}\cup 
\mathcal{A}_{6}\cup \mathcal{A}_{7}$\textit{\ }\ is semi-symplectic and%
\textit{\ \ \ }$N\in \mathcal{A}_{3}$\textit{\ }the above equation becomes

\begin{equation*}
\mathcal{T}\left( f\right) =\overset{m}{\underset{i=1}{\sum }}S^{N}\left(
\phi E_{i},\ E_{i}\right) =\frac{1}{n}\delta \sigma \overset{m}{\underset{i=1%
}{\sum }}h\left( E_{i\ },E_{i}\right) \gamma ,\ \ \ 
\end{equation*}%
where $\ $\ $\left( \delta \sigma \right) \left( q\right) \neq 0,$ \ \ $%
\forall \ q\in N$ $\ $by the assumption. From this, the result follows.

$2^{\circ })/\left( ii\right) $ \ $:$ \ \ Observe that when \ $f$\ \ is a
weakly conformal or horizontally weakly conformal, the function $\lambda $ \
becomes the conformal factor of \ $f$. \ So we have%
\begin{equation*}
h\left( E_{i\ },E_{i}\right) =\lambda ^{2}g\left( e_{i\ },e_{i}\right)
=\lambda ^{2}.
\end{equation*}%
Using this in the last equation in part \ $\left( 2^{\circ }\right) /\left(
i\right) $ \ we get \ \ \ \ \ \ \ \ \ \ \ \ \ \ \ \ \ \ \ \ \ \ \ \ \ \ \ \
\ \ \ \ \ \ \ \ \ \ \ \ \ \ \ \ \ \ \ \ \ \ \ \ \ \ \ \ \ \ \ \ \ \ \ \ \ \
\ \ \ \ \ \ \ \ \ \ \ \ \ \ \ \ \ \ \ \ \ \ \ \ \ \ \ \ \ 
\begin{equation*}
\mathcal{T}\left( f\right) =\frac{1}{n}\delta \sigma \overset{m}{\underset{%
i=1}{\sum }}h\left( E_{i\ },E_{i}\right) \gamma =\frac{m}{n}\lambda
^{2}\delta \sigma .
\end{equation*}
\ \ \ \ \ From this we deduce the following:

\textit{\ }$f$\textit{\ \ }is harmonic if and only if \textit{\ }$\lambda $
vanishes and therefore $f$\textit{\ \ }is constant.

\textbf{Remark }$\left( 4.3\right) :$ In $\left[ 29\right] ,$ Theorem $%
\left( 4.2\right) $ states that

\textit{Let }$\ f:M\rightarrow N$\textit{\ \ be a }$\left( \varphi ,\phi
\right) $\textit{-holomorphic map from a contact metric manifold into an
almost Kenmotsu manifold. Then }$f$ \textit{\ is harmonic if and only if it
is constant.}

In our work, Theorem $\left( 4.2\right) $ generalizes this result by
allowing the domain $M$ \ to be quasi-contact metric, quasi-$\mathcal{K}$%
-Sasakian, nearly-$\alpha $-Sasakian and trans-Sasakian as well as contact
metric manifold.

\textbf{Theorem }$\left( 4.3\right) :$ \ 

$1^{\circ })$\textit{\ \ Let }$M\in \mathcal{A}_{2}\cup \mathcal{A}_{3}$%
\textit{\ and \ \ }$N\in \mathcal{B}_{1}\cup \mathcal{C}_{1}.$\textit{Then\
\ }$f$\textit{\ \ is harmonic and }$D$\textit{-pluriharmonic.}$\ $ \textit{\
\ \ \ \ \ \ \ \ \ \ \ \ \ \ \ \ \ \ \ \ \ \ \ \ \ \ \ \ \ \ \ \ \ \ \ \ \ \
\ \ \ \ \ \ \ \ \ \ \ \ \ \ \ \ \ \ \ \ \ }

$2^{\circ })$\textit{\ \ \ Let \ }$M\in \mathcal{A}_{2}\cup \mathcal{A}_{3}$%
\textit{\ \ and \ }$N\in \mathcal{A}_{1}\cup \mathcal{A}_{4}.$\textit{\ \
Then\ }

\textit{\ }$\ i)$\textit{\ \ \ \ }$f$\textit{\ \ is harmonic if and only if }%
$\ f$\textit{\ \ is constant along \ }$\xi $, that is, \ $\lambda =0.$%
\textit{\ }

\textit{\ }$ii)$\textit{\ \ \ every }$\left( \pm \right) $\textit{holomorphic%
} \textit{weakly conformal or horizontally weakly conformal map\ is harmonic
if and only if \ it is constant}

$3^{\circ })$\textit{\ \ \ Let \ }$M\in \mathcal{A}_{2}\cup \mathcal{A}_{3}$%
\textit{\ \ and \ }$N\in \mathcal{A}_{3}.$

\textit{\ }$\ i)$\textit{\ \ \ Then }$\ f$\textit{\ \ is harmonic if and
only if \ }%
\begin{equation*}
\ \delta \sigma \overset{m}{\underset{i=1}{\sum }}h\left( E_{i},\
E_{i}\right) =n\lambda \left( \delta \eta \right) .
\end{equation*}

\textit{\ In particular, if }$M$\textit{\ is trans-Sasakian of type \ }$%
\left( \wp ,\theta \right) $\textit{\ then \ the }$\left( \pm \right) $%
\textit{holomorphic} \textit{\ }$f:M\rightarrow M$\textit{\ \ is harmonic if
and only if }%
\begin{equation*}
\overset{m}{\underset{i=1}{\sum }}h\left( E_{i},\ E_{i}\right) =\lambda m
\end{equation*}

$\ ii)$\textit{\ \ \ \ every }$\left( \pm \right) $holomorphic \textit{%
weakly conformal (resp: horizontally weakly conformal) map\ is harmonic if
and only if \ either it is constant or }%
\begin{equation*}
\lambda =\frac{n\delta \eta }{m\delta \sigma }\ \ \ \ \left( res:\lambda =%
\frac{\delta \eta }{\delta \sigma }\right)
\end{equation*}

\textit{\ In particular, if \ }$M$\textit{\ \ is trans-Sasakian of type \ }$%
\left( \wp ,\theta \right) $\textit{\ then every }$\left( \pm \right) $%
\textit{holomorphic} \textit{weakly conformal or horizontally weakly
conformal map\ \ }$f:M\rightarrow M$\textit{\ \ is harmonic if and only if
it is \ either an isometric immersion or constant.}

$4^{\circ })$\textit{\ \ Let \ }$M\in \mathcal{A}_{2}\cup \mathcal{A}_{3}$%
\textit{\ \ and \ }$N\in \mathcal{A}_{2}$

\textit{\ }$\ i)$\textit{\ \ \ Then \ }$f$ \ \textit{is harmonic if and only
if }%
\begin{equation*}
\sigma \left( \overset{m}{\underset{i=1}{\sum }}W_{E_{i}}\right) =\lambda
\delta \eta
\end{equation*}

\bigskip $\ ii)$\textit{\ \ \ \ every }$\left( \pm \right) $\textit{%
holomorphic \textbf{horizontally} weakly conformal map\ is a harmonic
morphism if and only if \ either it is constant or }%
\begin{equation*}
\lambda =\frac{\delta \eta }{\delta \sigma }
\end{equation*}

\textit{In particular, if }$M$\textit{\ \ is \textbf{nearly}-trans-Sasakian
of type \ }$\left( \wp ,\theta \right) $\textit{\ then every }$\left( \pm
\right) $\textit{holomorphic \textbf{horizontally} weakly conformal map\ \ }$%
f:M\rightarrow M$\textit{\ \ is harmonic if and only if it is \ either an
isometric immersion or constant.}

\textit{\ \ }\textbf{Proof : }Since \textbf{\ }$M\in \mathcal{A}_{2}\cup 
\mathcal{A}_{3},$ as it can be seen from the table that,\ \ $M$ \ is non-$%
\varphi $-involutive and it satisfies \ $\left( GC\right) .$ Also recall
that, by definition, \ $\wp \left( p\right) \neq 0$ \ and \ $\theta \left(
p\right) \neq 0,$ \ $\forall \ p\in M.$

$1^{\circ })$ \ : \ \ $N$\ \ is $\ \varphi $-involutive and quasi-symplectic
since $\ N$ \ $\in \mathcal{B}_{1}\cup \mathcal{C}_{1}$. So by Proposition $%
\left( 3.3\right) /\left( ii\right) ,$ \ we get \ $\lambda =0$ \ so that \ $%
\ f_{\ast }\xi =0$ \ and that $\ U_{\xi }=0$ \ since \ $M$ \ satisfies \ $%
\left( GC\right) .$ On the other hand, since \ $M\in \mathcal{A}_{2}\cup 
\mathcal{A}_{3},$%
\begin{equation*}
S^{M}\left( \varphi X,X\right) =\left\{ 
\begin{tabular}{cc}
$\eta \left( W_{X}\right) \xi $ & $M\in \mathcal{A}_{2}$ \\ 
$-2\theta g\left( X,X\right) \xi ,$ & $M\in \mathcal{A}_{3}$%
\end{tabular}%
\right\}
\end{equation*}%
and hence%
\begin{equation*}
\ f_{\ast }S^{M}\left( \varphi X,X\right) =\left\{ 
\begin{tabular}{cc}
$\eta \left( W_{X}\right) f_{\ast }\xi =0,$ & $M\in \mathcal{A}_{2}$ \\ 
$-2\theta g\left( X,X\right) f_{\ast }\xi =0,$ & $M\in \mathcal{A}_{3}$%
\end{tabular}%
\right\}
\end{equation*}%
$\forall \ X\in \Gamma \left( D^{M}\right) $. \ From this and by Lemma $%
\left( 3.2\right) /\left( 1^{\circ }-ii\right) ,$ we get

\begin{equation*}
U_{X}=S^{N}\left( \phi \ f_{\ast }X,\ \ f_{\ast }X\right) ,\ \ \ \forall \
X\in \Gamma \left( D^{M}\right) .
\end{equation*}
(Note here that \ $\ f_{\ast }X\in \Gamma \left( D^{N}\right) ,$ \ \ $\
\forall \ X\in \Gamma \left( D^{M}\right) ).$ But then, since \ $N$ $\ $is
quasi-symplectic we get%
\begin{equation*}
U_{X}=S^{N}\left( \phi \ f_{\ast }X,\ \ f_{\ast }X\right) =0,\ \ \ \forall \
X\in \Gamma \left( D^{M}\right)
\end{equation*}%
This gives that \ $f$ \ is $D$-pluriharmonic. Harmonicity follows from the
fact that \ $U_{\xi }=0$.

$2^{\circ })/(i)$ \ \ : \ Note here that by Proposition $\left( 3.3\right)
/\left( 2^{\circ }\right) $ \ we have that \ $\lambda $ \ is constant and
that \ $U_{\xi }=0$. Also noting that \ $N$ \ is a quasi-symplectic, from
Proposition $\left( 3.1\right) /\left( 1^{\circ }\right) $, we get 
\begin{equation*}
\mathcal{T}\left( f\right) =-f_{\ast }\left( \overset{m}{\underset{i=1}{\sum 
}}S^{M}\left( \varphi e_{i},\ e_{i}\right) \right) =-\left( \delta \eta
\right) f_{\ast }\xi
\end{equation*}%
with $\left( \delta \eta \right) \left( p\right) \neq 0,\ \ \forall p\in M.$
\ Consequently ,$f$\ \ is harmonic if and only if $f_{\ast }\xi =0.$

$2^{\circ })/(ii)$ \ \ : Note here that the function\ \ $\lambda $ is the
conformal factor (resp: dilation) of the weakly conformal (resp: \textbf{%
horizontally} weakly conformal) map $f$. \ But then, by Part $(2^{\circ
})/(i),$ \ we have \ "$\lambda =0$ and therefore $\ f$ \ is constant\ if and
only if $f$\ \ is harmonic."

$3^{\circ })/(i)$ \ \ : \ Since $\ N\in \mathcal{A}_{3},$ it satisfies $%
\left( GC\right) .$ So, Proposition $\left( 3.3\right) /\left( 2^{\circ
}\right) $ \ gives that \ $\lambda $ \ is constant and that \ $U_{\xi }=0$.
Therefore, from Proposition $\left( 3.1\right) /\left( 1^{\circ }\right) $,
we get%
\begin{equation*}
\mathcal{T}\left( f\right) =\overset{m}{\underset{i=1}{\sum }}S^{N}\left(
\phi E_{i},\ E_{i}\right) -\overset{m}{\underset{i=1}{\sum }}\ f_{\ast
}S^{M}\left( \varphi e_{i},e_{i}\right)
\end{equation*}%
But then, from the TABLE-II, we see that 
\begin{equation*}
S^{M}\left( \varphi e_{i},e_{i}\right) =\left( \delta \eta \right) \xi
\end{equation*}%
and%
\begin{equation*}
S^{N}\left( \phi E_{i},\ E_{i}\right) =\frac{1}{n}\left( \delta \sigma
\right) h\left( E_{i},\ E_{i}\right) \gamma .
\end{equation*}%
Thus the above harmonicity equation becomes%
\begin{eqnarray}
\mathcal{T}\left( f\right) &=&\left( \frac{1}{n}\left( \delta \sigma \right) 
\overset{m}{\underset{i=1}{\sum }}h\left( E_{i},\ E_{i}\right) \right)
\gamma -f_{\ast }\left( \delta \eta \right) \xi  \notag \\
&=&\left( \frac{1}{n}\left( \delta \sigma \right) \overset{m}{\underset{i=1}{%
\sum }}h\left( E_{i},\ E_{i}\right) -\left( \delta \eta \right) \lambda
\right) \gamma  \TCItag{$4.1$}
\end{eqnarray}%
from which the result follows.

In particular, \textit{if }$M$\textit{\ is a trans-Sasakian manifold of type
\ }$\left( \wp ,\theta \right) $\textit{\ then for the }$\left( \pm \right) $%
\textit{holomorphic} \textit{\ }$f:M\rightarrow M$\textit{\ the above \
equation simplifies to }%
\begin{equation*}
\mathcal{T}\left( f\right) =\left( \delta \eta \right) \left( \frac{1}{m}%
\overset{m}{\underset{i=1}{\sum }}h\left( E_{i},\ E_{i}\right) -\lambda
\right) \gamma .
\end{equation*}%
Thus the result follows.

$3^{\circ })/(ii)$ \ \ : \ As the function\ \ $\lambda $ is the conformal
factor (resp: dilation) of the weakly conformal (resp: \textbf{horizontally }%
weakly conformal) map $f,$ we have, for $f,e_{i}=E_{i}\neq 0$ 
\begin{equation*}
h\left( E_{i},\ E_{i}\right) =\lambda ^{2}g\left( e_{i},\ e_{i}\right)
=\lambda ^{2}
\end{equation*}%
so that the equation $\left( 4.1\right) $ becomes%
\begin{equation}
\mathcal{T}\left( f\right) =\lambda \left\{ \frac{m}{n}\left( \delta \sigma
\right) \lambda -\left( \delta \eta \right) \right\} \gamma  \tag{$4.2$}
\end{equation}%
when $f$ \ is weakly conformal and 
\begin{equation}
\mathcal{T}\left( f\right) =\lambda \left\{ \left( \delta \sigma \right)
\lambda -\left( \delta \eta \right) \right\} \gamma  \tag{$4.3$}
\end{equation}%
when $f$ \ is \textbf{horizontally }weakly conformal.

In particular, \textit{if }$M$\textit{\ is a trans-Sasakian manifold of type
\ }$\left( \wp ,\theta \right) $\textit{\ then for the }weakly conformal or 
\textbf{horizontally }weakly conformal map \textit{\ }$f:M\rightarrow M$%
\textit{\ \ we have}%
\begin{equation*}
\mathcal{T}\left( f\right) =\left( \delta \sigma \right) \lambda \left(
\lambda -1\right) \gamma .
\end{equation*}%
So, this and Equations $\left( 4.2\right) ,\left( 4.3\right) \ $give the
required results.

$4^{\circ })/(i)$ \ \ : \ By the same argument used in the proof of $\left(
3^{\circ }\right) )/(i),$ the harmonicity equation becomes 
\begin{equation*}
\mathcal{T}\left( f\right) =\overset{m}{\underset{i=1}{\sum }}S^{N}\left(
\phi E_{i},\ E_{i}\right) -\overset{m}{\underset{i=1}{\sum }}\ f_{\ast
}S^{M}\left( \varphi e_{i},e_{i}\right) .
\end{equation*}%
But then. from the TABLE-II. we see that%
\begin{equation*}
S^{M}\left( \varphi e_{i},e_{i}\right) =\left( \delta \eta \right) \xi
\end{equation*}%
and%
\begin{equation*}
S^{N}\left( \phi E_{i},\ E_{i}\right) =\sigma \left( W_{E_{i}}\right) \gamma
.
\end{equation*}%
Thus the above harmonicity equation becomes%
\begin{equation}
\mathcal{T}\left( f\right) =\left\{ \overset{m}{\underset{i=1}{\sum }}\sigma
\left( W_{E_{i}}\right) -\left( \delta \eta \right) \lambda \right\} \gamma 
\tag{$4.4$}
\end{equation}%
from which the result follows.

$4^{\circ })/(ii)$ \ \ : \ Observe here that the function\ \ $\lambda $ is
also the dilation of the \textbf{horizontally} weakly conformal map $f$ \
and that the set of sections%
\begin{equation*}
\left\{ \gamma ,\ v_{1},\cdots ,v_{n};\ \phi v_{1},\cdots ,\phi v_{n}\right\}
\end{equation*}%
becomes a local orthonormal frame field for $\Gamma \left( TN\right) ,$
where $E_{i}=f_{\ast }e_{i}=\lambda v_{i}.$ Now using the fact that 
\begin{equation*}
\sigma \left( W_{E_{i}}\right) =\lambda ^{2}\sigma \left( W_{v_{i}}\right)
\end{equation*}%
we get%
\begin{equation*}
\overset{m}{\underset{i=1}{\sum }}\sigma \left( W_{E_{i}}\right) =\lambda
^{2}\overset{n}{\underset{i=1}{\sum }}\sigma \left( W_{v_{i}}\right)
=\lambda ^{2}\left( \delta \sigma \right) .
\end{equation*}%
So, using this in $\left( 4.4\right) $ \ we get the equation $\left(
4.3\right) $, namely:%
\begin{equation*}
\mathcal{T}\left( f\right) =\lambda \left\{ \left( \delta \sigma \right)
\lambda -\left( \delta \eta \right) \right\} \gamma
\end{equation*}%
from which the result follows.

\textbf{Remark }$\left( 4.4\right) :$ Theorem\ $\left( 4.3\right) /\left(
1^{\circ }\right) $ recovers the results obtained\ in $\left( \left[ 16%
\right] \right) :$

$i)$\ $\ \ \ \left( \left[ 16\right] ,\text{ Theorem }2\text{ }\right) :$\ 
\textit{Let}%
\begin{equation*}
f:M=\left( M^{2m+1};\ \varphi ,\ \xi ,\ \eta ,\ g\right) \rightarrow \left(
H^{2n};\ J,\ G\right)
\end{equation*}%
\textit{\ \ be a \ }$\left( \varphi ,J\right) $\textit{-holomorphic map from
a nearly-trans-Sasakian manifold into a \textbf{quasi}-Kaehler one. Then \ }$%
f$\textit{\ \ is \ harmonic}

$ii)$ $\ \left( \left[ 16\right] ,\text{ Proposition }1\text{ }\right) :$ 
\textit{Let \ \ }%
\begin{equation*}
f:M=\left( M^{2m+1};\ \varphi ,\ \xi ,\ \eta ,\ g\right) \rightarrow \left(
H^{2n};\ J,\ G\right)
\end{equation*}%
\textit{\ \ be a \ }$\left( \varphi ,J\right) $\textit{-holomorphic map from
a nearly-trans-Sasakian manifold into a Kaehler one. Then \ }$f$\textit{\ \
is }$D$\textit{-pluriharmonic.}

Note here that \ Theorem\textbf{\ }$\left( 4.3\right) $ not only recovers
result in $\left( \left[ 16\right] ,\text{ Proposition }1\text{ }\right) $
but also improves it too, by allowing the target manifold $N$ to be \textbf{%
quasi}-Kaehler\textit{\ }(as well as Kaehler).

\textbf{Theorem }$\left( 4.4\right) :$\textit{\ \ Let \ }$M\in \mathcal{A}%
_{7}^{a},$\textit{\ that is, \ }$M$\textit{\ \ is an almost semi-Sasakian
manifold. }

$1^{\circ })$\textit{\ }$\ \ $\textit{If}$\ \ N\in \mathcal{B}_{1}\cup 
\mathcal{C}_{1}$\textit{\ \ then\ \ }$f$\textit{\ \ is harmonic. }

$2^{\circ })$\textit{\ }$\ $\textit{If }$\ N\in \mathcal{B}_{3}$\textit{\ \
then there is no non-constant }$\left( \pm \right) $\textit{holomorphic
harmonic map from }$M$\textit{\ into }$N.$

$3^{\circ })$\textit{\ }$\ $\textit{Let \ }$N\in \mathcal{A}_{3}$\textit{\ \
that is, \ }$N$\textit{\ \ is a trans-Sasakian manifold of type \ }$\left(
\wp ,\theta \right) $\textit{\ \ with \ }$\wp \left( q\right) =\left( \delta
\Omega _{N}\right) \left( q\right) \neq 0$\textit{\ \ and }$2n\theta \left(
q\right) =-\left( \delta \sigma \right) \left( q\right) \neq 0,$\textit{\ \
\ }$\forall \ q\in N.$\textit{\ \ Then}

$\ \ i)$\textit{\ \ \ }$f$\textit{\ \ is harmonic if and only if \ }$\ $%
\begin{equation*}
nd\lambda \left( \xi \right) +\left( \delta \sigma \right) \overset{m}{%
\underset{i=1}{\sum }}h\left( E_{i},\ E_{i}\right) =0
\end{equation*}

$\ ii)$\textit{\ \ any two of the following imply the third:}

\textit{\ \ \ }$\bullet $\textit{\ \ \ }$f$\textit{\ \ is harmonic}

\textit{\ \ \ }$\bullet $\textit{\ \ \ }$\lambda $\textit{\ \ is constant
along \ }$\xi $\textit{\ \ that is, \ }$d\lambda \left( \xi \right) =0$

$\ \ \ \bullet $\textit{\ \ }$f$\textit{\ \ is constant along \ }$D^{M}$

$iii)$\textit{\ \ there is no non-constant }$\left( \pm \right) $\textit{%
holomorphic weakly conformal harmonic map from }$M$\textit{\ into }$N.$

$iv)$\textit{\ when }$\dim M=m>n=\dim N,$\textit{\ any }$\left( \pm \right) $%
\textit{holomorphic horizontally weakly conformal }$f$\textit{\ \ is a
harmonic (and therefore harmonic morphism) if and only if\ }

\begin{equation*}
d\lambda \left( \xi \right) +\left( \delta \sigma \right) \lambda ^{2}=0.
\end{equation*}

\textit{Suppose further that either }$f$\textit{\ \ has minimal fibres or }$%
M $\textit{\ is also normal }$(i.e.\ \ M\in A_{7}\subset \mathcal{A}_{7}^{a})
$\textit{. Then }$f$\textit{\ \ is a harmonic morphism if and only if\ it is
constant.}

$4^{\circ })$\textit{\ }$\ $\textit{Let }$\ N\in A_{1}\cup A_{4}.$\textit{\
Then\ }

$\ \ \ \ \ \ i)$\textit{\ \ }$f$\textit{\ \ is harmonic if and only if }$\
\lambda $\textit{\ \ is constant along \ }$\xi ,$\textit{\ (that\ is,}$\ \
d\lambda \left( \xi \right) =0$\textit{)}$.$

$\ \ \ ii)\ $\textit{\ }$M$\textit{\ is also normal }$(i.e.\ \ M\in
A_{7}\subset \mathcal{A}_{7}^{a})$ then \textit{\ }$f$\textit{\ \ is
harmonic.}

In particular, \textit{every }$\left( \pm \right) $\textit{holomorphic map
from a semi-Sasakian manifold\textbf{\ }into a }$\alpha $-\textit{Sasakian
one is harmonic with }$\lambda $ constant.

$\ \ iii)\mathit{\ }$\textit{\ }$f$\textit{\ \ is also weakly conformal then
it is homothetic harmonic}$\mathit{\ \ \ }$\textit{.}$\ \ $ $\ \ $

$\ \ iv)$\textit{\ \ when }$\dim M=m\leq n=\dim N,$\textit{\ for a }$\left(
\pm \right) $\textit{holomorphic weakly conformal }$f$\textit{\ \ from }$M$%
\textit{\ into }$N$\textit{\ the following are equivalent:}

$\ \ \ \ \bullet $\textit{\ }$\ f$\textit{\ \ is harmonic }

$\ \ \ \ \bullet $\textit{\ }$\ \lambda \ \ $\textit{is constant and
therefore }$f$\textit{\ \ is homothetic }

$\ \ \ \ \bullet $\textit{\ }$\lambda $\textit{\ \ is constant along \ }$\xi
\ $\textit{\ }

\textit{Further If any one of those above happens then }$f$\textit{\ \ has a
minimal image.}

$\ $\textbf{Proof:}

Since \ $M$ \ is \textbf{almost} semi-Sasakian, from TABLE-II, we have%
\begin{equation*}
\ f_{\ast }\overset{m}{\underset{i=1}{\sum }}S^{M}\left( \varphi
e_{i},e_{i}\right) =-f_{\ast }\left( \nabla _{\xi }\xi \right)
\end{equation*}%
So, by Proposition $\left( 3.1\right) $ we get%
\begin{eqnarray}
\mathcal{T}\left( f\right) &=&U_{\xi }+\overset{m}{\underset{i=1}{\sum }}%
S^{N}\left( \phi E_{i},\ E_{i}\right) -\overset{m}{\underset{i=1}{\sum }}\
f_{\ast }S^{M}\left( \varphi e_{i},e_{i}\right)  \notag \\
&=&\nabla _{\left( f_{\ast }\xi \right) }f_{\ast }\xi +\overset{m}{\underset{%
i=1}{\sum }}S^{N}\left( \phi E_{i},\ E_{i}\right)  \TCItag{$4.5$}
\end{eqnarray}%
\ \ \ \ \ \ \ \ \ \ \ \ \ \ \ \ \ \ \ \ \ \ \ \ \ \ \ \ \ \ \ \ \ \ \ \ \ \
\ \ \ \ \ \ \ \ \ \ \ \ \ \ \ \ \ \ \ \ \ \ \ \ \ \ \ \ \ \ \ \ \ \ \ \ \ \
\ \ \ \ \ \ \ 

$1^{\circ }):$ \ Since \ $M$ \ is non semi-$\varphi $-involutive and \ $N$\
is $\varphi $-involutive when $N\in \mathcal{B}_{1}$, Corollary $\left(
3.1\right) /\left( 3^{\circ }\text{-}i\right) $ \ gives that\ $\lambda =0$ \
and therefore\ $\ f_{\ast }\left( \xi \right) =0.$ When $N\in \mathcal{C}%
_{1} $, we also get $f_{\ast }\left( \xi \right) =0,$ (by the $\left( \pm
\right) $\textit{holomorphicity of \ }$f$ ). Hence \ the equation $\left(
4.5\right) $ becomes%
\begin{equation*}
\mathcal{T}\left( f\right) =\overset{m}{\underset{i=1}{\sum }}S^{N}\left(
\phi E_{i},\ E_{i}\right) .
\end{equation*}%
But then, since $N\in \mathcal{B}_{1}\cup \mathcal{C}_{1},$ one has (see
TABLE-II), 
\begin{equation*}
S^{N}\left( \phi E_{i},\ E_{i}\right) =0,\ \ \forall \ i.
\end{equation*}%
Thus, harmonicity of \ $f$ \ follows.

$2^{\circ }):$ \ When $\ N\in \mathcal{B}_{3},$ \ by the same argument used
for $\ \left( 1^{\circ }\right) $, \ we see that $\ f_{\ast }\left( \xi
\right) =0$ \ and 
\begin{equation*}
\mathcal{T}\left( f\right) =\overset{m}{\underset{i=1}{\sum }}S^{N}\left(
\phi E_{i},\ E_{i}\right) .
\end{equation*}%
Also we have%
\begin{equation*}
S^{N}\left( \phi E_{i},\ E_{i}\right) =2\beta h\left( E_{i},\ E_{i}\right)
\gamma .
\end{equation*}%
Therefore one gets%
\begin{equation*}
\mathcal{T}\left( f\right) =2\beta \overset{m}{\underset{i=1}{\sum }}h\left(
E_{i},\ E_{i}\right) \gamma .
\end{equation*}%
so $f$ \ is harmonic if and only if 
\begin{equation*}
h\left( E_{i},\ E_{i}\right) =0;\ \ \ \ i=1,\cdots ,m\ \ \ 
\end{equation*}%
This gives the result.

$3^{\circ })/\left( i\right) :$ \ Note that%
\begin{equation*}
\ \left( \nabla _{\ \left( f_{\ast }\xi \right) }f_{\ast }\xi \right)
=d\lambda \left( \xi \right) \gamma .
\end{equation*}%
Since $\ N$ \ is a trans-Sasakian manifold and therefore satisfies \ $\left(
GC\right) .$ Also we have (see TABLE-II )%
\begin{equation*}
S^{N}\left( \phi E_{i},\ E_{i}\right) =\frac{1}{n}\left( \delta \sigma
\right) h\left( E_{i},\ E_{i}\right) \gamma .
\end{equation*}%
Therefore $\left( 4.5\right) $ gives that 
\begin{equation}
\mathcal{T}\left( f\right) =\left( d\lambda \left( \xi \right) +\frac{1}{n}%
\left( \delta \sigma \right) \overset{m}{\underset{i=1}{\sum }}h\left(
E_{i},\ E_{i}\right) \right) \gamma .  \tag{$4.6$}
\end{equation}%
Then the result follows.

$3^{\circ })/\left( ii\right) :$ \ \textit{\ }This is just another way of
interpreting the equation $\left( 4.6\right) .$

$3^{\circ })/\left( iii\right) :$ \ Let\textit{\ } $f$ \ also be a weakly
conformal harmonic map. For constancy of \ $f$ \ it is enough to show that
its conformal factor $\left\vert \lambda \right\vert $ \ vanishes. For this
note that 
\begin{equation*}
h\left( E_{i},\ E_{i}\right) =\lambda ^{2}g\left( e_{i},\ e_{i}\right)
=\lambda ^{2}.
\end{equation*}%
and \ $\left\vert \lambda \right\vert $ \ is also constant by the virtue of
Lemma $\left( 3.4\right) .$ So, using these and harmonicity of $f$ \ in $%
\left( 4.6\right) ,$ we get 
\begin{equation*}
\frac{m}{n}\lambda ^{2}\left( \delta \sigma \right) =0
\end{equation*}%
which gives that $\lambda =0.$

$3^{\circ })/\left( iv\right) :$ \ For a \textbf{horizontally }weakly
conformal map $f$ \ from \ $M$ \ into \ $N\in \mathcal{A}_{3},$ the function 
$\left\vert \lambda \right\vert $ becomes its dilation. W.l.o.g. choose an
orthonormal frame field%
\begin{equation*}
\left\{ \xi ,e_{1},\cdots ,e_{n},e_{n+1},\cdots ,e_{m};\varphi e_{1},\cdots
,\varphi e_{n},\varphi e_{n+1},\cdots ,\varphi e_{m}\right\}
\end{equation*}%
over $M$ \ with%
\begin{equation*}
\xi ,e_{1},\cdots ,e_{n}\in \mathcal{H=}\left( Ker\ f_{\ast }\right) ^{\perp
}.
\end{equation*}%
(We can make such a choice since $\xi \in \mathcal{H}$ and $\mathcal{H}$ is
closed under $\varphi $ \ by Lemma $\left( 3.1\right) /\left( 2^{\circ }%
\text{-}i\right) $)$.$ So note that since \ $f$ \ is \textbf{horizontally }%
weakly conformal and therefore $f_{\ast }\left( \mathcal{H}\right) =TN,$ we
have $f_{\ast }\left( \xi \right) =\lambda \gamma $ $\ $with $\ \lambda \neq
0.$ Then 
\begin{equation*}
\overset{m}{\underset{i=1}{\sum }}h\left( E_{i},\ E_{i}\right) =\lambda
^{2}n.
\end{equation*}%
So, using this in $\left( 4.6\right) $ one gets%
\begin{equation*}
\mathcal{T}\left( f\right) =\left( d\lambda \left( \xi \right) +\lambda
^{2}\left( \delta \sigma \right) \right) \gamma
\end{equation*}%
which gives the required result.

$4^{\circ })$ $\ $Since $\ N\in \mathcal{A}_{1}\cup \mathcal{A}_{4},$ \ we
see from \textbf{TABLE-II} that, $N$ \ is quasi-symplectic and satisfies $%
\left( GC\right) .$ Therefore\textit{\ }%
\begin{equation*}
\nabla _{\left( f_{\ast }\xi \right) }\mathit{\ }f_{\ast }\xi =d\lambda
\left( \xi \right) \gamma \text{ \ and \ }S^{N}\left( \phi E_{i},\
E_{i}\right) =0.
\end{equation*}%
Thus $\left( 4.5\right) $ gives us%
\begin{equation}
\mathcal{T}\left( f\right) =d\lambda \left( \xi \right) \gamma .  \tag{$4.7$}
\end{equation}

$4^{\circ })/\left( i\right) \ \ :$This follow immediately from $\left(
4.7\right) $

$4^{\circ })/\left( ii\right) \ :$This can be seen from Corollary $\left(
4.1\right) /\left( 1^{\circ }\text{-}ii\right) .$

$4^{\circ })/\left( iii\right) :$This can be seen from Corollary $\left(
4.1\right) /\left( 2^{\circ }\right) $

$4^{\circ })/\left( iv\right) :$This follows immediately from $\left(
4.7\right) $ together with Lemma $\left( 3.4\right) .$

\textbf{Corollary} $\left( 4.2\right) :$ \ 

$1^{\circ })$\textit{\ }$\ $\textit{Let \ }$f:M\rightarrow N$\textit{\ \ be
a \ }$\left( \pm \right) $-\textit{holomorphic map from an \textbf{almost}
semi-Sasakian manifold into a quasi-}$\mathcal{K}$\textit{-cosymplectic one.
Then }$f$\ \ \textit{is harmonic.}

$2^{\circ })$\textit{\ }$\ $\textit{Let \ }$f:M\rightarrow N$\textit{\ \ be
a }$\left( \pm \right) $\textit{holomorphic map from an \textbf{almost}
semi-Sasakian manifold into a }$\alpha $-\textit{Sasakian one. Then }$f$\ \ 
\textit{is harmonic if and only if \ }$d\lambda \left( \xi \right) =0.$

\textbf{Proof : \ }

$1^{\circ })$\textit{\ }$\ :$\textbf{\ }It follows directly from Theorem $%
\left( 4.4\right) /\ \left( 1^{\circ }\right) .$

$2^{\circ })$\textit{\ }$\ :$\textbf{\ }It follows directly from Theorem $%
\left( 4.4\right) /\left( \ 4^{\circ }\text{-}ii\right) .$

\textbf{Remark}$\ \left( 4.5\right) $ \ In $\left[ 8\right] ,$ Corollary $%
\left( 3.6\right) $ states that:

\textit{For semi-Sasakian }(in our terminology: \textbf{almost}
semi-Sasakian)\textit{\ manifold }$M$

$i)$\textit{\ }$\ $\textit{\ every holomorphic map from }$M$\textit{\ \ into
a quasi-}$\mathcal{K}$\textit{-cosymplectic manifold }$N$\textit{\ is
harmonic if and only if \ }$d\lambda \left( \xi \right) =0$\textit{.}

$ii)$\textit{\ }$\ $\textit{\ every holomorphic map from }$M$\textit{\ \
into a Sasakian manifold }$N$\textit{\ is harmonic if and only if \ }$%
d\lambda \left( \xi \right) =0.$

In our work

$\ \bullet $ \ Corollary $\left( 4.2\right) /\left( 1^{\circ }\right) $
improves the above result $\left( i\right) $ by showing that $\lambda =0$ \
under the circumstances and therefore removing the condition $\ "\ d\lambda
\left( \xi \right) =0\ "$ and consequently stating:

"\textit{every holomorphic map from }$M$\textit{\ \ into a quasi-}$\mathcal{K%
}$\textit{-cosymplectic manifold }$N$\textit{\ is harmonic."}

$\ \bullet $ \ Corollary $\left( 4.2\right) /\left( 2^{\circ }\right) $
recovers and generalizes the above result $\left( ii\right) .$

\textbf{Theorem }$\left( 4.5\right) :$ \ \textit{Let \ }$M\in \mathcal{A}%
_{4}.$\textit{\ \ }$\ $

$\ 1^{\circ })$\textit{\ \ If}$\ \ N\in \mathcal{A}_{1}\cup \mathcal{A}%
_{4}\cup \mathcal{B}_{1},$\textit{\ then}

\textit{\ \ \ \ }$i)$\textit{\ \ The following are equivalent}

$\ \ \ \ \ \ a^{\circ })$ $\ f$\textit{\ \ is harmonic }

$\ \ \ \ \ \ b^{\circ })$ $\ \lambda $\textit{\ \ is constant along \ }$\xi $

$\ \ \ \ \ \ c^{\circ })$ $\ \lambda $\textit{\ \ is constant }

\textit{\ \ \ }$ii)$\textit{\ \ }$f$\textit{\ \ is }$D$\textit{%
-pluriharmonic.}

$2^{\circ })$\textit{\ \ Let }$\ N\in \mathcal{A}_{1}\cup \mathcal{A}%
_{4}\cup \mathcal{B}_{1}$\ \textit{and \ }$f$\textit{\ \ is weakly
conformal. If }$\lambda $\textit{\ \ is constant along \ }$\xi $ \textit{\
then }$f$\textit{\ \ is homothetic with minimal image.}

$3^{\circ })$\textit{\ }$\ $

$\ \ \ \ i)$\textit{\ \ Let \ }$N\in \mathcal{A}_{3}\cup \mathcal{B}_{3}$\ \ 
\textit{then any two of the following imply the third:}

$\ \ \ \ \ \ \bullet $\textit{\ \ }$f$\textit{\ \ is harmonic}

$\ \ \ \ \ \ \bullet $\textit{\ \ }$\lambda $\textit{\ \ is constant along }$%
\xi $ (and therefore $\lambda $\textit{\ \ is constant)}

$\ \ \ \ \ \ \bullet $\textit{\ \ }$f$\textit{\ \ is constant along }$D^{M}$

\textit{\ \ \ }$ii)$\textit{\ \ There is no non-constant }$\left( \pm
\right) $-\textit{holomorphic \textbf{harmonic} weakly conformal map from }$%
M $\textit{\ into }$N$

$4^{\circ })$\textit{\ }$\ $

$\ \ \ i)$\textit{\ \ Let} $\ N\in \boldsymbol{C}_{1}$\textit{\ \ then }$f$%
\textit{\ \ is harmonic and }$D$\textit{-pluriharmonic.}

\textit{\ }$ii)$\textit{\ \ For }$N\in \mathcal{C}_{1},$ \textit{\ there is
no non-constant \ }$\left( \pm \right) $-\textit{holomorphic weakly
conformal map from }$M$\textit{\ into }$N$\textit{.}

$\ $\textbf{Proof: \ }

$1^{\circ }):$ $\ $Noting that $\ M$ $\ $and $\ N$ $\ $are both
quasi-symplectic, that is, 
\begin{equation*}
S^{N}\left( \phi f_{\ast }X,\ f_{\ast }X\right) =0\ \ \text{and }\
S^{M}\left( \varphi X,\ X\right) =0,\text{ \ \ }\forall \ X\in \Gamma \left(
D^{M}\right)
\end{equation*}%
we get 
\begin{equation*}
U_{X}=S^{N}\left( \phi f_{\ast }X,\ f_{\ast }X\right) -S^{M}\left( \varphi
X,\ X\right) =0,\ \ \forall \ X\in \Gamma \left( D^{M}\right) ,\ 
\end{equation*}%
from which $\left( 1^{\circ }\right) /\left( ii\right) ,$ (that is, $\ D$%
-pluriharmonicity) follows\textit{. }On the other\textit{\ }hand,
Proposition $\left( 3.1\right) $ gives that $\mathcal{T}\left( f\right)
=U_{\xi }.$ But then, since $M$ $\ $and $\ N$ $\ $both satisfy \ $\left(
GC\right) ,$ we get 
\begin{equation*}
\mathcal{T}\left( f\right) =d\lambda \left( \xi \right) \gamma ,
\end{equation*}%
from which \ the equivalence of $(a^{\circ })$ \ and $(b^{\circ })$\ of $\
\left( 1^{\circ }\right) /\left( i\right) $\ follows. For the equivalence of 
$(b^{\circ })$ \ and $(c^{\circ }),$ note that $M$\textit{\ \ and \ }$N$%
\textit{\ \ both satisfy \ }$\left( GC\right) $ and therefore, from
Proposition $\left( \left( 3.3\right) /\left( 1^{\circ }\right) \right) ,$%
\textit{\ we have }%
\begin{equation*}
d\lambda \left( X\right) =0,\ \ \forall \ X\in \Gamma \left( D^{M}\right)
\end{equation*}%
Thus, the result follows.

$2^{\circ }):$ $\ $Let $\lambda $ be constant along $\xi .$ Then, $\left(
1^{\circ }\right) /\left( i\right) $ gives that $\ f$ \ is harmonic. So, $f$
\ being homothetic with minimal image follows from Lemma $\left( 3.4\right)
. $

$3^{\circ })/\left( i\right) :$ Noting that $\ M$ $\ $and $\ N$ \ both
satisfy \ $\left( GC\right) $ \ and \ $M$ $\ $is quasi-symplectic,
Proposition $\left( 3.1\right) $ gives that%
\begin{equation}
\mathcal{T}\left( f\right) =d\lambda \left( \xi \right) \gamma +\overset{m}{%
\underset{i=1}{\sum }}S^{N}\left( \phi E_{i},\ E_{i}\right) .  \tag{$4.8$}
\end{equation}%
Since $N\in \mathcal{A}_{3}\cup \mathcal{B}_{3},$ we get \ 
\begin{equation*}
S^{N}\left( \phi E_{i},\ E_{i}\right) =\left\{ 
\begin{array}{ll}
\frac{1}{n}\left( \delta \sigma \right) h\left( E_{i},\ E_{i}\right) \gamma ,
& N\in \mathcal{A}_{3} \\ 
2\beta h\left( E_{i},\ E_{i}\right) \gamma , & N\in \mathcal{B}_{3}%
\end{array}%
\right\} .
\end{equation*}%
So, $\left( 4.8\right) $ becomes%
\begin{equation*}
\mathcal{T}\left( f\right) =\left\{ 
\begin{array}{ll}
\left( d\lambda \left( \xi \right) +\frac{1}{n}\left( \delta \sigma \right) 
\overset{m}{\underset{i=1}{\sum }}h\left( E_{i},\ E_{i}\right) \right)
\gamma , & N\in \mathcal{A}_{3} \\ 
\left( d\lambda \left( \xi \right) +2\beta \overset{m}{\underset{i=1}{\sum }}%
h\left( E_{i},\ E_{i}\right) \right) \gamma , & N\in \mathcal{B}_{3}%
\end{array}%
\right\}
\end{equation*}%
so that, the result \ follows. (Recall that \ $\delta \sigma \left( q\right)
\neq 0$ \ and \ $\beta \left( q\right) \neq 0;$ \ $\forall \ q\in N).$ \ \ \ 

$3^{\circ })/\left( ii\right) :$ Now suppose further that \ $f$ \ is weakly
conformal. Recall that if \ $f$ \ is harmonic then it is homothetic ( Lemma $%
\left( 3.4\right) $) with \textbf{constant} conformal factor \ $\lambda .$
So, for the harmonic map $f$ \ we have, when\ $\ N\in \mathcal{A}_{3},$%
\begin{eqnarray*}
\left( d\lambda \left( \xi \right) +\frac{1}{n}\left( \delta \sigma \right) 
\overset{m}{\underset{i=1}{\sum }}h\left( E_{i},\ E_{i}\right) \right) &=&%
\frac{1}{n}\left( \delta \sigma \right) \lambda ^{2}\overset{m}{\underset{i=1%
}{\sum }}g\left( e_{i},\ e_{i}\right) \\
&=&\frac{m}{n}\left( \delta \sigma \right) \lambda ^{2}=0,\ \text{\ }
\end{eqnarray*}%
and$\ $when$\ \ N\in \mathcal{B}_{4}$%
\begin{eqnarray*}
\left( d\lambda \left( \xi \right) +2\beta \overset{m}{\underset{i=1}{\sum }}%
h\left( E_{i},\ E_{i}\right) \right) &=&2\beta \lambda \overset{m}{\underset{%
i=1}{\sum }}g\left( e_{i},\ e_{i}\right) \\
&=&2m\beta \lambda =0,\ 
\end{eqnarray*}%
from which we deduce that \ $\lambda =0.$ But then any homothetic map with
vanishing conformal factor is constant. So, the required result \ follows.

$4^{\circ })/\left( i\right) $\textit{\ }$\ :$ Since $M$ $\ $is
quasi-symplectic and$\ \ N$ \ is quasi kaehler we have that $U_{X}=0,\ \
\forall \ X\in \Gamma \left( D^{M}\right) ,\ $that is, $f$\ \ is $D$%
-pluriharmonic. So, harmonic equation $\left( 3.1\right) $ gives that%
\begin{equation*}
\mathcal{T}\left( f\right) =U_{\xi }
\end{equation*}%
But then, since $\ M$ $\ $satisfies \ $\left( GC\right) $ \ and \ $f_{\ast
}\left( \xi \right) =0$ \ we get that $\ U_{\xi }=0,$which$\ $means that $f$%
\ \ is harmonic.

$4^{\circ })/\left( ii\right) $\textit{\ }$\ :$ \ Since \ \ $f_{\ast }\left(
\xi \right) =0,$\ (see Lemma $\left( 3.1\right) /\left( 2^{\circ }\right) ),$
\ $f_{\ast }$ can not be injective and therefore it has to be constant.

\textbf{Remark }$\left( 4.6\right) \ :$ \ Theorem $\left( 4.5\right) /\left(
\left( 1^{\circ }\right) /\left( i\right) \right) $ \ generalizes the result
given in $\left( \left[ 20\right] ,\text{ Theorem }\left( 1\right) \right) $
by allowing the target manifold $N$ \ to be quasi-$\mathcal{K}$-Sasakian,
quasi-contact metric, quasi-$\mathcal{K}$-cosymplectic and
nearly-cosymplectic as well as quasi-Sasakian.

\textbf{Theorem }$\left( 4.6\right) :$ \ \textit{\ \ \ }

$1^{\circ })$\textit{\ }$\ $\textit{Let\ \ }$M\in \mathcal{A}_{5}$\textit{\
\ and\ \ }$N\in \mathcal{A}_{1}\cup \mathcal{A}_{4}\cup \mathcal{B}_{1}\cup 
\mathcal{C}_{1}$\textit{\ \ then \ }$f$\textit{\ \ is harmonic if \ }$M$%
\textit{\ \ is semi-symplectic, (that is, }$\ \overset{m}{\underset{i=1}{%
\sum }}S^{M}\left( \varphi e_{i},\ e_{i}\right) =0).$

$2^{\circ })$\textit{\ }$\ $\textit{Let\ \ }$M\in \mathcal{A}_{5}$\textit{\
\ and\ \ }$N\in \mathcal{A}_{1}\cup \mathcal{A}_{4}\cup \mathcal{A}_{6}\cup 
\mathcal{A}_{7}\cup \mathcal{C}_{2}$\textit{\ \ then \ }$\left( \pm \right) $%
-\textit{holomorphic \textbf{horizontally} weakly conformal map }$f$\textit{%
\ \ from }$M$\textit{\ into }$N$\textit{\ \ is a harmonic morphism if \ }$M$%
\textit{\ \ is semi-symplectic}$.$

\textbf{Proof : \ }

$1^{\circ })$\textit{\ }$\ $Note that \ $M$ \ is a non-$\varphi $-involutive
manifold satisfying \ $(GC).$ On the other hand, \ 

$\ \bullet $ $\ N$ \ satisfies\ \ $(GC)$ when $N\in \mathcal{A}_{1}\cup 
\mathcal{A}_{4}\cup \mathcal{B}_{1}$ and therefore $\mathcal{U}_{\xi }=0$
(see Proposition $\left( 3.3\right) /\left( 2^{\circ }\right) $)

$\ \bullet $ $\ f_{\ast }\xi =0$ \ \ when $N\in \mathcal{C}_{1}$ (see Lemma $%
\left( 3.1\right) /\left( 2^{\circ }\right) $) and therefore $\mathcal{U}%
_{\xi }=0.$

Also, since \ $N$ \ is either quasi-symplectic or quasi-Kaehler, Proposition 
$\left( 3.1\right) $ gives that 
\begin{equation*}
\mathcal{T}\left( f\right) =-f_{\ast }\left( \overset{m}{\underset{i=1}{\sum 
}}S^{M}\left( \varphi e_{i},\ e_{i}\right) \right) ,
\end{equation*}%
from which the result follows.\ 

$2^{\circ })$\textit{\ }$\ $By the same argument used in $\left( 1^{\circ
}\right) $ we get $\mathcal{U}_{\xi }=0.$ On the other hand, (see Lemma $%
\left( 3.3\right) /\left( 4^{\circ }\right) $) 
\begin{equation*}
\overset{m}{\underset{i=1}{\sum }}S^{N}\left( \psi E_{i},\ E_{i}\right) =\mu
^{2}\overset{n}{\underset{i=1}{\sum }}S^{N}\left( \psi v_{i},\ v_{i}\right)
\end{equation*}%
for some local orthonormal frame field $\ $%
\begin{equation*}
\left\{ e_{1},\cdots ,e_{m};\ \varphi e_{1},\cdots ,\varphi e_{m}\right\}
\end{equation*}%
for $D^{M}$ such that the set 
\begin{equation*}
\left\{ v_{1},\cdots ,v_{n};\ \psi v_{1},\cdots ,\psi v_{n}\right\}
\end{equation*}%
forms a\ local orthonormal frame field for $D^{N}$ when $N\in \mathcal{A}%
_{1}\cup \mathcal{A}_{4}\cup \mathcal{A}_{6}\cup \mathcal{A}_{7}$ \ and \
for $TN$ \ when $N\in \mathcal{C}_{2}.$ \ Here $E_{i}=f_{\ast }e_{i}=\mu
v_{i},$ \ $\mu $ is the dilation\ and 
\begin{equation*}
\begin{tabular}{ll}
$\psi =$ & $\left\{ 
\begin{tabular}{l}
$J,$ \ \ if $\ N\in \mathcal{C}_{2}$ \\ 
$\phi ,$ \ \ otherwise%
\end{tabular}%
\right\} .$%
\end{tabular}%
\end{equation*}

So, since $N$ \ is either semi-symplectic or semi-Kaehler we have 
\begin{equation*}
\overset{m}{\underset{i=1}{\sum }}S^{N}\left( \psi E_{i},\ E_{i}\right) =0
\end{equation*}%
thus for $N\in \mathcal{A}_{1}\cup \mathcal{A}_{4}\cup \mathcal{A}_{6}\cup 
\mathcal{A}_{7}\cup \mathcal{C}_{2}$\ \ we get 
\begin{equation*}
\mathcal{T}\left( f\right) =-f_{\ast }\left( \overset{m}{\underset{i=1}{\sum 
}}S^{M}\left( \varphi e_{i},\ e_{i}\right) \right) ,
\end{equation*}%
from which the result follows.

\textbf{Theorem }$\left( 4.7\right) :$ \textit{\ Let \ }$M\in \mathcal{B}%
_{2}^{a}$, that is, $M$\textit{\ is \textbf{almost} semi-cosymplectic
manifold.}

$\ \ 1^{\circ })$\textit{\ }$\ $\textit{If\ \ }$N\in \mathcal{A}_{1}\cup 
\mathcal{A}_{2}\cup \mathcal{A}_{3}\cup \mathcal{A}_{5}\cup \mathcal{A}_{6},$%
\textit{\ then }

$\ \ \ i)\ $\textit{\ }$f$\textit{\ \ is constant along \ }$D^{M}.$\textit{\ 
}

$\ \ ii)$\textit{\ \ there is no non-constant }$\left( \pm \right) $-\textit{%
holomorphic, weakly conformal map from \ }$M$\textit{\ \ into \ }$N$\textit{.%
}

$\ 2^{\circ })$\textit{\ }$\ $\textit{If\ \ }$N\in \mathcal{A}_{1}\cup 
\mathcal{A}_{2}\cup \mathcal{A}_{3}\cup \mathcal{A}_{4}\cup \mathcal{A}%
_{5}\cup \mathcal{A}_{6},$\textit{\ then \ }$f$\textit{\ \ is harmonic if
and only if \ }$\lambda $\textit{\ \ is constant along \ }$\xi .\ \ $

$3^{\circ })$\textit{\ }$\ $\textit{If\ \ }$N\in \mathcal{B}_{1}$\textit{\ \
then}$\ \ f$\textit{\ \ is harmonic if and only if \ }$\lambda $\textit{\ \
is constant along \ }$\xi .$\textit{\ }

\textit{In particular, if }$\lambda =0$\textit{\ } \textit{(or equivalently }%
$f_{\ast }\xi =0$\textit{) and }$M\in \mathcal{B}_{1}$\textit{\ but not
nearly-cosymplectic then }$\ f$\textit{\ \ is also pluriharmonic.\ }

$\ \ \ 4^{\circ })$\textit{\ }$\ $\textit{If }$\ N\in \mathcal{B}_{1}$%
\textit{\ \ and}$\ \ f$\textit{\ \ is weakly conformal then }

$\ \ i)$ $\ f$\textit{\ \ is harmonic homothetic with minimal image if \ \ }$%
\lambda $\textit{\ \ is constant along \ }$\xi .$

\textit{\ }$ii)\ $\textit{\ }$\lambda $\textit{\ \ is constant if \ }$f$%
\textit{\ \ is harmonic.}

$5^{\circ })$\textit{\ }$\ $\textit{If\ }$\ N\in \mathcal{B}_{3}$\textit{\ \
then }

$\ i)$\textit{\ }$\ f$\textit{\ \ is harmonic if and only if}%
\begin{equation*}
d\lambda \left( \xi \right) -2\beta _{N}\overset{m}{\underset{i=1}{\sum }}%
h\left( E_{i},\ E_{i}\right) =0
\end{equation*}

$\ \ \ ii)$\textit{\ \ }

$\ \ a^{\circ }):$\textit{\ \ \ Any two of the following imply the third:}

$\ $\textit{\ \ \ \ \ \ \ }$\bullet $\textit{\ }$\ f$\textit{\ \ is harmonic 
}

$\ \ \ \ \ \ \ \ \bullet \ \lambda $\textit{\ \ is constant along \ }$\xi .$

$\ \ \ \ \ \ \ \ \bullet $\textit{\ \ }$f$\textit{\ \ is }$\ $\textit{%
constant along \ }$D^{M}\ $

\textit{\ \ \ }$b^{\circ }):$\textit{\ \ Let \ }$\lambda =0$ \textit{(or
equivalently }$f_{\ast }\xi =0$\textit{). Then \ }$f$\ \textit{\ is harmonic
if and only if it is }$\ $\textit{constant.}

$\ iii)$\textit{\ \ \ Let \ }$f$\textit{\ \ be weakly conformal. Then \ }$f$%
\textit{\ \ is harmonic if and only if it is constant \ }

$6^{\circ })$\textit{\ }$\ $\textit{If\ }$\ N\in \mathcal{C}_{^{1}}$\textit{%
\ then \ }$f$\textit{\ \ is harmonic.}$\ $

\textit{In particular, if }$M\in \mathcal{B}_{1}$\textit{\ but not
nearly-cosymplectic then\ any }$\left( \pm \right) $\textit{-holomorphic map
from \ }$M$\textit{\ \ into \ }$N$\textit{\ \ is pluriharmonic (and thus
harmonic). }$\ \ \ $

$7^{\circ })$\textit{\ }$\ $\textit{If\ }$\ \ N\ $\textit{\ is almost
Hermitian manifold then there is no non-constant }$\left( \pm \right) $%
\textit{-holomorphic weakly conformal map from \ }$M$\textit{\ \ into \ }$N.$

\textbf{Proof : }

$1^{\circ })/(i):$ \textit{\ }It follows directly from Proposition $\left(
3.2\right) /\left( 1^{\circ }\right) .$

$1^{\circ })/(ii):$ \ $f_{\ast }$\ can not be injective since \ $f$\ \ is
constant along \ $D^{M}$ and therefore it can not be weakly conformal. So, \ 
$f$\ \ has to be constant.

$2^{\circ })$\textit{\ }If $\ M$ $\in $\ $\mathcal{B}_{2}^{a},$ then $\ 
\overset{m}{\underset{i=1}{\sum }}S^{M}\left( \varphi e_{i},\ e_{i}\right)
=-\nabla _{\xi }\xi .$ So we have, by Proposition $\left( 3.1\right) ,$ that

\begin{equation}
\mathcal{T}\left( f\right) =\nabla _{\left( f_{\ast }\xi \right) }\left(
f_{\ast }\xi \right) +\overset{m}{\underset{i=1}{\sum }}S^{N}\left( \phi
E_{i},\ E_{i}\right) .  \tag{$4.9$}
\end{equation}%
But then, $S^{N}\left( \phi E_{i},\ E_{i}\right) =0,\ \forall \ i$ \ when $%
N\in \mathcal{A}_{4}$ (see TABLE-II) and also $S^{N}\left( \phi E_{i},\
E_{i}\right) =0,\ \forall \ i$ \ when $N\in \mathcal{A}_{1}\cup \mathcal{A}%
_{2}\cup \mathcal{A}_{3}\cup \mathcal{A}_{5}\cup \mathcal{A}_{6}$ \ since $\
f$\textit{\ \ is constant along \ }$D^{M}$ \ by Part $(1^{\circ })/(i).$
Therefore the equation $\left( 4.9\right) $ reduces to 
\begin{equation*}
\mathcal{T}\left( f\right) =\nabla _{\left( f_{\ast }\xi \right) }\left(
f_{\ast }\xi \right) .
\end{equation*}%
On the other hand, note that $N$ satisfies $\left( GC\right) $ and therefore
we have \ 
\begin{equation*}
\nabla _{\left( f_{\ast }\xi \right) }\left( f_{\ast }\xi \right) =d\lambda
\left( \xi \right) \gamma
\end{equation*}%
and hence 
\begin{equation*}
\mathcal{T}\left( f\right) =d\lambda \left( \xi \right) \gamma .
\end{equation*}%
This gives the result.$.$

$3^{\circ }):$ \ Since $N\in $ $\mathcal{B}_{1},$ it\ satisfies $\left(
GC\right) $ and it is also quasi-symplectic. So, we have $\ $%
\begin{equation*}
\nabla _{\left( f_{\ast }\xi \right) }\left( f_{\ast }\xi \right) =d\lambda
\left( \xi \right) \gamma \text{ \ \ and \ }S^{N}\left( \phi E_{i},\
E_{i}\right) =0,\text{\ \ \ }\forall \ i.
\end{equation*}%
Using these in $\left( 4.9\right) $ we get 
\begin{equation*}
\mathcal{T}\left( f\right) =d\lambda \left( \xi \right) \gamma ,
\end{equation*}%
from which the result follows.

In particular, asume now \ $f_{\ast }\xi =0$ (so that $\lambda =0$) and $%
M\in \mathcal{B}_{1}$\ but not nearly-cosymplectic. Then to show that $\ f$\
\ is pluriharmonic, it is enough to deal with the case where $M$ \ is quasi-$%
\mathcal{K}$-cosymplectic, (see DIAGRAM-II). For this, set \ $k=\eta \left(
Y\right) $ then we have

\textbf{Claim: }\textit{On a quasi-}$\mathcal{K}$\textit{-cosymplectic
manifold \ }$M$\textit{\ one has}%
\begin{equation*}
S^{M}\left( \varphi Y,\ Y\right) +k\left( \nabla _{Y}\xi \right) =0,\ \
\forall \ Y\in \Gamma \left( TM\right) .
\end{equation*}

Indeed, on an almost contact metric manifold, observe that $\forall \ Y\in
\Gamma \left( TM\right) $

\begin{eqnarray*}
P\left( Y,\ \varphi Y\right) &=&\left( \nabla _{Y}\varphi \right) \varphi
Y+\left( \nabla _{\left( \varphi Y\right) }\varphi \right) \varphi ^{2}Y \\
&=&\left( \nabla _{Y}\varphi \right) \varphi Y+\left( \nabla _{\left(
\varphi Y\right) }\varphi \right) \left( -Y+k\xi \right) \\
&=&\left( \nabla _{\left( \varphi Y\right) }\varphi \right) \left( k\xi
\right) -S^{M}\left( \varphi Y,\ Y\right)
\end{eqnarray*}%
But, since 
\begin{equation*}
\left( \nabla _{\left( \varphi Y\right) }\varphi \right) \left( k\xi \right)
=k\left( \nabla _{\left( \varphi Y\right) }\varphi \right) \xi =-k\varphi
\left( \nabla _{\left( \varphi Y\right) }\xi \right) ,
\end{equation*}%
we get%
\begin{equation}
P\left( Y,\ \varphi Y\right) =-\left\{ S^{M}\left( \varphi Y,\ Y\right)
+k\varphi \left( \nabla _{\left( \varphi Y\right) }\xi \right) \right\} . 
\tag{$4.10$}
\end{equation}%
Also 
\begin{equation}
\text{\ }P\left( Y,\ \xi \right) :=\left( \nabla _{Y}\varphi \right) \xi
+\left( \nabla \varphi _{Y}\varphi \right) \left( \varphi \xi \right)
=-\varphi \left( \nabla _{Y}\xi \right) .  \tag{$4.11$}
\end{equation}%
On the other hand, on a quasi-$\mathcal{K}$-cosymplectic manifold \ $M,$
from its definition, one has:%
\begin{equation*}
P\left( Y,\ \varphi Y\right) =\eta \left( \varphi Y\right) \nabla _{\left(
\varphi Y\right) }\xi =0.
\end{equation*}%
So, $\left( 4.10\right) $ gives%
\begin{equation}
S^{M}\left( \varphi Y,\ Y\right) +k\varphi \left( \nabla _{\left( \varphi
Y\right) }\xi \right) =0.  \tag{$4.12$}
\end{equation}%
Also 
\begin{equation*}
P\left( Y,\ \xi \right) =\eta \left( \xi \right) \nabla _{\left( \varphi
Y\right) }\xi =\nabla _{\left( \varphi Y\right) }\xi
\end{equation*}%
So this, together with $\left( 4.11\right) ,$ gives%
\begin{equation*}
-\varphi \left( \nabla _{Y}\xi \right) =\nabla _{\left( \varphi Y\right)
}\xi .
\end{equation*}%
Applying $\varphi $ and noting that $\ \nabla _{Y}\xi \in \Gamma \left(
D^{M}\right) ,$ \ this gives 
\begin{equation}
-\varphi ^{2}\left( \nabla _{Y}\xi \right) =\nabla _{Y}\xi =\varphi \nabla
_{\left( \varphi Y\right) }\xi .  \tag{$4.13$}
\end{equation}%
So, using $\left( 4.13\right) $ \ in $\left( 4.12\right) $ we get 
\begin{equation*}
S^{M}\left( \varphi Y,\ Y\right) +k\left( \nabla _{Y}\xi \right) =0,\ \
\forall \ Y\in \Gamma \left( TM\right) .
\end{equation*}%
which completes the proof of the claim.

Now using the assumptions that $f_{\ast }\xi =0$ $\ $and $\ \left( \pm
\right) $-holomorphicity of $\ f$, \ we see that $\forall \ Y\in \Gamma
\left( TM\right) \ $ \ \ \ 
\begin{equation*}
\sigma \left( f_{\ast }Y\right) =0\text{ \ since}\ f_{\ast }Y\ \in \Gamma
\left( D^{N}\right) \text{\ and\ \ }\ f_{\ast }\left( \nabla _{Y}\left( k\xi
\right) \right) =kf_{\ast }\left( \nabla _{Y}\xi \right)
\end{equation*}%
So this and Lemma $\left( 3.2\right) /\left( 1^{\circ }\text{-}ii\right) $
give that,$\ \forall \ Y\in \Gamma \left( TM\right) $

\begin{equation*}
U\left( Y,Y\right) =S^{N}\left( \phi f_{\ast }Y,\ f_{\ast }Y\right) -f_{\ast
}\left\{ S^{M}\left( \varphi Y,\ Y\right) \ +k\left( \nabla _{Y}\xi \right)
\right\}
\end{equation*}%
Using the claim this gives \ 
\begin{equation*}
U\left( Y,Y\right) =S^{N}\left( \phi f_{\ast }Y,\ f_{\ast }Y\right) .
\end{equation*}%
But then, since $N\in \mathcal{B}_{1}$ and thus it is quasi-symplectic, (see
TABLE II ), $U\left( Y,Y\right) =0,\ $\ $\forall \ Y\in \Gamma \left(
TM\right) .$ That is, $f$ \ is pluriharmonic. This ends the proof.

$4^{\circ })$ We know from \ $(3^{\circ })$ that 
\begin{equation*}
\mathcal{T}\left( f\right) =d\lambda \left( \xi \right) \gamma .
\end{equation*}%
So, $\left( i\right) $ and $\left( ii\right) $ follows from Lemma $\left(
3.4\right) .$

$5^{\circ })$\ \ Since $N\ $\ satisfies $\left( GC\right) ,$ \textit{we have}%
\begin{equation*}
\nabla _{\left( f_{\ast }\xi \right) }\left( f_{\ast }\xi \right) =d\lambda
\left( \xi \right) \gamma .
\end{equation*}%
On the other hand, since $N\in \mathcal{B}_{3}$ \ 
\begin{equation*}
S^{N}\left( \phi E_{i},\ E_{i}\right) =-2\beta _{N}h\left( E_{i},\
E_{i}\right) \gamma ,\text{ \ \ }\forall \ i=1,\cdots ,m
\end{equation*}%
So, the equation $\left( 4.9\right) $ ,which is valid for such an $f$ \
under consideration,\ becomes%
\begin{equation}
\mathcal{T}\left( f\right) =\left( d\lambda \left( \xi \right) -2\beta _{N}%
\overset{m}{\underset{i=1}{\sum }}h\left( E_{i},\ E_{i}\right) \right)
\gamma .  \tag{$4.14$}
\end{equation}%
Parts $\left( 5^{\circ }\right) /$\ $\left( \left( i\right) ,\ \left(
ii\right) /\mathit{\ }a^{\circ }\right) $ follows from \ $\left( 4.14\right)
.$ Also $\left( 5^{\circ }\right) /$\ $\left( \left( ii\right) /\mathit{\ }%
b^{\circ }\right) $ follows from $\left( ii\right) /\mathit{\ }a^{\circ }.$

$iii)\ :$ \ \textit{\ }When\ \ $f$ \ is weakly conformal, note that 
\begin{equation*}
\overset{m}{\underset{i=1}{\sum }}h\left( E_{i},\ E_{i}\right) =m\lambda
^{2}.
\end{equation*}%
Also, from Lemma $\left( 3.4\right) $ the function $\lambda $ ( which is
also the conformal factor of \ $f$ ) is constant and therefore $\ d\lambda
\left( \xi \right) =0.$ So, from $\left( 4.14\right) $ we get%
\begin{equation*}
\mathcal{T}\left( f\right) =2m\beta _{N}\lambda ^{2}\gamma ,
\end{equation*}%
from which we have that $f$ \ is harmonic if and only if \ $\lambda $
vanishes. But then, any weakly conformal map with vanishing conformal factor
is constant. This completes the proof.

$6^{\circ })$ \ Note that since $\ N=H^{2n}\in \mathcal{C}_{1}$ one has $%
f_{\ast }\xi =0\ \ $and $H$ is quasi-Kaehler, that is,$\ $%
\begin{equation*}
S^{N}\left( JE_{i},\ E_{i}\right) =0,\ \forall \ i=1,\cdots ,m.
\end{equation*}%
\ So, the result\ follows from the equation $\left( 4.9\right) .\ $

In particular, pluriharmonicity of \ $f$\ \ from a quasi-$\mathcal{K}$%
-cosymplectic $M$ into a\textit{\ }quasi-Kaehler \ $N$ \ follows by the same
argument used in $\left( 2^{\circ }\right) $ with some minor adjustments.$\ $

$7^{\circ })\ \ \ f_{\ast }$ cannot be injective since $\ $ $f_{\ast }\xi
=0, $ and therefore a non-constant $f$ cannot be weakly conformal.

\textbf{Remark} \ $\left( 4.7\right) :$ \ \ \ \ \ \ \ \ \ \ \ \ \ \ \ \ \ \
\ \ \ \ \ \ \ \ \ \ \ \ \ \ \ \ \ \ \ \ \ \ \ \ \ \ \ \ \ \ \ \ \ \ \ \ \ \
\ \ \ \ \ \ \ \ \ \ \ \ \ \ \ \ \ \ \ \ \ \ \ \ \ \ \ \ \ \ \ \ \ \ \ \ \ \
\ \ \ \ \ \ \ \ \ \ \ \ \ \ \ \ \ \ \ \ \ \ \ \ \ \ \ \ \ \ \ \ \ \ \ \ \ \
\ \ \ \ \ \ \ \ \ \ \ \ \ \ \ \ \ \ \ \ \ \ \ \ \ \ \ \ \ \ \ \ \ \ \ \ \ \
\ \ \ \ \ \ \ \ \ \ \ \ \ \ \ \ \ \ \ \ \ \ \ \ \ \ \ \ \ \ \ \ \ \ \ \ \ \
\ \ \ \ \ \ \ \ \ \ \ \ \ \ \ \ \ \ \ \ \ \ \ \ \ \ \ \ \ \ \ \ \ \ \ \ \ \
\ \ \ \ \ \ \ \ \ \ \ \ \ \ \ \ \ \ \ \ \ \ \ \ \ \ \ \ \ \ \ \ \ \ \ \ \ \
\ \ \ \ \ \ \ \ \ \ \ \ \ \ \ \ \ \ \ \ \ \ \ \ \ \ \ \ \ \ \ \ \ \ \ \ \ \
\ \ \ \ \ \ \ \ \ \ \ \ \ \ \ \ \ \ \ \ \ \ \ \ \ \ \ \ \ \ \ \ \ \ \ \ \ \
\ \ \ \ \ \ \ \ \ \ \ \ \ \ \ \ \ \ \ \ \ \ \ \ \ \ \ \ \ \ \ \ \ \ \ \ \ \
\ \ \ \ \ \ \ \ \ \ \ \ \ \ \ \ \ \ \ \ \ \ \ \ \ \ \ \ \ \ \ \ \ \ \ \ \ \
\ \ \ \ \ \ \ \ \ \ \ \ \ \ \ \ \ \ \ \ \ \ \ \ \ \ \ \ \ \ \ \ \ \ \ \ \ \
\ \ \ \ \ \ \ \ \ \ \ \ 

$1^{\circ })$ In $\left[ 8\right] $ Corollary $\left( 3.5\right) $ states
that:

$\ \bullet $\ \ \textit{Let \ }$f:M\rightarrow N$\textit{\ \ be a
holomorphic map from a semi-cosymplectic (in our terminology: almost
semi-cosymplectic) manifold into a quasi-}$\mathcal{K}$\textit{-cosymplectic
manifold or Sasakian manifold. Then }$f$\textit{\ \ is harmonic if and only
if \ }$d\lambda \left( \xi \right) =0.$

In our work, Theorem $\left( 4.7\right) /\left( 3^{\circ }\text{ \textit{and 
}}2^{\circ }\right) $\ \ generalize this result by allowing that the target
manifold,

$\ \bullet \ N$ \ to be nearly-cosymplectic as well as quasi-$\mathcal{K}$%
-cosymplectic and

$\ \bullet \ N\ \ $to be in $\mathcal{A}_{1}\cup \mathcal{A}_{2}\cup 
\mathcal{A}_{3}\cup \mathcal{A}_{4}\ \cup \mathcal{A}_{5}\cup \mathcal{A}%
_{6}\ $ \ which includes the case where $N$ is Sasakian.

(see \textbf{Diagram-I }and\textbf{\ \ Diagram-II}):

$2^{\circ })$ In $\left[ 5\right] $ Proposition $\left( 3.1\right) $ states
that:

$\ \bullet $\ \ Any $\left( \varphi ,\ J\right) $-holomorphic mapping from a
cosymplectic manifold \ $M$\ \ into a Kaehler manifold \ $N$\ \ is
pluriharmonic (and thus harmonic), $\ $

In our work, Theorem $\left( 4.7\right) /\left( 6^{\circ }\right) \ $\ \
generalizes this result by allowing the domain manifold $M$ to be quasi-$%
\mathcal{K}$-cosymplectic (which covers the nearly-$\mathcal{K}$%
-cosymplectic, almost cosymplectic as well as cosymplectic case, see \textbf{%
Diagram-II}) and the target manifold $N$ to be quasi-Kaehler (which covers
the almost Kaehler and nearly Kaehler cases as well as Kaehler cases, see 
\textbf{Diagram-III}). In this context, see also Theorem $\left( 4.7\right)
/\left( 3^{\circ }\right) .$

\textbf{Theorem }$\left( 4.8\right) :$ \ \textit{Let }\ $M\in \mathcal{B}%
_{3} $ \ ($i.e.$ $M$ \ \ \textit{is an almost }$\beta _{M}$\textit{-Kenmotsu
manifold)}.\textit{\ \ \ }$\ $

$1^{\circ })$ $\ $\textit{If }$\ N\in \mathcal{A}_{1}\cup \mathcal{A}%
_{2}\cup \mathcal{A}_{3}\cup \mathcal{A}_{4}\cup \mathcal{A}_{5}\cup 
\mathcal{A}_{6}\cup \mathcal{B}_{1},$

\ \ \ $i)$ \textit{then }$\ f$ \ \textit{is harmonic if and only if}%
\begin{equation*}
d\lambda \left( \xi \right) +2\lambda m\beta _{M}=0
\end{equation*}

$\ \ ii)$ $\ $\textit{and suppose further that} $\ \beta _{M}=\alpha \ \in 
%TCIMACRO{\U{211d} }%
%BeginExpansion
\mathbb{R}
%EndExpansion
^{+},$ $\ $\textit{then} $\ f$ \ \textit{is harmonic if and only if } $%
\forall \ p\in M,$ there is a curve $\ $%
\begin{equation*}
\vartheta \ {\large =\ }\vartheta _{p}{\large :}\left( -{\large \varepsilon
,\ \ \varepsilon }\right) \subset 
%TCIMACRO{\U{211d} }%
%BeginExpansion
\mathbb{R}
%EndExpansion
\rightarrow M
\end{equation*}%
\ with \ $\vartheta \left( 0\right) =p$ \ and \ $\frac{{\LARGE d}{\large %
\vartheta }}{{\LARGE dt}}\left( 0\right) =\xi _{p}$ \ satisfying%
\begin{equation*}
\lambda \circ {\normalsize \vartheta }\left( t\right) =\lambda \left(
p\right) {\Large e}^{-2{\large m\alpha t}},
\end{equation*}

$2^{\circ })$ $\ $\textit{If }$\ N\in \mathcal{A}_{1}\cup \mathcal{A}%
_{2}\cup \mathcal{A}_{3}\cup \mathcal{A}_{4}\cup \mathcal{A}_{5}\cup 
\mathcal{A}_{6}\cup \mathcal{B}_{1},\ $\textit{then, there is no
non-constant }$\left( \pm \right) $\textit{-holomorphic\ weakly conformal
harmonic map from \ }$M$\textit{\ \ into }$N$\textit{\ .}

$3^{\circ })$ $\ $\textit{If }$\ N\in \mathcal{B}_{3}$

$\ \ \ i)$ \textit{then} \ $f$ \ \textit{is harmonic if and only if \ }$%
\forall \ p\in M$\textit{\ \ }%
\begin{equation*}
\left( d\lambda \right) _{p}\left( \xi \right) +2\lambda _{p}m\beta
_{M}\left( p\right) -2\beta _{N}\left( q\right) \overset{m}{\underset{i=1}{%
\sum }}h_{q}\left( E_{i},\ E_{i}\right) =0
\end{equation*}%
or equivalently 
\begin{equation*}
\left( d\lambda \right) _{p}\left( \xi \right) +2\lambda _{p}m\beta
_{M}\left( p\right) -2\beta _{N}\left( q\right) \left( trace_{g}\left(
f^{\ast }h\right) -\lambda ^{2}\right) =0.
\end{equation*}%
where \ $q=f\left( p\right) $

\ \ $ii)\ $any two of the following imply the third:

$\ \ \ \ \ \ \bullet $ $\ f$ \ \textit{is harmonic }

\ \ \ \ \ \ $\bullet $ $\ d\lambda \left( \xi \right) +2\lambda m\beta
_{M}=0 $

\ \ \ \ \ \ $\bullet $ \ $f$ \ \textit{is constant along} \ $D^{M}$

$4^{\circ })$ \textit{If \ }$N\in \mathcal{B}_{3}$\textit{\ \textit{and}} 
\textit{\ }$f$\ \textit{\ is a \textbf{non-constant} weakly conformal map,
then the following are equivalent:}

$\ \bullet $ \ \textit{\ }$f$\ \ \textit{is harmonic }

$\ \bullet $ \ \textit{\ }$f$\ \ \textit{is homothetic and is of minimal
image with constant conformal factor }%
\begin{equation*}
\lambda \left( p\right) =\frac{\beta _{M}\left( p\right) }{\beta _{N}\left(
q\right) }=\lambda _{0},\ \ \forall \ p\in M
\end{equation*}%
where \ $q=f\left( p\right) .$\textit{\ }

$\ \bullet $ $\ \ \left( d\lambda \left( \xi \right) +2m\lambda \beta
_{M}-2m\lambda ^{2}\beta _{N}\right) =0.$

\textit{\ In particular, if \ }$\beta _{M}=\alpha =\beta _{N},$\textit{\ for
some number }$\alpha >0,$\textit{\ then \ }$f$\textit{\ \ is harmonic if and
only if it is an isometric minimal immersion.}

$5^{\circ })$ $\ $\textit{If }$\ N=H\in \mathcal{C}_{1}$ \ \textit{then} \ $%
f $ \ \textit{is harmonic.}

$6^{\circ })\ \ $\textit{If\ \ }$N=H\in \mathcal{C}_{2}$ \textit{\ then} \ 
\textit{there is \textbf{no} non-constant }$\left( \pm \right) $\textit{%
-holomorphic weakly conformal map from \ an almost }$\ \beta _{M}$-\textit{%
Kenmotsu manifold }$M$\textit{\ \ into \ semi-Kaehler manifold }$H$\textit{. 
}

\textbf{Proof :}

$\left( 1^{\circ }\right) /(i)\ :$ Note that since $\ M\in \mathcal{B}_{3},$
\ we have%
\begin{equation*}
\overset{m}{\underset{i=1}{\sum }}S^{M}\left( \varphi e_{i},\ e_{i}\right)
=-2m\beta _{M}\xi ,
\end{equation*}%
and\ since $\ M$ and $\ N$ \ satisfy \ $(GC)$ \ we have%
\begin{equation*}
U_{\xi }=\left( d\lambda \left( \xi \right) \right) \gamma .
\end{equation*}%
So, Proposition $\left( 3.1\right) $ gives%
\begin{equation}
\mathcal{T}\left( f\right) =\left( d\lambda \left( \xi \right) +2m\lambda
\beta _{M}\right) \gamma +\overset{m}{\underset{i=1}{\sum }}S^{N}\left( \phi
E_{i},\ E_{i}\right) .  \tag{$4.15$}
\end{equation}%
Now\ note that $S^{N}\left( \phi E_{i},\ E_{i}\right) =0,$ when$\ \ N\in 
\mathcal{A}_{1}\cup \mathcal{A}_{4}\cup \mathcal{B}_{1}.$ When $\ N\in 
\mathcal{A}_{2}\cup \mathcal{A}_{3}\cup \mathcal{A}_{5}\cup \mathcal{A}_{6},$
we also have $S^{N}\left( \phi E_{i},\ E_{i}\right) =0$ by\ Proposition $%
\left( 3.2\right) ,$ since $\ N$ \ is \textbf{non}-$\varphi $-involutive
while $\ M$ \ is $\varphi $-involutive and therefore $E_{i}=f_{\ast
}e_{i}=0. $ \ Thus\ the equation $\left( 4.15\right) $ \ becomes \ \ $\ $ $\
\ \ $ \ \ 
\begin{equation}
\mathcal{T}\left( f\right) =\left( d\lambda \left( \xi \right) +2m\lambda
\beta _{M}\right) \gamma  \tag{$4.16$}
\end{equation}%
This gives the result.$\ \ $

\bigskip $\left( 1^{\circ }\right) /(ii)\ :$ \ From part $\left( i\right) $
we have that, for \ $\beta _{M}=\alpha >0$,

$\bullet $ \ $f$ \ is harmonic if and only if 
\begin{equation*}
d\lambda \left( \xi \right) +2m\lambda \alpha =0.
\end{equation*}%
But then solving this differential equation gives the result.

$2^{\circ }):$ \ Let \ $f$\ \ be also weakly conformal. Since $f$ \ is $%
\left( \pm \right) $-holomorphic, \ $\lambda $ becomes its conformal factor.
\ So, by the virtue of Lemma $\left( 3.4\right) ,$ if \ $f$\ \ is harmonic
then its conformal factor\ $\lambda $ is constant and therefore $d\lambda
\left( \xi \right) =0.$ So, $\left( 4.16\right) $ \ becomes 
\begin{equation*}
\mathcal{T}\left( f\right) =2m\beta \lambda \gamma ,
\end{equation*}%
we see that \ $\lambda \ $\ vanishes if $f$\ \ is harmonic. But then,
vanishing of conformal factor $\lambda $ gives the constancy of $\ f$\ .

$3^{\circ })\ :$ \ 

$\ i):$ Since $M,N\in \mathcal{B}_{3},$

$\bullet $ \ they both satisfy \ $\left( GC\right) $ \ and therefore $U_{\xi
}=d\lambda \left( \xi \right) ,$

$\bullet $ \ we get 
\begin{equation*}
f_{\ast }\left( \overset{m}{\underset{i=1}{\sum }}S^{M}\left( \varphi
e_{i},\ e_{i}\right) \right) =-2m\beta _{M}f_{\ast }\left( \xi \right)
=-2m\left( \lambda \beta _{M}\right) \gamma \text{\ }
\end{equation*}%
and 
\begin{equation*}
\overset{m}{\underset{i=1}{\sum }}S^{N}\left( \phi E_{i},\ E_{i}\right)
=-2\beta _{N}\left( \overset{m}{\underset{i=1}{\sum }}h\left( E_{i},\
E_{i}\right) \right) \gamma =-2\beta _{N}\left( trace_{g}\left( f^{\ast
}h\right) -\lambda ^{2}\right) \gamma
\end{equation*}%
So, Proposition $\left( 3.1\right) $ gives%
\begin{equation}
\mathcal{T}\left( f\right) =\left( d\lambda \left( \xi \right) +2m\beta
_{M}\lambda -2\beta _{N}\overset{m}{\underset{i=1}{\sum }}h\left( E_{i},\
E_{i}\right) \right) \gamma  \tag{$4.17$}
\end{equation}%
or equivalently 
\begin{equation*}
\mathcal{T}\left( f\right) =\left\{ d\lambda \left( \xi \right) +2m\beta
_{M}\lambda -2\beta _{N}\left( trace_{g}\left( f^{\ast }h\right) -\lambda
^{2}\right) \right\} \gamma ,
\end{equation*}%
from which the result follows.

$\ ii):$ This is just another interpretation of the result in $\left(
i\right) .$

$4^{\circ }):$ \ \ Let \ $N\in \mathcal{B}_{3}$ $\ $and $f$\ \ be also a
non-constant weakly conformal map. \ Then the equation $\left( 4.17\right) $
is valid, that is, 
\begin{equation*}
\mathcal{T}\left( f\right) =\left( d\lambda \left( \xi \right) +2m\beta
_{M}\lambda -2\beta _{N}\overset{m}{\underset{i=1}{\sum }}h\left( E_{i},\
E_{i}\right) \right) \gamma
\end{equation*}%
But then, since 
\begin{equation*}
\overset{m}{\underset{i=1}{\sum }}h\left( E_{i},\ E_{i}\right) =\overset{m}{%
\underset{i=1}{\sum }}\lambda ^{2}g\left( e_{i},\ e_{i}\right) =m\lambda
^{2}.
\end{equation*}%
we get 
\begin{equation*}
\mathcal{T}\left( f\right) =\left\{ d\lambda \left( \xi \right) +2m\beta
_{M}\lambda -2m\lambda ^{2}\beta _{N}\right\} \gamma ,
\end{equation*}%
from which the equivalence of first and third statements follow. The
equivalence of first and second follows from Lemma $\left( 3.4\right) .$
Also observe that the conformal factor $\lambda $ is constant since\ \textit{%
\ }$f$\ \ \textit{is homothetic, and therefore }$d\lambda \left( \xi \right)
=0.$ So the third statement gives that\textit{\ } 
\begin{equation*}
\lambda \left( p\right) =\frac{\beta _{M}\left( p\right) }{\beta _{N}\left(
q\right) }=\lambda _{0},\ \ \forall \ p\in M.
\end{equation*}

In particular, if \ $\beta _{M}=\alpha =\beta _{N}$ \ then \ $\lambda =1$
and therefore $\ f$ \ becomes isometric minimal immersion.

$5^{\circ })$ \ $:$ \medskip Note that $\ f_{\ast }\left( \xi \right) =0$ \
and \ $S^{N}\left( \phi E_{i},\ E_{i}\right) =0,\ \ \forall \ i$ \ since \ $%
N\in \mathcal{C}_{1}$, \ and \ $\nabla _{\xi }\xi =0$ $\ $since\ $\ M$ \
satisfies $\left( GC\right) .$ Then $U_{\xi }=0$ \ and therefore, by
Proposition $\left( 3.1\right) ,$ \ we have 
\begin{equation*}
\mathcal{T}\left( f\right) =-f_{\ast }\overset{m}{\underset{i=1}{\sum }}%
S^{M}\left( \varphi e_{i},\ e_{i}\right) =-\left( 2m\beta _{M}\right)
f_{\ast }\xi =0.
\end{equation*}%
This completes the proof.

$6^{\circ })$ \ $:$ This is just Lemma $\left( 3.3\right) /\left( ii\right)
. $

\textbf{Remark} $(4.8):$ For a $\left( \varphi ,\phi \right) $\textit{%
-holomorphic map }$f$\textit{\ }$:M$\textit{\ }$\rightarrow N$\textit{\ \
between almost contact metric manifolds}

$1^{\circ })$\ In $\left( \left[ 29\right] ,\text{Theorem}\left( 4.1\right)
\right) $ states that

$\ \ \bullet $ If $M$ and $N$ are both almost Kenmotsu then \textit{\ }$f$%
\textit{\ \ is harmonic if and only if \ }%
\begin{equation*}
trace_{g}\left( f^{\ast }h\right) =\lambda ^{2}+2m\lambda +d\lambda \left(
\xi \right) .
\end{equation*}

However, in our work, Theorem $\left( 4.8\right) /\ \left( 3^{\circ }\text{-}%
i\right) $ gives an alternative result, namely:

$\ \ \bullet $ If $M$ and $N$ are both almost Kenmotsu then \textit{\ }$f$%
\textit{\ \ is harmonic if and only if \ }%
\begin{equation*}
trace_{g}\left( f^{\ast }h\right) =\lambda ^{2}+m\lambda +\frac{1}{2}%
d\lambda \left( \xi \right) .
\end{equation*}

$2^{\circ })$\ In $\left( \left[ 29\right] ,\text{Theorem}\left( 4.3\right)
\right) $ states that

$\ \ \bullet $ \ \ \textit{Any }$\left( \varphi ,\phi \right) $\textit{%
-holomorphic map }$f$\textit{\ }$:M$\textit{\ }$\rightarrow N$\textit{\ \
from an almost Kenmotsu manifold into a contact metric manifold\ is harmonic
if and only if \ }%
\begin{equation*}
d\lambda \left( \xi \right) =0.
\end{equation*}

But again, in our work, Theorem $\left( 4.8\right) /\ \left( 1^{\circ }\text{%
-}i\right) $ gives a different result, namely:

$\ \ \bullet $ \ \ \textit{Any }$\left( \varphi ,\phi \right) $\textit{%
-holomorphic map }$f$\textit{\ }$:M$\textit{\ }$\rightarrow N$\textit{\ \
from an almost Kenmotsu manifold into a contact metric manifold\ is harmonic
if and only if \ }%
\begin{equation*}
d\lambda \left( \xi \right) +2m\lambda =0.
\end{equation*}

$3^{\circ })$\ In $\left( \left[ 28\right] ,\text{Theorem}\left( 3.1\right)
\right) $ states that

$\ \ \bullet $ \ \ \textit{Any }$\left( \varphi ,J\right) $\textit{%
-holomorphic map }$f$\textit{\ }$:M$\textit{\ }$\rightarrow H$\textit{\ \
from a Kenmotsu manifold into a Kaehler one\ is harmonic.}

Our Theorem $\left( 4.8\right) /\ \left( 5^{\circ }\right) $\textit{\ }%
generalizes this result by allowing the domain manifold $M$ to be \textbf{%
almost} $\beta _{M}$-Kenmotsu (which covers the Kenmotsu cases as well as
the $\beta _{M}$-Kenmotsu or $CR$-integrable almost $\beta _{M}$-Kenmotsu
cases) and the target manifold\ $H$ to be quasi Kaehler (which covers
Kaehler cases as well as nearly-Kaehler and almost-Kaehler cases).

\textbf{Theorem }$\left( 4.9\right) :$

$1^{\circ }\mathbf{\ })$\textit{\ \ Let \ }$M=H=\left( H,\ J,\ G\right) \in 
\mathcal{C}_{2}.$\textit{\ \ }$\ $

$\ i)\ \ $\textit{If \ }$N\in \mathcal{A}_{1}\cup \mathcal{A}_{4}^{a}\cup 
\mathcal{B}_{1}\cup \mathcal{C}_{1}$\textit{\ \ then \ }$f$\textit{\ \ is
harmonic.}

In particular, if $\ M=H\in \mathcal{C}_{1}$ \textit{and \ }$N\in \mathcal{A}%
_{1}\cup \mathcal{A}_{4}^{a}\cup \mathcal{B}_{1}\cup \mathcal{C}_{1}$ 
\textit{then \ }$f$\textit{\ \ is pluriharmonic (and therefore harmonic).}

$ii)$\textit{\ }$\ $\textit{If \ }$N\in \mathcal{A}_{3}\cup \mathcal{B}_{3}$%
\textit{\ then \ }$f$\textit{\ \ is harmonic if and only if it is constant. }

$iii)$\textit{\ }$\ $\textit{If \ }$N\in \mathcal{B}_{2}$ \ \textit{(i.e. \ }%
$N$\textit{\ \ is a semi-cosymplectic \textbf{\ }) with \ }$\dim N=1+\dim H$ 
\textit{and} \textit{\ }$f$\textit{\ \ is weakly conformal then \ }$f$%
\textit{\ is harmonic. }

\textbf{Proof : }$\mathbf{\ }$Note that $\ \overset{m}{\underset{i=1}{\sum }}%
S^{H}\left( Je_{i},\ e_{i}\right) =0,$\ when \ $M\in \mathcal{C}_{2}$ and$\
\ S^{H}\left( Je_{i},\ e_{i}\right) =0,\ \forall \ i\ \ \ $when \ $M\in 
\mathcal{C}_{1}.$

$\bigskip 1^{\circ })/\left( i\right) :$ \ Note also that $\ S^{N}\left(
\phi E_{i},\ E_{i}\right) =0,\ \forall \ i\ ;\ \ $if \ $N\in \mathcal{A}%
_{1}\cup \mathcal{A}_{4}^{a}\cup \mathcal{B}_{1}$ and $\ S^{N}\left(
J_{1}E_{i},\ E_{i}\right) =0,$ $\forall \ i\ ;\ \ $if \ $N=H\in \mathcal{C}%
_{1}.$ So, 
\begin{equation*}
S^{N}\left( \psi E_{i},\ \ E_{i}\right) =0,\forall \ i\ 
\end{equation*}%
where

\begin{equation*}
\ 
\begin{array}{ll}
\psi = & \left\{ 
\begin{array}{l}
\phi ,\text{ for \ }N\in \mathcal{A}_{1}\cup \mathcal{A}_{4}^{a}\cup 
\mathcal{B}_{1} \\ 
J_{1},\text{ for \ }N=H_{1}=\left( H_{1},\ J_{1},\ G_{1}\right) \in \mathcal{%
C}_{1}%
\end{array}%
\right\} .%
\end{array}%
\end{equation*}%
So, from Proposition $\left( 3.1\right) /\left( 3^{\circ }\right) $ \ we have%
\begin{equation*}
\mathcal{T}\left( f\right) =\ \overset{m}{\underset{i=1}{\sum }}\left\{
S^{N}\left( \psi E_{i},\ \ E_{i}\right) -S^{H}\left( Je_{i},\ e_{i}\right)
\right\} ,
\end{equation*}%
from which, harmonicity (and in particular pluriharmonicity) of $f$ \ will
follow.

$1^{\circ })/\left( ii\right) :$ \ Since \ $\overset{m}{\underset{i=1}{\sum }%
}S^{H}\left( Je_{i},\ e_{i}\right) =0,$\ \ we have 
\begin{equation}
\mathcal{T}\left( f\right) =\ \overset{m}{\underset{i=1}{\sum }}S^{N}\left(
\phi E_{i},\ E_{i}\right) .  \tag{$4.18$}
\end{equation}%
by Proposition $\left( 3.1\right) /\left( 3^{\circ }\right) .$ But then, \ 
\begin{equation*}
S^{N}\left( \phi E_{i},\ E_{i}\right) =\left\{ 
\begin{tabular}{ll}
$\frac{\delta \eta }{m}h\left( E_{i},\ E_{i}\right) \gamma ,$ & when $N\in 
\mathcal{A}_{3}$ \\ 
$-2\beta _{N}h\left( E_{i},\ E_{i}\right) \gamma ,$ & when $N\in \mathcal{B}%
_{3}$%
\end{tabular}%
\right\} \ 
\end{equation*}%
So, the equation $\left( 4.18\right) $ \ gives%
\begin{equation*}
\mathcal{T}\left( f\right) =\left\{ 
\begin{tabular}{ll}
$\frac{\delta \eta }{m}\overset{m}{\underset{i=1}{\sum }}h\left( E_{i},\
E_{i}\right) \gamma ,$ & when $N\in \mathcal{A}_{3}$ \\ 
$-2\beta _{N}\overset{m}{\underset{i=1}{\sum }}h\left( E_{i},\ E_{i}\right)
\gamma ,$ & when $N\in \mathcal{B}_{3}$%
\end{tabular}%
\right\}
\end{equation*}%
from which the result follows.

$1^{\circ })/\left( iii\right) :$ Observe that%
\begin{eqnarray*}
\ f_{\ast }\left( TH\right) &=&D^{N}=span\left\{ f_{\ast }e_{1},\cdots ,\
f_{\ast }e_{m};\ f_{\ast }\varphi e_{1},\cdots ,\ f_{\ast }\varphi
e_{m}\right\} \\
&=&span\left\{ f_{\ast }e_{1},\cdots ,\ f_{\ast }e_{m};\ \phi f_{\ast
}e_{1},\cdots ,\phi \ f_{\ast }e_{m}\right\} \\
&=&span\left\{ \ \lambda v_{1},\cdots ,\ \lambda v_{m};\lambda \phi
v_{1},\cdots ,\ \lambda \phi v_{m}\right\}
\end{eqnarray*}%
\ since \ $f_{\ast }$ \ is injective\ (as $f$ \ is weakly conformal and $%
\lambda $ \ becomes the conformality factor) and by the assumptions that \ $%
f $ \ is $\left( \pm \right) \left( J,\phi \right) $-holomorphic and $\dim
N=1+\dim H$. \ Here $\left\{ \ v_{1},\cdots ,\ v_{m};\phi v_{1},\cdots ,\
\phi v_{m}\right\} $\ forms a local orthonormal frame field for $D^{N},$ $\ $%
where \ $f_{\ast }e_{i}=E_{i}=\lambda v_{i}$. On the other hand, we see that 
$\ $%
\begin{equation*}
\overset{m}{\underset{i=1}{\sum }}S^{N}\left( \phi E_{i},\ E_{i}\right)
=\lambda ^{2}\overset{m}{\underset{i=1}{\sum }}S^{N}\left( \phi v_{i},\
v_{i}\right) =0,
\end{equation*}
\ since \ $N\in \mathcal{B}_{2}$\textit{\ (i.e. \ }$N$\textit{\ \ is
semi-cosymplectic ). }So we get 
\begin{equation*}
\mathcal{T}\left( f\right) =\overset{m}{\underset{i=1}{\sum }}S^{N}\left(
\phi E_{i},\ E_{i}\right) -\ f_{\ast }\overset{m}{\underset{i=1}{\sum }}%
S^{H}\left( Je_{i},\ e_{i}\right) =0.
\end{equation*}%
That is, \ $f$ \ is harmonic.

\textbf{Remark }$\left( 4.9\right) $ $:$

In $\left[ 28\right] ,$Theorem $\left( 3.2\right) $ states that

$\ \ \bullet $ \ \ \textit{Any }$\left( J,\phi \right) $\textit{-holomorphic
map }$f$\textit{\ }$:M$\textit{\ }$\rightarrow N$\textit{\ \ from a Kaehler
manifold into a Kenmotsu manifold\ is harmonic if and only if \ it is
constant.}

Our Theorem $\left( 4.9\right) /\left( 1^{\circ }\text{- }ii\right) \ $\ \
generalizes this result by allowing the domain manifold $M$ to be
semi-Kaehler (which covers the Kaehler cases as well as quasi-Kaehler or
nearly-Kaehler or almost-Kaehler cases) and the target manifold $N$ to be 
\textbf{almost} $\beta _{N}$-Kenmotsu (which covers the Kenmotsu cases as
well as the $\beta _{N}$-Kenmotsu or $CR$-integrable almost $\beta _{N}$%
-Kenmotsu cases).

\textbf{Theorem }$\left( 4.10\right) :$

$1^{\circ }\mathbf{\ })$\textit{\ \ }\ \textit{Let }$f$\textit{\ }$:M$%
\textit{\ }$\rightarrow N$\textit{\ \ be a }$\left( \pm \right) $-\textit{%
holomorphic \textbf{horizontally} \textbf{weakly} conformal map with
dilation }$\mu $\textit{.}

$\ \ \ i)$ \ \textit{If \ }$M\in $\textit{\ }$\mathcal{A}_{1}\cup \mathcal{A}%
_{6}\cup \mathcal{A}_{7}$\textit{\ and \ }$N\in \mathcal{A}_{1}\cup \mathcal{%
A}_{4}\cup \mathcal{A}_{6}\cup \mathcal{A}_{7}$\textit{\ \ then \ }$f$%
\textit{\ \ is a horizontally homothetic harmonic morphism with minimal
fibres}

$\ \ ii)$\textit{\ \ If \ }$M\in \mathcal{A}_{4}\cup \mathcal{A}_{7}^{a}$ \
and $N\in \mathcal{A}_{1}\cup \mathcal{A}_{4}\cup \mathcal{A}_{6}\cup 
\mathcal{A}_{7}^{a}$ \ then \textit{\ \ }$f$\textit{\ \ is a harmonic
morphism if and only if \ }$d\lambda \left( \xi \right) =0.$

\textit{In particular,}

$\ a^{\circ })$ \textit{If } $f$ \ \textit{is horizontally homothetic ( In
particular, horizontally conformal Riemannian submersion) then if it is a
harmonic morphism with minimal fibres. Conversally, If } $f$ \ \textit{is a
harmonic morphism with minimal fibres then it is horizontally homothetic. }

$\ b^{\circ })$ \ \textit{If \ }$M\in A_{4}$\textit{\ and \ }$N\in A_{1}\cup
A_{4}\cup A_{6}\cup A_{7}$\textit{\ \ then }$\ f$\textit{\ \ is horizontally 
\textbf{homothetic} harmonic morphism (and therefore it has\ minimal fibres)
if and only if \ }$d\lambda \left( \xi \right) =0.$

$iii\mathbf{\ })$\textit{\ \ If \ }$M\in \mathcal{B}_{2}^{a}$ \ and $N\in 
\mathcal{B}_{2}^{a}$\ (that is, $M$ \ and $N$ \ are both \textbf{almost}
semi-cosymplectic manifolds) then \textit{\ \ }$f$\textit{\ \ is a harmonic
morphism if and only if \ }$d\lambda \left( \xi \right) =0.$

\textit{In particular, }

$\ a^{\circ })$ \textit{If } $f$ \ \textit{is horizontally homothetic ( In
particular, horizontally conformal Riemannian submersion) then if it is a
harmonic morphism with minimal fibres. Conversally, If } $f$ \ \textit{is a
harmonic morphism with minimal fibres then it is horizontally homothetic. \ }

$\ b^{\circ })$ \ \textit{If \ }$M\in \mathcal{B}_{1}\cup \mathcal{B}_{2}$%
\textit{\ and }$\ N\in \mathcal{B}_{1}\cup \mathcal{B}_{2}$\textit{\ \ then }%
$\ f$\textit{\ \ is horizontally \textbf{homothetic} harmonic morphism (and
therefore it has\ minimal fibres) if and only if \ }$d\lambda \left( \xi
\right) =0.$

$iv)$\textit{\ \ If \ }$M\in \mathcal{B}_{2}^{a}$\textit{\ \ and }$N\in 
\mathcal{B}_{3}$\textit{\ }(that is, $M$ is an \textbf{almost}
semi-cosymplectic manifold and $N$ is an \textbf{almost }$\beta _{N}$%
-Kenmotsu manifold)\textit{\ then \ \ }$f$\textit{\ \ is a harmonic morphism
if and only if \ }%
\begin{equation*}
d\lambda \left( \xi \right) -2n\lambda ^{2}\beta _{N}=0.
\end{equation*}

$\ v)$\textit{\ \ If \ }$M\in \mathcal{B}_{3}$\textit{\ \ and }$N\in $%
\textit{\ }$\mathcal{B}_{2}^{a}$\textit{\ }(that is, $M$ is an \textbf{%
almost }$\beta _{N}$-Kenmotsu manifold and $N$ is an \textbf{almost}
semi-cosymplectic manifold ) \textit{then \ \ }$f$\textit{\ \ is a harmonic
morphism if and only if \ }%
\begin{equation*}
d\lambda \left( \xi \right) +2m\lambda \beta _{M}=0.
\end{equation*}

$vi\ )$\textit{\ \ Let \ }$M,\ N\in \mathcal{B}_{3},$\textit{\ \ that is, }$%
M $\textit{\ \ and }$N$\textit{\ are almost }$\beta _{M}$-\textit{Kenmotsu
and} \textit{almost} $\beta _{N}$-\textit{Kenmotsu manifolds of dimensions }$%
m$\textit{\ and }$n$\textit{\ respectively. Then}

$\ \ a^{\circ })$ \ $f$\textit{\ \ is a harmonic morphism if \ and only if}%
\begin{equation*}
d\lambda _{p}\left( \xi \right) +2\lambda \left[ m\beta _{M}\left( p\right)
-n\lambda _{p}\beta _{N}\left( q\right) \right] =0,\text{ \ \ }\forall \text{%
\ }p\in M\text{ \ \textit{and} \ }q=f\left( p\right) .
\end{equation*}%
\textit{Suppose further that }$f$\textit{\ \ has minimal fibres then }$f$%
\textit{\ \ is harmonic (so is a harmonic morphism) if and only if }$f$%
\textit{\ \ is horizontally homothetic whose dilation is given by \ \ }%
\begin{equation*}
\mu \left( p\right) =\left\vert \lambda \left( p\right) \right\vert
=\left\vert \frac{m\beta _{M}\left( p\right) }{n\beta _{N}\left( q\right) }%
\right\vert ,\ \ \ \ p\in M
\end{equation*}%
\textit{(which is constant along horizontal curves).}

$\ \ b^{\circ })$ \ \textit{The two of the followig imply the third:}

$\ \ \ \ \ \bullet $\textit{\ \ \ }$f$\textit{\ \ is a harmonic morphism}

$\ \ \ \ \ \bullet $\textit{\ \ \ }$\lambda \left( p\right) =\frac{m\beta
_{M}\left( p\right) }{n\beta _{N}\left( q\right) },$\textit{\ \ \ }$\forall $%
\textit{\ }$p\in M$\textit{\ \ and \ }$q=f\left( p\right) $

$\ \ \ \ \ \bullet $\textit{\ \ \ }$d\lambda _{p}\left( \xi \right) =0,$%
\textit{\ \ \ \ \ }$\forall $\textit{\ }$p\in M$

\textit{\ }$\ c^{\circ })$ \ \textit{Suppose }$\beta _{M}$\textit{\ \ and \ }%
$\beta _{N}$\textit{\ \ are both \textbf{constant} funtions. Then followig
are equivalent:}

$\ \ \ \ \ \bullet $\textit{\ \ \ }$f$\textit{\ \ is a harmonic morphism
with minimal fibres}

$\ \ \ \ \ \bullet $\textit{\ \ \ }$f$\textit{\ \ is a horizontally
homothetic with dilation}%
\begin{equation*}
\mu =\mathit{\ }\left\vert \lambda \right\vert =\left\vert \frac{m\beta _{M}%
}{n\beta _{N}}\right\vert
\end{equation*}

$2^{\circ }\mathbf{\ })$\textit{\ }\ \textit{Let }$f$\textit{\ }$:M$\textit{%
\ }$\rightarrow \left( H_{1}^{2n},\ J_{1},\ G_{1}\right) $\textit{\ \ be a }$%
\left( \pm \right) $\textit{-holomorphic horizontally weakly conformal map
in to an almost Hermitian manifold,\ where \ }$M\in \mathcal{A}_{1}\cup 
\mathcal{A}_{2}\cup \mathcal{A}_{3}\cup \mathcal{A}_{4}\cup \mathcal{A}%
_{6}\cup \mathcal{A}_{7}^{a}\cup \mathcal{B}_{2}^{a}\cup \mathcal{B}_{3}\cup 
\mathcal{C}_{2}.$ \textit{\ Then \ }$f$\textit{\ \ is a harmonic morphism if
and only if \ }$H_{1}^{2n}$\textit{\ is semi-Kaehler.}

\textbf{Proof:}

$1^{\circ }\mathbf{\ })$ Combining Proposition $\left( 3.1\right) /1^{\circ
} $ \ and Lemma $\left( 3.3\right) /2^{\circ }$ \ we have%
\begin{equation*}
\mu =\left\vert \lambda \right\vert
\end{equation*}%
and that 
\begin{eqnarray}
\mathcal{T}\left( f\right) &=&U_{\xi }+\overset{m}{\underset{i=1}{\sum }}%
\left\{ S^{N}\left( \phi E_{i},\ E_{i}\right) -\ f_{\ast }S^{M}\left(
\varphi e_{i},e_{i}\right) \right\}  \notag \\
&=&d\lambda \left( \xi \right) \gamma -\ f_{\ast }\left\{ \left( \nabla
_{\xi }^{M}\xi \right) +\underset{i=1}{\overset{m}{\sum }}S^{M}\left(
\varphi e_{i},e_{i}\right) \right\}  \TCItag{$4.19$} \\
&&+\lambda ^{2}\left\{ \left( \nabla _{\gamma }^{N}\gamma \right) +\underset{%
i=1}{\overset{n}{\sum }}S^{N}\left( \phi v_{i},\ v_{i}\right) \right\} 
\notag
\end{eqnarray}%
for some orthonormal frame field $\left\{ e_{1},\cdots ,e_{m};\varphi
e_{1},\cdots ,\varphi e_{m}\right\} $ \ for $D^{M}.$ \ Here 
\begin{equation*}
\ f_{\ast }e_{i}=\left\{ 
\begin{array}{c}
\mu v_{i}=\left\vert \lambda \right\vert v_{i},\ \ \ i=1,\cdots ,n\ \ \ \ 
\\ 
0,\ \ \ i>n\ \ \ \ \ \ \ \ \ \ 
\end{array}%
\right\}
\end{equation*}%
and \ $\left\{ v_{1},\cdots ,v_{n};\phi v_{1},\cdots ,\phi v_{n}\right\} $ \
is a local orthonormal frame field for $D^{N},$ \ (see Lemma $\left(
3.3\right) /4^{\circ }\left( i\right) ).$

$\left( i\right) :$ \ Since \ $M$ and \ $N$ both satisfy $\left( GC\right) $
\ and \ $M$ \ is non-semi-$\varphi $-involutive, Proposition $\left( \left(
3.3\right) /2^{\circ }\right) $ gives that\ \ $\left\vert \lambda
\right\vert =\mu $ \ is constant and that \ $\mathcal{U}_{\xi }=0.$ \ So $\
f $\ \ is horizontally homothetic and by the virtue of \ $\left( 4.19\right)
,$ \ it satisfies

\begin{equation*}
\mathcal{T}\left( f\right) =\lambda ^{2}\underset{i=1}{\overset{n}{\sum }}%
S^{N}\left( \phi v_{i},\ v_{i}\right) -f_{\ast }\underset{i=1}{\overset{m}{%
\sum }}S^{M}\left( \varphi e_{i},e_{i}\right) .
\end{equation*}%
But then, $\mathcal{T}\left( f\right) =0$ \ since,\ from the \textbf{TABLE-II%
}, one sees that

\begin{equation*}
\underset{i=1}{\overset{n}{\sum }}S^{N}\left( \phi v_{i},\ v_{i}\right) =0%
\text{ \ \ \ \ }and\ \ \ \underset{i=1}{\overset{m}{\sum }}S^{M}\left(
\varphi e_{i},e_{i}\right) =0.\ \ \ \ 
\end{equation*}%
That is, \ $f$ \ is horizontally homothetic harmonic morphism. Minimality of
the fibres follows from Lemma $\left( 3.5\right) .$

$\left( ii\right) :$ \ Observe that (see \textbf{TABLE-II)}

\begin{equation*}
\left( \nabla _{\xi }^{M}\xi \right) +\underset{i=1}{\overset{m}{\sum }}%
S^{M}\left( \varphi e_{i},e_{i}\right) =0
\end{equation*}%
and

\begin{equation*}
\left( \nabla _{\gamma }^{N}\gamma \right) +\underset{i=1}{\overset{n}{\sum }%
}S^{N}\left( \phi v_{i},\ v_{i}\right) =0.
\end{equation*}%
Then $\left( 4.19\right) $ gives that 
\begin{equation*}
\mathcal{T}\left( f\right) =d\lambda \left( \xi \right) \gamma ,
\end{equation*}%
so the required result follows.

In particular,

$\ a^{\circ })$ Suppose $f$ \ is horizontally homothetic. Then $d\lambda
\left( \xi \right) =0$, so by $(ii)$ \ and Lemma $(3.5),$\textit{\ \ }$f$%
\textit{\ \ }is a harmonic morphism with minimal fibres. The converse is
just Lemma $(3.5).$

$b^{\circ })$ Suppose \ $f$\ \ is horizontally homothetic harmonic morphism.
Then clearly $d\lambda \left( \xi \right) =0.$ (Minimality of fibres is due
to Lemma $(3.5)).$ Conversaly suppose $\ d\lambda \left( \xi \right) =0.$
So, $f$\textit{\ \ }is a harmonic morphism. On the other hand, since \ $M\in 
\mathcal{A}_{4}$ and \ $N\in \mathcal{A}_{1}\cup \mathcal{A}_{4}\cup 
\mathcal{A}_{6}\cup \mathcal{A}_{7}$ we see that \ $M$ $\ $and$\ N$ \ both
satisfy $\left( GC\right) $ \ and therefore, Proposition $\left( \left(
3.3\right) /1^{\circ }\right) ,$ \ $\left\vert \lambda \right\vert =\mu $ \
is constant along \ $D^{M}.$ This together with \ $d\lambda \left( \xi
\right) =0$ \ give us that \ $\lambda $\ \ is constant. That is, $f$ \ is
horizontally homothetic harmonic morphism.

\bigskip $\left( iii\right) :$ \ By the same argument used in $\left(
1^{\circ }/\left( ii\right) \right) $ we get the result.

$\left( iv\right) :$ \ From the \textbf{TABLE-II}, one sees that%
\begin{equation*}
\left( \nabla _{\xi }^{M}\xi \right) +\underset{i=1}{\overset{m}{\sum }}%
S^{M}\left( \varphi e_{i},e_{i}\right) =0
\end{equation*}%
since $M\in \mathcal{B}_{2}^{a}.$ \ Also%
\begin{equation*}
\left( \nabla _{\gamma }^{N}\gamma \right) =0\text{ \ \ \ \ and \ \ \ }%
\underset{i=1}{\overset{n}{\sum }}S^{N}\left( \phi v_{i},\ v_{i}\right)
=-2n\beta _{N}\gamma .
\end{equation*}%
since $N\in \mathcal{B}_{3}.$ \ Therefore, $\left( 4.19\right) $ \ gives us
that%
\begin{equation*}
\mathcal{T}\left( f\right) =\left( d\lambda \left( \xi \right) -2n\lambda
^{2}\beta _{N}\right) \gamma .
\end{equation*}%
So the result follows.

\bigskip $\left( v\right) :$ \ \ From the \textbf{TABLE-II}, one sees that%
\begin{equation*}
\left( \nabla _{\xi }^{M}\xi \right) =0\text{ \ \ \ \ and \ \ \ }\underset{%
i=1}{\overset{m}{\sum }}S^{M}\left( \varphi e_{i},e_{i}\right) =-2m\beta
_{M}\xi .
\end{equation*}%
since $M\in \mathcal{B}_{3}.$ \ Also%
\begin{equation*}
\nabla _{\gamma }^{N}\gamma +\underset{i=1}{\overset{n}{\sum }}S^{N}\left(
\phi v_{i},\ v_{i}\right) =0.
\end{equation*}%
since $N\in \mathcal{B}_{2}^{a}.$ \ Therefore, $\left( 4.19\right) $ \ gives
us that%
\begin{equation*}
\mathcal{T}\left( f\right) =\left( d\lambda \left( \xi \right) +2m\lambda
\beta _{M}\right) \gamma .
\end{equation*}%
So the result follows.

$\left( vi\right) :$ \ \ \ Since \ $M$ and \ $N$ both satisfy $\left(
GC\right) ,$\ from the \textbf{TABLE-II}, we see that$\ \ \left( 4.19\right) 
$ \ gives

\begin{eqnarray}
\mathcal{T}\left( f\right) &=&d\lambda \left( \xi \right) \gamma +\ f_{\ast
}\left( 2m\beta _{M}\xi \right) -2n\lambda ^{2}\beta _{N}\gamma  \notag \\
&=&\left\{ d\lambda \left( \xi \right) +2\lambda \left( m\beta _{M}-n\lambda
\beta _{N}\right) \right\} \gamma  \TCItag{$4.20$}
\end{eqnarray}

$\ \ \ \left( vi\right) /\left( a^{\circ }\right) :$ The first part of this
is immediate from $\left( 4.20\right) .$ For the second part, suppose
further that $f$ \ has minimal fibres. Now if $f$\ \ is harmonic (so is a
harmonic morphism) then, Lemma $\left( 3.5\right) $ \ gives that \ $f$ \ is
horizontally homothetic and so that the dilation $\mu =\left\vert \lambda
\right\vert $ \ is constant along the horizotal curves. Therefore we have $%
d\lambda \left( \xi \right) =0$ \ since $\xi $ \ is horizontal by Lemma $%
\left( 3.3\right) /\left( 2^{\circ }\text{-}i\right) .$ \ But then $\left(
4.20\right) $ together with the harmonicity gives 
\begin{equation*}
m\beta _{M}\left( p\right) -2n\lambda \beta _{N}\left( q\right) =0.
\end{equation*}%
Conversaly, if $f$ \ is horizontally homothetic with dilation \textit{\ \ }%
\begin{equation*}
\mu \left( p\right) =\left\vert \lambda \right\vert \left( p\right)
=\left\vert \frac{m\beta _{M}\left( p\right) }{n\beta _{N}\left( q\right) }%
\right\vert ,\ \ \ \ p\in M
\end{equation*}%
then we get that $\ d\lambda \left( \xi \right) =0$ since $\xi $ \ is
horizontal\ and \ $2\lambda \left( m\beta _{M}\left( p\right) -2n\lambda
\beta _{N}\left( q\right) \right) =0.$ \ So the harmonicity follows from $%
\left( 4.20\right) .$

$\ \ \ \left( vi\right) /\left( b^{\circ }\right) :$ It is immediate from $%
\left( 4.20\right) .$

$\ \ \ \left( vi\right) /\left( c^{\circ }\right) :$ $\ $\ If\textit{\ \ }$f$%
\textit{\ \ }is a harmonic morphism with minimal fibres, then by Lemma $%
\left( 3.5\right) ,$ it\textit{\ }is horizontally homothetic. So, $\mu
=\left\vert \lambda \right\vert $ is horizontally constant and therefore $%
d\lambda \left( \xi \right) =0$ $\ $since $\xi $ \ is horizontal. \ But then
harmonicity of \ $f$ \ and $\left( 4.20\right) $ \ gives that

\begin{equation*}
m\beta _{M}-2n\lambda \beta _{N}=0
\end{equation*}%
so that $\mu =\lambda =\left\vert \frac{m\beta _{M}}{n\beta _{N}}\right\vert
.$ \ Conversaly, If\textit{\ \ }$f$\textit{\ \ }is horizontally homothetic 
\textit{with dilation }$\mu =\lambda =\left\vert \frac{m\beta _{M}}{n\beta
_{N}}\right\vert $ then $\ d\lambda \left( \xi \right) =0.$ So, from $\left(
4.20\right) $ \ we get%
\begin{equation*}
\mathcal{T}\left( f\right) =2\lambda \left( m\beta _{M}-2n\lambda \beta
_{N}\right) \gamma .
\end{equation*}%
But then this gives that $\mathcal{T}\left( f\right) =0,$ (that is, \ $f$ \
is a harmonic morphism)$.$ Also we see from Lemma $\left( 3.5\right) $ that,
every horizontally homothetic harmonic map has minimal fibres. This
completes the proof.

$2^{\circ }\mathbf{\ })$\textit{\ }\ Combining Proposition $\left(
3.1\right) /\left( 2^{\circ }\right) $ \ and Lemma $\left( 3.3\right)
/\left( 3^{\circ }\right) $ \ we have

If $M\in \mathcal{A}_{1}\cup \mathcal{A}_{2}\cup \mathcal{A}_{3}\cup 
\mathcal{A}_{4}\cup \mathcal{A}_{6}\cup \mathcal{A}_{7}^{a}\cup \mathcal{B}%
_{2}^{a}\cup \mathcal{B}_{3}$ \ then 
\begin{equation}
\mathcal{T}\left( f\right) =\mu ^{2}\overset{n}{\underset{i=1}{\sum }}%
S^{H_{1}}\left( J_{1}v_{i},\ v_{i}\right) -\ f_{\ast }\left\{ \left( \nabla
_{\xi }^{M}\xi \right) +\overset{m}{\underset{i=1}{\sum }}S^{M}\left(
\varphi e_{i},e_{i}\right) \right\}  \tag{$4.21$}
\end{equation}%
for some orthonormal frame field $\left\{ e_{1},\cdots ,e_{m};\varphi
e_{1},\cdots ,\varphi e_{m}\right\} $ \ for $D^{M}.$

If $M\in \mathcal{C}_{2}$ then 
\begin{equation}
\mathcal{T}\left( f\right) =\mu ^{2}\overset{n}{\underset{i=1}{\sum }}%
S^{H_{1}}\left( J_{1}v_{i},\ v_{i}\right) -\ f_{\ast }\overset{m}{\underset{%
i=1}{\sum }}S^{H}\left( Je_{i},e_{i}\right) .  \tag{$4.22$}
\end{equation}%
for some orthonormal frame field $\left\{ e_{1},\cdots ,e_{m};Je_{1},\cdots
,Je_{m}\right\} $ \ for \ $TH^{2m}$

\ Here 
\begin{equation*}
\ f_{\ast }e_{i}=\left\{ 
\begin{array}{c}
\mu v_{i}=\left\vert \lambda \right\vert v_{i},\ \ \ i=1,\cdots ,n\ \ \ \ 
\\ 
0,\ \ \ i>n\ \ \ \ \ \ \ \ \ \ 
\end{array}%
\right\}
\end{equation*}%
and \ $\left\{ v_{1},\cdots ,v_{n};J_{1}v_{1},\cdots ,J_{1}v_{n}\right\} $ \
is a local orthonormal frame field for $TH_{1},$ \ (see Lemma $\left(
3.3\right) /4^{\circ }\left( i,ii\right) ).$

Observe now that

$\bullet $ \ when $M\in \mathcal{A}_{1}\cup \mathcal{A}_{2}\cup \mathcal{A}%
_{4}\cup \mathcal{A}_{6}$ \ we have $\ $%
\begin{equation*}
\nabla _{\xi }^{M}\xi =0\text{ \ and \ }\overset{m}{\underset{i=1}{\sum }}%
S^{M}\left( \varphi e_{i},e_{i}\right) =0
\end{equation*}%
since $M$ \ is semi-symplectic and satisfies $\left( GC\right) .$

$\bullet $ \ when $M\in \mathcal{A}_{7}^{a}\cup \mathcal{B}_{2}^{a}$, ($ie.$
it is almost semi-Sasakian or almost semi-cosymplectic) we have%
\begin{equation*}
\nabla _{\xi }^{M}\xi +\overset{m}{\underset{i=1}{\sum }}S^{M}\left( \varphi
e_{i},e_{i}\right) =0.
\end{equation*}

$\bullet $ \ when $M\in \mathcal{B}_{3}\cup \mathcal{A}_{3}$, we have $\ \
\nabla _{\xi }^{M}\xi =0$ \ \ and \ 
\begin{equation*}
\overset{m}{\underset{i=1}{\sum }}S^{M}\left( \varphi e_{i},e_{i}\right)
=\left\{ 
\begin{array}{ll}
-2m\beta _{M}\xi , & M\in \mathcal{B}_{3} \\ 
\left( \delta \eta \right) \xi , & M\in \mathcal{A}_{3}%
\end{array}%
\right\} \ \ \ 
\end{equation*}%
So, since $\ f_{\ast }\xi =0,$ $\ $%
\begin{eqnarray*}
f_{\ast }\left\{ \left( \nabla _{\xi }^{M}\xi \right) +\overset{m}{\underset{%
i=1}{\sum }}S^{M}\left( \varphi e_{i},e_{i}\right) \right\} &=&\left\{ 
\begin{array}{ll}
-2m\beta _{M}f_{\ast }\xi , & M\in \mathcal{B}_{3} \\ 
\left( \delta \eta \right) f_{\ast }\xi , & M\in \mathcal{A}_{3}%
\end{array}%
\right\} \\
&=&0
\end{eqnarray*}%
Finally,

$\bullet $ \ when $M=\left( H^{2n},\ J,\ G\right) \in \mathcal{C}_{2}$, we
have

\begin{equation*}
\overset{m}{\underset{i=1}{\sum }}S^{H}\left( Jv_{i},\ v_{i}\right) =0
\end{equation*}%
So, in all the above cases the equation $\left( 4.21\right) $ and $\left(
4.22\right) $ reduce to%
\begin{equation*}
\mathcal{T}\left( f\right) =\mu ^{2}\overset{n}{\underset{i=1}{\sum }}%
S^{H_{1}}\left( J_{1}v_{i},\ v_{i}\right)
\end{equation*}%
Hence the result follows.

\textbf{Remark }$\left( 4.10\right) :$ In $\left[ 6\right] $, Theorem $%
\left( 4.1\right) $ states that

$\ \bullet $ \textit{Any horizontally conformal }$\left( \varphi ,\phi
\right) $\textit{-holomorphic (Riemannian) submersion between quasi-}$%
\mathcal{K}$\textit{-cosymlectic manifolds is harmonic with minimal fibres
and so in particular a harmonic morphism. \ }

Our Theorem $\left( 4.10\right) /\left( 1^{\circ }\text{-}ii\left( a^{\circ
}\right) \right) $ generalizes this result.

\end{document}